\documentclass{amsart}
\usepackage{amsthm, amscd}
\usepackage{amsmath}
\usepackage{amssymb}
\usepackage{amsfonts}
\usepackage{mathrsfs}
\usepackage{esint}
\usepackage{pdfsync}
\usepackage[textsize=tiny]{todonotes}
\usepackage{tikz, float}
\usepackage{pgf}
\usepackage{pgfplots}
\usetikzlibrary{angles,quotes}
\usepackage{graphicx}
\newtheorem{thm}{Theorem}[section]
\newtheorem{Prop}[thm]{Proposition}
\newtheorem{Def}[thm]{Definition}
\newtheorem{Lem}[thm]{Lemma}
\newtheorem{Cor}[thm]{Corollary}
\newtheorem{rk}[thm]{Remark}

\newcommand{\nn}{\mathbb{N}}
\newcommand{\R}{\mathbb{R}}
\newcommand{\rr}{\mathbb{R}}
\newcommand{\cc}{\mathbb{C}}
\newcommand{\CC}{\mathbb{C}}
\newcommand{\eps}{\varepsilon}

\newcommand{\real}{\mathbb{R}}
\def\un{{\mathrm{1~\hspace{-1.4ex}l}}}

\def\un{{\mathrm{1~\hspace{-1.4ex}l}}}
\numberwithin{equation}{section}

\author{Karine Beauchard, Michela Egidi \& Karel Pravda-Starov}

\title[Conditions for null-controllability  with moving control]{Geometric conditions for the null-controllability of hypoelliptic quadratic parabolic equations with moving control supports}

\address{\noindent \textsc{Karine Beauchard, Univ Rennes, CNRS, IRMAR - UMR 6625, F-35000 Rennes, France
}}
\email{karine.beauchard@ens-rennes.fr}

\address{\noindent \textsc{Michela Egidi, 
Ruhr Universit\"at Bochum, 
Fakult\"at f\"ur Mathematik, 
Gebaude IB 3/55, Universit\"atsstr. 150, 44780 Bochum, Germany
}}
\email{michela.egidi@ruhr-uni-bochum.de}

\address{\noindent \textsc{Karel Pravda-Starov, Univ Rennes, CNRS, IRMAR - UMR 6625, F-35000 Rennes, France
}}
\email{karel.pravda-starov@univ-rennes1.fr}


\thanks{The first and third authors express their gratefulness to the Centre de Math\'ematiques Henri Lebesgue for the very stimulating scientific environment.
The second author acknowlegdes the support of the DFG Project Ve 253/7-1 ``Multiscale Version of the Logvinenko-Sereda Theorem''.}

\begin{document}

\begin{abstract}
We study the null-controllability of some hypoelliptic quadratic parabolic equations posed on the whole Euclidean space with moving control supports, and provide necessary or sufficient geometric conditions on the moving control supports to ensure null-controllability. The first class of equations is the one associated to non-autonomous Ornstein-Uhlenbeck operators satisfying a generalized Kalman rank condition. 
In particular, when the moving control supports comply with the flow associated to the transport part of the Ornstein-Uhlenbeck operators, a necessary and sufficient condition for null-controllability on the moving control supports  is established. The second class of equations is the class of accretive non-selfadjoint quadratic operators with zero singular spaces for which some sufficient geometric conditions on the moving control supports are also given to ensure null-controllability.   
\end{abstract}

\maketitle


\section{Introduction}

We study parabolic equations posed on the whole Euclidean space $\rr^d$,
\begin{equation}\label{syst_general1}
\left\lbrace \begin{array}{ll}
(\partial_t + P)f(t,x)=\un_{\omega(t)}(x)u(t,x), \quad &  x \in \mathbb{R}^d, t>0, \\
f|_{t=0}=f_0 \in L^2(\rr^d),                                       &  
\end{array}\right.
\end{equation}
and controlled by a source term $u$ locally distributed in a time-dependent control subset $\omega(t) \subset \mathbb{R}^d$. The controllability of partial differential equations with moving control subsets is at the core of current investigations and the topics of several recent works \cite{Chaves-Rosier-Zuazua, LR-L-T-T,Rosier_2013}.

We consider in this work two specific classes of hypoelliptic quadratic parabolic equations, and we aim at pointing out necessary or sufficient geometric conditions on the moving control subsets $(\omega(t))_{t \in I}$ to ensure null-controllability. 
In order to ensure the well-posedness of the evolution equations (\ref{syst_general1}), the moving control subsets are assumed to satisfy the following measurability property:

\medskip

\begin{Def}[Moving control support] \label{Def:MCS}
Let $\Omega$ be an open subset of $\mathbb{R}^d$ and $I$ be an interval of  $\rr$.
A moving control support on $I$ in $\Omega$ is a family $(\omega(t))_{t \in I}$ of subsets of $\mathbb{R}^d$ such that the map
$(t,x) \in I \times \Omega \mapsto \un_{\omega(t)}(x)$ is measurable,
where $\un_{\omega(t)}$ denotes the characteristic function of the moving set $\omega(t)$.
\end{Def}

\medskip

The Cauchy problems (\ref{syst_general1}) studied in this work will all be well-posed in the space $C^0([0,T],L^2(\mathbb{R}^d))$, for any initial datum $f_0 \in L^2(\mathbb{R}^d)$ and control function $u \in L^2((0,T)\times\mathbb{R}^d)$. The previous definition of moving control support $(\omega(t))_{t \in I}$ does not rule out the case when the control subsets are empty $\omega(t)=\emptyset$ at certain times. Let us notice in particular that if $(\omega(t))_{t \in I}$ is a moving control support in $\rr^d$, and $E$ is a measurable subset of $I$, then the family of subsets $(\widetilde{\omega}(t))_{t\in I}$ defined as $\widetilde{\omega}(t)=\omega(t)$ if $t\in E$, or  $\widetilde{\omega}(t)=\emptyset$ if $t \in I \setminus E$, also defines a moving control support in $\mathbb{R}^d$, since $\un_{\widetilde{\omega}(t)}(x)=\un_{E}(t) \un_{\omega(t)}(x)$ is a product of measurable functions.
 
 \medskip
 
\begin{Def}[Null-controllability and control cost]
Let $T>0$ and $(\omega(t))_{t \in [0,T]}$ be a moving control support on the time interval $[0,T]$ in $\mathbb{R}^d$.
Equation \eqref{syst_general1} is said to be null-controllable on the time interval $[0,T]$ if,
for any initial datum $f_0 \in L^{2}(\mathbb{R}^d)$,
there exists a control function $u \in L^2((0,T)\times\mathbb{R}^d)$,
such that the solution of \eqref{syst_general1} satisfies $f(T,\cdot)=0$. If the equation is null-controllable, the control cost is defined as the smallest positive constant $C_T>0$ such that any initial datum $f_0 \in L^{2}(\mathbb{R}^d)$
can be steered to zero by means of a control function $u \in L^2((0,T)\times\mathbb{R}^d)$ satisfying
$$\| u \|_{L^2((0,T)\times\rr^d)} \leq C_T \| f_0\|_{L^2(\rr^d)}.$$
\end{Def}

\medskip

We consider in this work two specific classes of hypoelliptic quadratic parabolic equations (\ref{syst_general1}). The first one is the class of evolution equations associated to non-autonomous Ornstein-Uhlenbeck operators satisfying a generalized Kalman rank condition described in Section~\ref{OUclass}. The null-controllability with fixed control subsets $\omega(t)=\omega_0$ of these non-autonomous equations was studied by the first and third authors in~\cite{KB_KPS_JMAA}, and the following sufficient condition for null-controllability on fixed open control subsets was established in \cite{KB_KPS_JMAA} (Theorem 1.3),
\begin{equation}\label{lop20}
\exists r,\delta>0, \forall y \in \mathbb{R}^d, \exists y' \in \omega_0, \quad B_d(y',r)\subset\omega_0 \text{ and } |y-y'|<\delta,
\end{equation}
where $B_d(y',r)$ denotes the open Euclidean ball centered at $y'$ with radius $r$. The second class studied in this article is the class of evolution equations associated to accretive non-selfadjoint quadratic operators with zero singular spaces described in Section~\ref{quadclass}. The null-controllability with fixed control subsets $\omega(t)=\omega_0$ of these hypoelliptic equations was studied by the first and third authors in~\cite{KB_KPS_JEP} (Theorem~1.4), and was shown to hold for any fixed control subset satisfying the very same geometric condition (\ref{lop20}). The results in~\cite{KB_KPS_JMAA,KB_KPS_JEP} were some first steps outlining and providing preliminary insights on the geometry of the control subsets needed to get null-controllability for these two classes of hypoelliptic non-selfadjoint evolution equations. However, the geometric condition (\ref{lop20}) was not expected to be sharp to ensure null-controllability, and a new breakthrough was then made by Veseli\'{c} and the second author in~\cite{Egidi_Veselic}, who established that the following notion of thickness is a necessary and sufficient condition on fixed control subsets to ensure the null-controllability of the heat equation posed on the whole Euclidean space $\rr^d$ in some positive time, as well as in any positive time: 

\medskip

\begin{Def}[$(\delta,\alpha)$-thick set]\label{Def:thick}
Let $0<\delta \leq 1$ and $\alpha=(\alpha_1,...,\alpha_d)\in(0,+\infty)^{d}$.
A measurable subset $S \subset \mathbb{R}^d$ is a $(\delta,\alpha)$-thick set in $\mathbb{R}^d$ if the following estimate holds:
$$\forall x \in \mathbb{R}^d, \quad \lambda\big( S \cap(x+[0,\alpha_1]\times...\times[0,\alpha_d])\big) \geq \delta \prod_{j=1}^d\alpha_j,$$
with $\lambda$ being the Lebesgue measure on $\rr^d$.
A measurable subset $S \subset \mathbb{R}^d$ is said to be thick in $\rr^d$
if $S$ is a  $(\delta,\alpha)$-thick subset in $\mathbb{R}^d$,
for some $0<\delta \leq 1$ and $\alpha\in(0,+\infty)^{d}$.
\end{Def}

\medskip

The very same result about the heat equation was obtained independently by Wang, Wang, Zhang and Zhang in~\cite{Gengsheng_Wang2}. As the heat equation is a particular example of Ornstein-Uhlenbeck equations which obviously enjoys the very specific features of being autonomous, elliptic, as well as having only a diffusive structure with no transport part, it is then natural to wonder to which extent this thickness condition is also relevant to ensure the null-controllability of general Ornstein-Uhlenbeck equations, and whether the null-controllability results obtained   in~\cite{KB_KPS_JMAA,KB_KPS_JEP} actually extend for thick control subsets. We therefore aim in this work at sharply understanding which geometry on the control subsets rules the null-controllability of Ornstein-Ulhenbeck equations, and how possible transport phenomena induced by these equations interplay with the geometry of the control subsets. This geometry of the control subsets naturally ends up to be time-dependent because of the transport phenomena associated to general Ornstein-Ulhenbeck equations. The results given in Section~\ref{OUclass} provide some necessary or sufficient geometric conditions on the moving control subsets to ensure null-controllability. In the particular case when the moving control subsets comply with the flow associated to the transport part of the Ornstein-Uhlenbeck operators, a necessary and sufficient condition related to the thickness property of the moving control subsets is derived. In a second part of this work, we consider evolution equations associated to accretive non-selfadjoint quadratic operators with zero singular spaces. The results given in Section~\ref{quadclass} show how those established in~\cite{KB_KPS_JEP} can actually extend in the framework of moving control subsets. For this second class of hypoelliptic quadratic parabolic equations, we unfortunately provide only a sufficient condition for null-controllability. It would be of course most interesting to also derive a necessary condition for null-controllability as for Ornstein-Ulhenbeck equations. However, the transport phenomena at play for this second class of equations are far more complex and this topic is not studied in this work.

\subsection{Null-controllability of non-autonomous Ornstein-Uhlenbeck equations}\label{OUclass}

\subsubsection{Non-autonomous Ornstein-Uhlenbeck operators and generalized Kalman rank condition}

We consider the evolution equation 
\begin{equation}\label{syst_general}
\left\lbrace \begin{array}{l}
\partial_t f -  \frac{1}{2}\textrm{Tr}(A(t)A(t)^T\nabla_x^2f) - \langle B(t)x, \nabla_x f\rangle
=\un_{\omega(t)}(x)u,  \\
f|_{t=0}=f_0 \in L^2(\rr^d),
\end{array}\right.
\end{equation}
associated to the non-autonomous Ornstein-Uhlenbeck operator
\begin{multline}\label{NAOU}
P(t)=\frac{1}{2}\textrm{Tr}\big(A(t)A(t)^T\nabla_x^2\big) + \big\langle B(t)x, \nabla_x \big\rangle\\
=\frac{1}{2}\sum_{i,j,k=1}^{d}a_{i,k}(t)a_{j,k}(t)\partial_{x_i,x_j}^2+\sum_{i,j=1}^d b_{i,j}(t)x_j\partial_{x_i}, \quad   x \in \mathbb{R}^d,
\end{multline}
where $$A=(a_{i,j})_{1 \leq i,j \leq d},\quad B=(b_{i,j})_{1 \leq i,j \leq d} \in C^{\infty}(I,M_{d}(\rr)),$$
are smooth mappings with values in real $d \times d$ matrices, with $I$ being an open interval of $\rr$ containing zero,
$A(t)^T$ standing for the transpose matrix of $A(t)$.
The well-posedness of the Cauchy problem (\ref{syst_general}) is proved in appendix (Section~\ref{appensix:well-pos}).

In order to ensure appropriate smoothing properties in Gevrey spaces, which are a key ingredient to establish null-controllability thanks to an adapted Lebeau-Robbiano method when using the abstract observability result in Theorem~\ref{Meta_thm_AdaptedLRmethod}, we assume that the following generalized Kalman rank condition holds:

 To that end, we define by induction the sequence of smooth mappings $(\tilde{A}_k)_{k \geq 0} \in C^{\infty}(I,M_{d}(\rr))^{\mathbb{N}}$ by
\begin{equation}\label{asd2}
\forall t \in I, \quad \tilde{A}_0(t)=A(t),
\end{equation}
\begin{equation}\label{asd3}
\forall k \geq 0, \forall t \in I, \quad \tilde{A}_{k+1}(t)= \frac{d}{dt}\tilde{A}_{k}(t)+B(t)\tilde{A}_{k}(t).
\end{equation}

\medskip

\begin{Def}[Generalized Kalman rank condition]
The generalized Kalman rank condition is said to hold at some time $T>0$, if $T\in I$ and
\begin{equation} \label{Kalman_time}
\emph{\text{Span}}\{\tilde{A}_k(T)x :\ x\in\mathbb{R}^d, k \geq 0\}=\mathbb{R}^d.
\end{equation}
\end{Def}

\medskip

The condition (\ref{Kalman_time}) was shown by Chang~\cite{chang}, and by Silverman and Meadows~\cite{silverman} to be sufficient for the controllability of the linear control system $\dot{x}=-B(t)x+A(t)u$ on the interval $I$. As noticed in~\cite[p. 11]{coron_book}, the two following vector spaces 
$$\textrm{Span}\{\tilde{A}_{k}(T)x : \ x \in \rr^d, k \geq 0\} \neq \textrm{Span}\{\tilde{A}_{k}(T)x : \ x \in \rr^d, 0 \leq k \leq d-1\},$$
are in general distinct, as contrary to the constant case studied in Section~\ref{constantOU}, the Cayley-Hamilton theorem does not apply. However, it was proved by Coron in~\cite{coron_book} (Proposition~1.19) that when the condition (\ref{Kalman_time}) holds at some time $T \in I$, then there exists a positive constant $\eps>0$ such that 
\begin{equation}\label{Kalman_time11}
\forall t \in I \cap (T-\eps,T+\eps)\setminus \{T\}, \quad \textrm{Span}\{\tilde{A}_{k}(t)x : \ x \in \rr^d, 0 \leq k \leq d-1\}=\rr^d.
\end{equation}
This assertion (\ref{Kalman_time11}) can be reformulated as
\begin{equation}\label{Kalman_time2}
\forall t \in I \cap (T-\eps,T+\eps)\setminus \{T\}, \quad \textrm{Rank}[\tilde{A}_{0}(t),\tilde{A}_{1}(t),\dots,\tilde{A}_{d-1}(t)]=d,
\end{equation}
where $[\tilde{A}_{0}(t),\tilde{A}_{1}(t),\dots,\tilde{A}_{d-1}(t)]$ is the $d\times d^2$ matrix obtained by writing consecutively the columns of the matrices $\tilde{A}_{j}(t)$. The above formula directly relate to the classical Kalman rank condition (\ref{kal1}) appearing in the autonomous case.

\subsubsection{Geometric conditions for null-controllability}
The main result of null controllability for non-autonomous Ornstein-Uhlenbeck equations (\ref{syst_general}) contained in this work is the following one:

\medskip

\begin{thm} \label{thm_OU_CS_CN}
Let $T>0$ and $(\omega(t))_{t \in [0,T]}$ be a moving control support in $\mathbb{R}^d$.
We assume that the generalized Kalman rank condition \emph{(\ref{Kalman_time})} holds at time $T$.
\begin{enumerate}
\item[$(i)$] \emph{\textbf{(Sufficient condition).}}
Let $\delta>0$, $\alpha \in (0,+\infty)^d$ and $E$ be a measurable subset of $[0,T]$ satisfying
\begin{equation} \label{hyp:E}
\exists 0<r_0 \leq T,  \forall  0< r \leq r_0 , \quad \lambda(E\cap[T-r,T])>0.
\end{equation}
If  $\omega(t)$ is a $(\delta,\alpha)$-thick subset in $\mathbb{R}^d$ for all $t \in E$,
then the non-autonomous Ornstein-Uhlenbeck equation \eqref{syst_general} is null-controllable on $[0,T]$ from the moving control support $(\omega(t))_{t \in [0,T]}$.
\item[$(ii)$] \emph{\textbf{(Necessary condition).}}
If the non-autonomous Ornstein-Uhlenbeck equation \emph{(\ref{syst_general})} is null-controllable on $[0,T]$ from $(\omega(t))_{t \in [0,T]}$,
then the moving control support satisfies the following integral thickness condition on $[0,T]$,
\begin{equation} \label{Integral_relative_density}
\exists r, \delta>0, \forall x \in \mathbb{R}^d, \quad \int_0^T \lambda \big(B_d(x,r) \cap R(0,T-t) \omega(t) \big) dt \geq \delta>0,
\end{equation}
where $B_d(x,r)$ denotes the open Euclidean ball centered at $x$ with radius $r$, and where $R$ stands for the resolvent of the time-varying linear system 
$$\dot{X}(t)=B(T-t)X(t),$$ 
that is, the solution of the system
\begin{equation} \label{Resolvante}
\left\lbrace \begin{array}{l}
\frac{\partial R}{\partial t_1} (t_1,t_0)=B(T-t_1)R(t_1,t_0), \\
R(t_0,t_0)=I_d.
\end{array}\right.
\end{equation}
\item[$(iii)$] In particular, if $\omega(t)=R(T-t,0)\omega_0$, with $\omega_0$ a fixed subset of $\mathbb{R}^d$,
then the non-autonomous Ornstein-Uhlenbeck equation \eqref{syst_general} is null-controllable on $[0,T]$ from $(\omega(t))_{t \in [0,T]}$, if and only if $\omega_0$ is a thick subset in $\rr^d$.
\end{enumerate}
\end{thm}

\medskip

The assertion $(i)$ in Theorem~\ref{thm_OU_CS_CN} extends the result of \cite{KB_KPS_JMAA} (Theorem~1.3) to non necessarily open and possibly moving control subsets since the geometric assumption (\ref{lop20}) readily implies the thickness condition. We also notice that in the assumptions of Theorem~\ref{thm_OU_CS_CN}, the generalized Kalman rank condition (\ref{Kalman_time}) is only supposed to hold at time $T$, and is allowed to fail for smaller times in the limit of the constraints highlighted by (\ref{Kalman_time11}). Regarding the assumption (\ref{hyp:E}) in assertion $(i)$, and for a given subset $E \subset \mathbb{R}$ with positive Lebesgue measure, we observe that condition (\ref{hyp:E}) holds for almost every time $T \in E$, as it holds in particular for any Lebesgue point $T$ in $E$, that is, points satisfying the condition
\begin{equation}\label{pp1}
\lim_{\substack{r \rightarrow 0\\ r>0}}\frac{1}{2r} \lambda(E\cap[T-r,T+r])=1.
\end{equation}
Indeed, condition (\ref{pp1}) readily implies that 
$$\lim_{\substack{r \rightarrow 0\\ r>0}}\frac{1}{r} \lambda(E\cap[T-r,T])=1.$$
The proof of assertion $(i)$ relies on an adapted Lebeau-Robbiano method and an abstract observability result proved in Theorem~\ref{Meta_thm_AdaptedLRmethod}, which is applied to the adjoint problem. The result of Theorem~\ref{Meta_thm_AdaptedLRmethod} extends the abstract observability result established in~\cite{KB_KPS_JEP} (Theorem~2.1) to the non-autonomous case and under weaker dissipation assumptions allowing a controlled blow-up for small times in the dissipation estimates. This generalization with weaker dissipation assumptions is motivated by various study cases, and is actually needed in Theorem~\ref{Main_result_quad} (even in the framework of fixed control subsets) to derive the null-controllability of 
evolution equations associated to accretive non-selfadjoint quadratic operators with zero singular spaces.

The assertion $(ii)$ in Theorem~\ref{thm_OU_CS_CN}  provides a necessary condition for the null-controllability, which does take into account the transport phenomena induced by the drift term in the Ornstein-Uhlenbeck equations associated to the drift matrices~$B$. Notice that this condition implies in particular that any point $x\in\mathbb{R}^d$ is at distance less than $r$ of the set $R(0,T-t) \omega(t)$
on a time subset whose Lebesgue measure is bounded from below as
$$\forall x\in\mathbb{R}^d, \quad \lambda(I(x)) \geq \frac{\delta}{C_d r^d},$$ 
with $C_d>0$ being the measure of the unit open Euclidean ball in $\rr^d$, and 
$$I(x)=\{ t\in[0,T]:\  B_d(x,r) \cap R(0,T-t) \omega(t) \neq \emptyset \}.$$
The necessary condition $(ii)$ is derived by trying out gaussian explicit solutions in the observability estimate for the adjoint system.
We observe that assertion $(ii)$ turns out to be useful to produce negative null-controllability results and can allow to find out cases when null-controllability may require a positive minimal time as the integral thickness condition can fail for small times while holding for large ones. We refer the reader to the various study cases in Section~\ref{sec:Ex}.

When the moving control subsets comply with the flow associated to the transport part of the Ornstein-Uhlenbeck equations, that is, when $\omega(t)=R(T-t,0)\omega_0$ for all $0 \leq t \leq T$, with $\omega_0$ a fixed subset of $\mathbb{R}^d$, assertion $(iii)$ in Theorem~\ref{thm_OU_CS_CN} provides a necessary and sufficient condition, namely the thickness condition on $\omega_0$ for null-controllability to hold on $[0,T]$. This result allows to directly recover the necessary and sufficient condition for null-controllability of the heat equation established independently in~\cite{Egidi_Veselic} and ~\cite{Gengsheng_Wang2}, as in this case $B=0$ and $\omega(t)=R(T-t,0)\omega_0=\omega_0$ for all $0 \leq t \leq T$.
Regarding the proof of assertion $(iii)$, the necessary result is a direct consequence of $(ii)$ as condition (\ref{Integral_relative_density}) implies that there exist $\delta, r>0$ such that for all $x \in \rr^d$,
$$0<\delta \leq \int_0^T \lambda \big(B_d(x,r) \cap R(0,T-t) \omega(t) \big) dt=T\lambda \big(B_d(x,r) \cap \omega_0 \big),$$
which by equivalence of norms in finite dimension, implies the thickness property of the subset $\omega_0$ in $\rr^d$. On the other hand, if the subset $\omega_0$ is thick in $\rr^d$, there exist some positive constant $\delta_0, r_0>0$ such that 
$$\forall x \in \rr^d, \quad \lambda\big(B_d(x,r_0) \cap \omega_0\big) \geq \delta_0>0.$$
The control subsets $\omega(t)=R(T-t,0)\omega_0$ are then $(\delta,\alpha)$-thick subsets in $\mathbb{R}^d$ for all $0 \leq t \leq T$, since for all $0 \leq t \leq T$ and $x \in \rr^d$,
\begin{align*}
\lambda\big(\{x+[-r,r]^d\} \cap \omega(t)\big) \geq & \ \lambda\big(B_d(x,r) \cap \omega(t)\big)\\
= & \ |\det R(T-t,0)|\lambda\big(\big(R(0,T-t)B_d(x,r)\big) \cap \omega_0\big) \\
\geq & \ c_0 
\lambda\big(B_d(R(0,T-t)x,r_0) \cap \omega_0\big) \geq c_0 \delta_0
\end{align*}
with
$$r=r_0 \sup_{t \in [0,T]}\|R(T-t,0)\|_2>0, \ c_0=\inf_{t \in [0,T]}|\det R(T-t,0)|>0, \ \delta=\frac{c_0 \delta_0}{(2r)^d}>0,$$
and $\alpha=(2r,..., 2r) \in (0,+\infty)^d$. This proves that assertion $(iii)$ is a direct consequence of assertions $(i)$ and $(ii)$ in 
Theorem~\ref{thm_OU_CS_CN}.

The following section provides some applications of Theorem~\ref{thm_OU_CS_CN} with several study cases of Ornstein-Uhlenbeck equations including the detailed analysis of the general autonomous case.
We prove for instance that the Kolmogorov equation 
\begin{equation} \label{Kolm_intro}
(\partial_t + v \partial_x - \partial_v^2 )f(t,x,v)=\un_{\omega}(x,v) u(t,x,v),\quad (t,x,v) \in (0,T) \times\mathbb{R} \times\mathbb{R},
\end{equation}
\begin{itemize}
\item[$(i)$]  is not null-controllable on $[0,T]$, for any arbitrary $T>0$, when the control acts on the fixed control subset $\omega$ made of parallel vertical strips 
$$\begin{array}{c}
\omega = \Big( \underset{n \in \mathbb{Z}}{\bigcup} (n-\eps_{|n|} , n+\eps_{|n|} ) \Big) \times\mathbb{R} \subset \rr_x \times \rr_v, 
\end{array}$$
whose width $0 \leq \eps_n \leq 1$ defines $(\eps_n)_{n \geq 0}$ a non-increasing sequence vanishing at infinity $\lim_{n \to +\infty} \eps_n=0$
\item[$(ii)$]  is not null-controllable on $[0,T]$, for any arbitrary $T>0$, when the control acts on the cone 
\begin{equation} \label{Cone_0_Kolm_trans}
\omega = \big\{ (x, v=\alpha x ) \in \rr^2:\   -\tan \theta_0< \alpha< 0 \big\} \subset \rr_x \times \rr_v \text{ with } 0<\theta_0 <\frac{\pi}{2}
\end{equation}
\item[$(iii)$]  is not null-controllable on $[0,T]$ with $T \leq \frac{2}{\tan \theta_0}$, when the control acts on the cone 
\begin{equation} \label{Cone_Kolm_trans}
\omega=\big\{ (x, v=\alpha x ) \in \rr^2: \  -\tan \theta_0 <\alpha< \tan \theta_0\big\} \subset \rr_x \times \rr_v, 
\end{equation}
with  $0<\theta_0 <\frac{\pi}{2}$
\end{itemize}
We also prove that the Kolmogorov equation with the non-degenerate quadratic external potential $V(x)=\frac{1}{2}x^2$,
\begin{equation} \label{Kolm_quad_1}
(\partial_t + v \partial_x - x \partial_v - \partial_v^2 )f(t,x,v)= \un_{\omega}(x,v) u(t,x,v) ,\quad (t,x,v) \in (0,T) \times\mathbb{R} \times\mathbb{R},
\end{equation}
\begin{itemize}
\item[$(i)$] is not null-controllable on $[0,T]$, for any arbitrary $T>0$, when the control acts on a strip shaped control subset $\omega=\mathbb{R} \times (-L,L) \subset \rr_x \times \rr_v$, with $L>0$
\item[$(ii)$]  is not null-controllable on $[0,T]$ with $T<\pi-\theta_0$ when the control acts on the cone
\begin{equation} \label{Cone_Kolm_rot}
\omega=\big\{(x, v=\alpha x ) \in \rr^2:\    0<\alpha<\tan \theta_0\big\} \subset \rr_x \times \rr_v \text{ with }  0<\theta_0<\frac{\pi}{4}
\end{equation}
\end{itemize}

We also provide in Section~\ref{sec:Ex} examples of moving control supports $(\omega(t))_{t\in I}$ that satisfy the integral thickness condition (\ref{Integral_relative_density}) without being thick subsets in $\rr^d$.
For instance, we prove that:
\begin{itemize}
\item[$(i)$] The cone (\ref{Cone_Kolm_trans}) satisfies the integral thickness condition (\ref{Integral_relative_density}) associated to the Kolmogorov equation (\ref{Kolm_intro}) if and only if $T>\frac{2}{\tan \theta_0}$
\item[$(ii)$] The cone (\ref{Cone_Kolm_rot}) satisfies the integral thickness condition (\ref{Integral_relative_density}) associated to the Kolmogorov equation 
with the quadratic external potential
(\ref{Kolm_quad_1}) when $T> \pi-\theta_0$
\item[$(iii)$] The moving control support $(\omega(t))_{t \in [0,T]}$ defined by 
$$\omega(t) = \omega  \sqrt{1+2\mu t}, \textrm{ where } \omega=[-1,1] \cup \underset{n \geq 1}{\bigcup} (n^2,n^2+n) \cup (-n^2-n,-n^2),$$
with $\mu > 0$, is not thick in $\rr$ for any $0 \leq t \leq T$, but does
satisfy the integral thickness condition (\ref{Integral_relative_density}) on $[0,T]$ associated to the one-dimensional heat equation
$$(\partial_t-\partial_x^2) f(t,x) =\un_{\omega(t)}(x)u(t,x), \quad (t,x) \in (0,T) \times \mathbb{R},$$
for any positive time $T>0$
\end{itemize}

\subsection{Evolution equations associated to accretive non-selfadjoint quadratic operators with zero singular spaces}\label{quadclass}

 \subsubsection{Miscellaneous facts about quadratic operators}
 
Quadratic operators are pseudodifferential operators defined in the Weyl quantization
\begin{equation}\label{3}
q^w(x,D_x) f(x) =\frac{1}{(2\pi)^d}\int_{\rr^{2d}}{e^{i(x-y) \cdot \xi}q\Big(\frac{x+y}{2},\xi\Big)f(y)dyd\xi}, 
\end{equation}
by symbols $q(x,\xi)$, with $(x,\xi) \in \rr^{d} \times \rr^d$, $d \geq 1$, which are complex-valued quadratic forms 
\begin{eqnarray*}
q : \rr_x^d \times \rr_{\xi}^d &\rightarrow& \cc\\
 (x,\xi) & \mapsto & q(x,\xi).
\end{eqnarray*}
These operators are actually differential operators with simple and fully explicit expression since the Weyl quantization of the quadratic symbol 
$x^{\alpha} \xi^{\beta}$, with $(\alpha,\beta) \in \nn^{2d}$, $|\alpha+\beta| = 2$, is given by the differential operator
$$\frac{x^{\alpha}D_x^{\beta}+D_x^{\beta} x^{\alpha}}{2}, \quad D_x=i^{-1}\partial_x.$$
Notice that these operators are non-selfadjoint as soon as their Weyl symbols have a non-zero imaginary part.
The maximal closed realization of the quadratic operator $q^w(x,D_x)$ on $L^2(\rr^d)$, that is, the operator equipped with the domain
\begin{equation}\label{dom1}
D(q^w)=\big\{f \in L^2(\rr^d) : \ q^w(x,D_x)f \in L^2(\rr^d)\big\},
\end{equation}
where $q^w(x,D_x)f$ is defined in the distribution sense, is known to coincide with the graph closure of its restriction to the Schwartz space~\cite{mehler} (pp.~425-426),
$$q^w(x,D_x) : \mathscr{S}(\rr^d) \rightarrow \mathscr{S}(\rr^d).$$
Classically, to any quadratic form $q : \rr_x^d \times \rr_{\xi}^d \rightarrow \cc$ defined on the phase space is associated
a matrix $F \in M_{2d}(\CC)$ called its Hamilton map, or its fundamental matrix, which is defined as the unique matrix satisfying the identity
\begin{equation}\label{10}
\forall  (x,\xi) \in \R^{2d},\forall (y,\eta) \in \R^{2d}, \quad q((x,\xi),(y,\eta))=\sigma((x,\xi),F(y,\eta)),
\end{equation}
where $q(\cdot,\cdot)$ is the polarized form associated with the quadratic form $q$, and where $\sigma$ stands for the standard symplectic form
\begin{equation}\label{11}
\sigma((x,\xi),(y,\eta))=\sum_{j=1}^d(\xi_j y_j-x_j \eta_j),
\end{equation}
with $x=(x_1,...,x_d)$, $y=(y_1,....,y_d)$, $\xi=(\xi_1,...,\xi_d)$ and $\eta=(\eta_1,...,\eta_d) \in \cc^d$.
We observe from the definition that 
$$F=\frac{1}{2}\left(\begin{array}{cc}
\nabla_{\xi}\nabla_x q & \nabla_{\xi}^2q  \\
-\nabla_x^2q & -\nabla_{x}\nabla_{\xi} q 
\end{array} \right),$$
where the matrices $\nabla_x^2q=(a_{i,j})_{1 \leq i,j \leq d}$,  $\nabla_{\xi}^2q=(b_{i,j})_{1 \leq i,j \leq d}$, $\nabla_{\xi}\nabla_x q =(c_{i,j})_{1 \leq i,j \leq d}$,
$\nabla_{x}\nabla_{\xi} q=(d_{i,j})_{1 \leq i,j \leq d}$ are defined by the entries
$$a_{i,j}=\partial_{x_i,x_j}^2 q, \quad b_{i,j}=\partial_{\xi_i,\xi_j}^2q, \quad c_{i,j}=\partial_{\xi_i,x_j}^2q, \quad d_{i,j}=\partial_{x_i,\xi_j}^2q.$$
The notion of singular space introduced in~\cite{kps2} by Hitrik and the third author is defined as the following finite intersection of kernels
\begin{equation}\label{h1bis}
S=\Big( \bigcap_{j=0}^{2d-1}\textrm{Ker}
\big[\textrm{Re }F(\textrm{Im }F)^j \big]\Big)\cap \rr^{2d},
\end{equation}
where $\textrm{Re }F$ and $\textrm{Im }F$ stand respectively for the real and imaginary parts of the Hamilton map $F$ associated with the quadratic symbol $q$,
$$\textrm{Re }F=\frac{1}{2}(F+\overline{F}), \quad \textrm{Im }F=\frac{1}{2i}(F-\overline{F}).$$
As pointed out in \cite{kps2,short,HPSVII,kps11,kps21,rodwahl,viola0}, the notion of singular space plays a basic role in the understanding of the spectral and hypoelliptic properties of (possibly non-elliptic) quadratic operators, as well as the spectral and pseudospectral properties of certain classes of degenerate doubly characteristic pseudodifferential operators~\cite{kps3,kps4,viola1,viola2}. In particular, the work~\cite{kps2} (Theorem~1.2.2) provides a complete description for the spectrum of any non-elliptic quadratic operator $q^w(x,D_x)$ whose Weyl symbol $q$ has a non-negative real part $\textrm{Re }q \geq 0$, and satisfies a condition of partial ellipticity along its singular space~$S$,
\begin{equation}\label{sm2}
(x,\xi) \in S, \ q(x,\xi)=0 \Rightarrow (x,\xi)=0. 
\end{equation}
Under these assumptions, the spectrum of the quadratic operator $q^w(x,D_x)$ is shown to be composed of a countable number of eigenvalues with finite algebraic multiplicities and the structure of this spectrum is similar to the one known for elliptic quadratic operators~\cite{sjostrand}. This condition of partial ellipticity is generally weaker than the condition of ellipticity, $S \subsetneq \rr^{2d}$, and allows one to deal with more degenerate situations. An important class of quadratic operators satisfying condition (\ref{sm2}) are those with zero singular spaces $S=\{0\}$. In this case, the condition of partial ellipticity trivially holds.
More specifically, these quadratic operators have been shown in \cite{kps21} (Theorem~1.2.1) to be hypoelliptic and to enjoy global subelliptic estimates of the type
\begin{equation}\label{lol1}
\exists C>0, \forall u \in \mathscr{S}(\rr^d), \ \|\langle(x,D_x)\rangle^{2(1-\delta)} u\|_{L^2} \leq C(\|q^w(x,D_x) u\|_{L^2}+\|u\|_{L^2}),
\end{equation}
where $\|\cdot\|_{L^2}=\|\cdot\|_{L^2(\rr^d)}$ and $\langle(x,D_x)\rangle^{2}=1+|x|^2+|D_x|^2$, with a sharp loss of derivatives $0 \leq \delta<1$ with respect to the elliptic case (case $\delta=0$), which can be explicitly derived from the structure of the singular space.

In this work, we study the class of quadratic operators whose Weyl symbols have non-negative real parts $\textrm{Re }q \geq 0$, and zero singular spaces $S=\{0\}$. 
These quadratic operators are known~\cite{kps2} (Theorem~1.2.1) to generate contraction semigroups $(e^{-tq^w})_{t \geq 0}$ on $L^2(\rr^d)$, which are smoothing in the Schwartz space for any positive time
$$\forall t>0, \forall f \in L^2(\rr^d), \quad e^{-t q^w}f \in \mathscr{S}(\rr^d).$$
In the recent work~\cite{HPSVII} (Theorem~1.2), these regularizing properties were sharpened and 
these contraction semigroups were shown to be actually smoothing for any positive time in the Gelfand-Shilov space $S_{1/2}^{1/2}(\rr^d)$: $\exists C>0$, $\exists t_0 > 0$, $\forall f\in L^2(\real^d)$, $\forall \alpha, \beta \in \nn^d$, $\forall 0<t \leq t_0$,
\begin{equation}\label{eq1.10}
\|x^{\alpha}\partial_x^{\beta}(e^{-tq^w}f)\|_{L^{\infty}(\rr^d)} \leq \frac{C^{1+|\alpha|+|\beta|}}{t^{\frac{2k_0+1}{2}(|\alpha|+|\beta|+2n+s)}}(\alpha!)^{1/2}(\beta!)^{1/2}\|f\|_{L^2(\rr^d)},
\end{equation}
where $s$ is a fixed integer verifying $s > d/2$, and where $0 \leq k_0 \leq 2d-1$ is the smallest integer satisfying 
\begin{equation}\label{h1bis2}
\Big( \bigcap_{j=0}^{k_0}\textrm{Ker}\big[\textrm{Re }F(\textrm{Im }F)^j \big]\Big)\cap \rr^{2d}=\{0\}.
\end{equation}
The definition and few facts about the Gelfand-Shilov regularity are recalled in appendix (Section~\ref{GSreg}). Notice that the definition (\ref{h1bis2}) makes sense as the singular space is zero $S=\{0\}$. An interesting example of accretive quadratic operator with zero singular space is the Kramers-Fokker-Planck operator acting on $L^2(\rr_{x,v}^2)$,
\begin{equation} \label{eq:KFP}
K=-\Delta_v+\frac{v^2}{4}+v\partial_x-\nabla_xV(x)\partial_v, \quad (x,v) \in \rr^{2},
\end{equation}
with a non-degenerate quadratic external potential $V(x)=\frac{1}{2}x^2$ for which $k_0=1$ in this case.
We refer the reader to the works \cite{KB_KPS_JEP,kps11,kps21} for other examples of accretive quadratic operators with zero singular spaces, as e.g. hypoelliptic Ornstein-Uhlenbeck or Fokker-Planck operators acting on $L^2$-spaces weighted by invariant measures~\cite{KB_KPS_JEP}, quadratic operators appearing in finite-dimensional Markovian approximations of the non-Markovian generalized Langevin equation~\cite{kps11}, or models of chain of oscillators coupled at heat baths at each side~\cite{KB_KPS_JEP,kps11}.

\subsubsection{Sufficient geometric condition for null-controllability}

The main result regarding evolution equations associated to accretive non-selfadjoint quadratic operators with zero singular spaces is the following sufficient geometric condition for null-controllability with moving control supports:

\medskip

\begin{thm}\label{Main_result_quad}
Let $q : \rr_{x}^{d} \times \rr_{\xi}^d \rightarrow \cc$ be a complex-valued quadratic form with a non-negative real part $\emph{\textrm{Re }}q \geq 0$, and a zero singular space $S=\{0\}$. Let $T>0$, $\delta>0$, $\alpha \in (0,+\infty)^d$, $(\omega(t))_{t \in [0,T]}$ be a moving control support in $\mathbb{R}^d$ 
and $E$ be a measurable subset of $[0,T]$ with positive Lebesgue measure $\lambda(E)>0$.
If  $\omega(t)$ is a $(\delta,\alpha)$-thick subset of $\mathbb{R}^d$ for all $t \in E$, then the parabolic equation 
\begin{equation} \label{Eq:quad}
\left\lbrace \begin{array}{ll}
\partial_tf(t,x) + q^w(x,D_x)f(t,x)=\un_{\omega(t)}(x)u(t,x), \quad &  x \in \mathbb{R}^d, \\
f|_{t=0}=f_0 \in L^2(\rr^d),                                       &  
\end{array}\right.
\end{equation}
is null-controllable on $[0,T]$.
\end{thm}

\medskip

The null-controllability from fixed control subsets of evolution equations associated to accretive non-selfadjoint quadratic operators with zero singular spaces was studied in the previous works~\cite{KB_KPS_JEP,KB_KPS_PJ}, and was shown to hold for any fixed control subset satisfying the geometric assumption (\ref{lop20}) in~\cite{KB_KPS_JEP} (Theorem~1.4), and for any fixed thick subset in~\cite{KB_KPS_PJ} (Theorem~2.2). The result of Theorem~\ref{Main_result_quad} therefore extends these previous results in the framework of moving control supports. In the case of fixed control subsets, let us mention that Theorem~\ref{Main_result_quad} actually provides an alternative proof to the one given in~\cite{KB_KPS_PJ}. Indeed, the proof of Theorem~\ref{Main_result_quad} relies on an abstract observability result (Theorem~\ref{Meta_thm_AdaptedLRmethod}), that extends the abstract observability result established in~\cite{KB_KPS_JEP} (Theorem~2.1) from which is derived Theorem~2.2 in~\cite{KB_KPS_PJ}. However, the abstract observability results are applied with different families of orthogonal projections. In the work~\cite{KB_KPS_PJ}, the orthogonal projections at play are the projections onto the first Hermite modes, whereas in the present work, the orthogonal projections are frequency cutoff projections. The spectral estimates (\ref{Meta_thm_IS}) and the dissipation estimates (\ref{Meta_thm_dissip}) are therefore of different kinds in the two proofs, and are linked to smoothing effects in different types of regularity. In~\cite{KB_KPS_PJ}, the dissipation estimates are derived from a Gelfand-Shilov smoothing effect, that is, some Gevrey type smoothing effects both for the solutions and their Fourier transforms. In the present work, the dissipation estimates are derived from a weaker smoothing effect only given by a Gevrey type smoothing effect for the solutions. Regarding the spectral estimates (\ref{Meta_thm_IS}), we use in this work the quantitative version of the Logvinenko-Sereda Theorem established by Kovrijkine~\cite{Kovrijkine}, whereas it was needed in~\cite{KB_KPS_PJ} to derive an adapted version of the Logvinenko-Sereda Theorem for finite combinations of Hermite functions. The proof given in the present work is then somehow more natural than the one derived in~\cite{KB_KPS_PJ}. However, it is actually most interesting to be able to use two different approaches for establishing null-controllability. On one hand, one can indeed consider more degenerate cases of evolution equations associated to accretive non-selfadjoint quadratic operators with possibly non zero singular spaces. This question was recently addressed by Alphonse in~\cite{alphonse} (Theorem~1.12), who obtained some null-controllability results from fixed thick subsets for evolution equations associated to certain classes of quadratic operators with non zero singular spaces that enjoy only partial Gelfand-Shilov smoothing effects.
Theorem~1.12 in~\cite{alphonse} is actually derived from the abstract observability result (Theorem~\ref{Meta_thm_AdaptedLRmethod}) established in the present work and used with frequency cutoff projections. On the other hand, the other approach based on the Gelfand-Shilov smoothing effects and projections onto the first Hermite modes which applies only for evolution equations associated to quadratic operators with zero singular spaces can be push further by taking advantage of the up to now  unused exponential decay of the solutions in order to weaken the thickness assumption of the control support. This is the topic of a work in preparation by Martin and the third author~\cite{jeremy}.

\subsection{Structure of this article}
Section~\ref{sec:Ex} is devoted to provide some applications of Theorem~\ref{thm_OU_CS_CN} and to study specific non-autonomous  Ornstein-Uhlenbeck equations with various moving control supports including the general autonomous case and various examples of non-thick control supports that still satisfy the integral thickness condition (\ref{Integral_relative_density}) for some positive time. The main result in Section~\ref{Sec:LR} is Theorem~\ref{Meta_thm_AdaptedLRmethod} which extends the abstract observability result established in~\cite{KB_KPS_JEP} (Theorem~2.1) to the non-autonomous case and under weaker dissipation assumptions allowing a controlled blow-up for small times in the dissipation estimates.  This abstract result is key in the proofs of Theorems~\ref{thm_OU_CS_CN} and~\ref{Main_result_quad}. Theorem~\ref{thm_OU_CS_CN} is established in Section~\ref{Sec:Proof}, whereas the proof of Theorem~\ref{Main_result_quad} is given in Section~\ref{sec:quad}. The appendix in Section~\ref{appendix} recalls various results needed in the core of this work as the well-posedness of the homogeneous and inhomogeneous Cauchy problems associated to non-autonomous Ornstein-Uhlenbeck equations in Section~\ref{appensix:well-pos}; the Hilbert uniqueness method in the framework of moving control supports in Section~\ref{Appendix:HUM}; some key uncertainty principles related to thick subsets in $\rr^d$ in Section~\ref{Sec:rds}; basic facts and estimates on Hermite functions in Section~\ref{6.sec.harmo}; and the definition and various characterizations of the Gelfand-Shilov regularity in Section~\ref{GSreg}.

\section{Some study cases of Ornstein-Uhlenbeck equations} \label{sec:Ex}

This section provides some applications of Theorem~\ref{thm_OU_CS_CN} with several study cases of Ornstein-Uhlenbeck equations.

\subsection{Heat equation}

By applying Theorem \ref{thm_OU_CS_CN} to the heat equation, that is, the autonomous Ornstein-Uhlenbeck equation with $A(t)=\sqrt{2}I_d$ and $B(t)=0$,
\begin{equation} \label{heat_eq}
\left\lbrace \begin{array}{ll}
(\partial_t - \Delta_x)f(t,x)= \un_\omega(x) u(t,x), \quad (t,x) \in (0,T)\times\mathbb{R}^d, \\
f|_{t=0}=f_0 \in L^2(\mathbb{R}^d),
\end{array}\right.
\end{equation}
we recover the following result established for fixed control subsets $\omega \subset \rr^d$ in~\cite{Egidi_Veselic,Gengsheng_Wang2}:

\medskip

\begin{Cor}\label{cvb1}
Let $\omega$ be a measurable subset of $\mathbb{R}^d$. The following assertions are equivalent:
\begin{enumerate}
\item[$(i)$] The subset $\omega$ is thick in $\rr^d$
\item[$(ii)$] The heat equation \emph{(\ref{heat_eq})} is null-controllable from $\omega$ for some time $T>0$
\item[$(iii)$] The heat equation \emph{(\ref{heat_eq})} is null-controllable from $\omega$ in any time $T>0$
\end{enumerate}
\end{Cor}

\medskip

\begin{proof}
We first notice that the generalized Kalman rank condition (\ref{Kalman_time}) holds at any positive time $T>0$, as $\tilde{A}_0(T)=\sqrt{2}I_d$.
Corollary~\ref{cvb1} is then a direct consequence of assertion $(iii)$ in Theorem~\ref{thm_OU_CS_CN}, as here $B=0$ and $\omega(t)=R(T-t,0)\omega=\omega$ for all $0 \leq t \leq T$.
\end{proof}

\subsection{Abstract autonomous hypoelliptic Ornstein-Uhlenbeck equations}\label{constantOU}

We consider Ornstein-Uhlenbeck equations on $\mathbb{R}^d$ in the autonomous case
\begin{equation}  \label{OU_eq}
\left\lbrace \begin{array}{l}
\partial_t f(t,x) - \frac{1}{2}\textrm{Tr}[Q\nabla_x^2 f(t,x)] - \langle B x , \nabla_x f(t,x) \rangle = \un_{\omega(t)} u(t,x),\\
f|_{t=0}=f_0 \in L^2(\mathbb{R}^d),
\end{array}\right.
\end{equation}
with $Q, B \in M_d(\mathbb{R})$, where $Q$ is a symmetric positive semidefinite matrix. We assume that generalized Kalman rank condition holds at some time $T>0$, that is, 
\begin{equation} \label{Kalman_time34}
\text{Span}\{B^kQ^{\frac{1}{2}}x :\ x\in\mathbb{R}^d, k \geq 0\}=\mathbb{R}^d,
\end{equation}
with $Q^{\frac{1}{2}}$ the symmetric positive semidefinite matrix given by the square root of~$Q$.
In the autonomous case, notice that if the generalized Kalman rank condition holds at some positive time, then it holds at any positive time.
According to the Cayley-Hamilton theorem, condition (\ref{Kalman_time34}) is equivalent to the classical Kalman rank condition
\begin{equation}\label{kal1}
\textrm{Rank}[B|Q^{\frac{1}{2}}]=d,
\end{equation}
where
$$[B|Q^{\frac{1}{2}}]=[Q^{\frac{1}{2}},BQ^{\frac{1}{2}},\dots, B^{d-1}Q^{\frac{1}{2}}],$$
is the $d\times d^2$ matrix obtained by writing consecutively the columns of the matrices $B^jQ^{\frac{1}{2}}$.
By applying Theorem~\ref{thm_OU_CS_CN}, we obtain the following result in the autonomous case:

\medskip

\begin{Cor} \label{Cor_OUa}
Let $T>0$ and $(\omega(t))_{t \in [0,T]}$ be a moving control support in $\mathbb{R}^d$.
We assume that the Kalman rank condition \emph{(\ref{kal1})} holds.
\begin{enumerate}
\item[$(i)$] \emph{\textbf{(Sufficient condition).}}
Let  $\delta>0$, $\alpha \in (0,+\infty)^d$ and $E$ be a measurable subset of $[0,T]$ with positive Lebesgue measure $\lambda(E)>0$.
If  $\omega(t)$ is a $(\delta,\alpha)$-thick subset in $\mathbb{R}^d$ for all $t \in E$,
then the autonomous Ornstein-Uhlenbeck equation \emph{(\ref{OU_eq})} is null-controllable on $[0,T]$ from the moving control support $(\omega(t))_{t \in [0,T]}$.
\item[$(ii)$] \emph{\textbf{(Necessary condition).}}  If the autonomous Ornstein-Uhlenbeck equation \emph{(\ref{OU_eq})} is null-controllable on $[0,T]$ from the moving control support $(\omega(t))_{t \in [0,T]}$, then the moving control support satisfies the following integral thickness condition on $[0,T]$, 
\begin{equation}\label{lop2}
\exists r,\delta>0, \forall x \in \mathbb{R}^d, \quad \int_0^T \lambda \big(B_d(x,r) \cap e^{(t-T)B} \omega(t) \big) dt \geq \delta>0,
\end{equation}
where $B_d(x,r)$ denotes the open Euclidean ball centered at $x$ with radius $r$.
\item[$(iii)$] In particular, if $\omega(t)=e^{(T-t)B}\omega_0$, with $\omega_0$ a fixed subset of $\mathbb{R}^d$,
then the following assertions are equivalent:
\begin{itemize}
\item[$a)$] The subset $\omega_0$ is thick in $\mathbb{R}^d$
\item[$b)$] The autonomous Ornstein-Uhlenbeck equation \eqref{OU_eq} is null-controllable on $[0,T]$ from the moving control support $(\omega(t))_{t \in [0,T]}$ for some positive time~$T>0$
\item[$c)$] The autonomous Ornstein-Uhlenbeck equation \eqref{OU_eq} is null-controllable on $[0,T]$ from the moving control support $(\omega(t))_{t \in [0,T]}$ in any positive time~$T>0$
\end{itemize}
\end{enumerate}
\end{Cor}

\medskip

\begin{proof}
The assertions $(ii)$ and $(iii)$ of Corollary~\ref{Cor_OUa} are a direct rephrasing of the corresponding statements in Theorem~\ref{thm_OU_CS_CN}, since here $R(t_1,t_0)=e^{(t_1-t_0)B}$, with $t_1,t_0 \in \rr$. It therefore remains to prove assertion $(i)$. If $E$ is a subset of $[0,T]$ with positive Lebesgue measure, we can find $0<T'\leq T$ such that 
$$\exists 0<r_0 \leq T',  \forall  0< r \leq r_0 , \quad \lambda(E\cap[T'-r,T'])>0,$$
as it is sufficient to take any Lebesgue point $T'$ of $E$. The Kalman rank condition (\ref{kal1}) implies that the generalized Kalman rank condition (\ref{Kalman_time}) holds at time~$T'$. The assertion $(i)$ in Theorem~\ref{thm_OU_CS_CN} provides the null-controllability 
of the autonomous Ornstein-Uhlenbeck equation \eqref{OU_eq}
on $[0,T']$ from the moving control support $(\omega(t))_{t \in [0,T']}$. We therefore obtain the null-controllability on $[0,T]$  from the moving control support $(\omega(t))_{t \in [0,T]}$ by extending the control function to zero on $[T',T]$.
\end{proof}

\subsection{Translation example with the Kolmogorov equation}

The Kolmogorov equation
\begin{equation} \label{Kolm}
(\partial_t + v \partial_x - \partial_v^2 )f(t,x,v)=\un_{\omega}(x,v) u(t,x,v),\quad (t,x,v) \in (0,T) \times\mathbb{R} \times\mathbb{R},
\end{equation}
with $T>0$, is an autonomous Ornstein-Uhlenbeck equation (\ref{OU_eq}) with the matrices
$$B=\left( \begin{array}{cc}
0 & -1 \\ 0 & 0
\end{array}\right), \quad
Q=\left( \begin{array}{cc}
0 & 0 \\ 0 & 2
\end{array}\right).$$
The Kalman rank condition holds since
$$[Q^{\frac{1}{2}},BQ^{\frac{1}{2}}]=
\left(\begin{array}{cccc}
0 & 0 & 0 & -\sqrt{2} \\
0 & \sqrt{2} & 0 & 0
\end{array}\right).$$
Notice that in this case, the characteristics associated to the drift matrix $B$ follow horizontal lines and are translations along the $x$-axis since
\begin{equation}\label{translation_exp}
\forall  t>0 , \quad e^{ -t B}=\left(\begin{array}{cc}
1 & t \\
0 & 1
\end{array}\right).
\end{equation}
The following result holds true:

\medskip

\begin{Prop}\label{lop3}\ 
\begin{enumerate}
\item[$(i)$] If $\omega$ is a thick subset in $\mathbb{R}^2_{x,v}$, then the Kolmogorov equation \emph{(\ref{Kolm})} is null-controllable on $[0,T]$ from the fixed control subset $\omega$ for any positive time $T>0$
\item[$(ii)$] If $\omega$ is made of parallel vertical strips whose width vanishes to zero at infinity
$$\begin{array}{c}
\omega = \Big( \underset{n \in \mathbb{Z}}{\bigcup} (n-\eps_{|n|} , n+\eps_{|n|} ) \Big) \times\mathbb{R}_v, 
\text{ where } 0 \leq \eps_n \leq 1 \text{ and } \underset{n \to +\infty}{\lim} \eps_n=0,
\end{array}$$
with $(\eps_n)_{n \geq 0}$ being a non-increasing sequence,
then, for any arbitrary time $T>0$, the Kolmogorov equation \emph{(\ref{Kolm})} is not null-controllable on $[0,T]$ from the fixed control subset $\omega$
\item[$(iii)$] If $\omega$ is a cone of the type
$$\omega=\big\{(x, v=\alpha x) \in \rr^2_{x,v}:\   - \tan \theta_0<\alpha< 0\big\} \text{ with }  0<\theta_0<\frac{\pi}{2},$$
then, for any arbitrary time $T>0$, the Kolmogorov equation \emph{(\ref{Kolm})} is not null-controllable on $[0,T]$ from the fixed control subset $\omega$
\item[$(iv)$] If $\omega$ is a cone of the type
$$\omega=\big\{(x, v=\alpha x ) \in \rr^2_{x,v}:\  -\tan \theta_0<\alpha<\tan \theta_0\big\} \text{ with } 0<\theta_0<\frac{\pi}{2},$$
then: 
\begin{itemize}
\item[$1)$] The integral thickness condition \emph{(\ref{lop2})} on $[0,T]$ holds if and only if $T>\frac{2}{\tan \theta_0}$
\item[$2)$] The Kolmogorov equation \emph{(\ref{Kolm})} is not null-controllable on $[0,T]$ from the fixed control subset $\omega$ when
$$0<T \leq \frac{2}{\tan \theta_0}$$
\item[$3)$] The null-controllability of the Kolmogorov equation \emph{(\ref{Kolm})} on $[0,T]$ from the fixed control subset $\omega$ when
$$T>\frac{2}{\tan \theta_0},$$ 
is not covered by the previous results and is an open problem
\end{itemize}
\end{enumerate}
\end{Prop}

\medskip

\begin{proof} The assertion $(i)$ in Proposition~\ref{lop3} is a consequence of $(i)$ in Corollary~\ref{Cor_OUa}. Regarding the proof of assertion $(ii)$ in Proposition~\ref{lop3}, we first observe that the strip-shaped control subset $\omega$ is not thick in $\mathbb{R}^2$, since
$$\forall N \geq 1, \forall n \geq N, \ \lambda\big(  \{(n,0)+[-N,N]^2\} \cap \omega \big) \leq \sum_{j=n-N}^{n+N} 4 \eps_{j} N\leq 4N(2N+1)\eps_{n-N},$$
tends to zero when $n \rightarrow +\infty$. The sufficient condition given by assertion $(i)$ in Proposition~\ref{lop3} therefore does not hold. 
According to (\ref{translation_exp}), the subset $e^{-t B}\omega$ is composed by parallel strips with the angle $\arctan t$ with respect to the $v$-axis as shown in Figure~\ref{Fig1}. It follows that 
for all $t \geq 0$, $N \geq 1$ and $n \geq ([t]+2)N+1$, 
\begin{multline*}
 \lambda\big(\{(n,0)+[-N,N]^2\}\cap e^{ -t B}\omega\big) \leq \sum_{j=n-([t]+2)N-1}^{n+([t]+2)N+1} 4 \eps_{j} N \\ \leq 4N\big(2([t]+2)N+3\big)\eps_{n-([t]+2)N-1},
\end{multline*}
where $[\cdot]$ denotes the floor function. The necessary condition for the null controllability (\ref{lop2}) therefore does not hold on $[0,T]$ for any $T>0$, as
for all $T>0$, $N \geq 1$ and $n \geq ([T]+2)N+1$   
$$0 \leq \int_0^T \lambda \big(\{(n,0)+[-N,N]^2\} \cap e^{-t B}\omega\big) dt \leq 4TN\big(2([T]+2)N+3\big)\eps_{n-([T]+2)N-1},$$
tends to zero when $n \rightarrow +\infty$. For any arbitrary $T>0$, the Kolmogorov equation (\ref{Kolm}) is therefore not null-controllable on $[0,T]$ from the fixed control subset $\omega$.
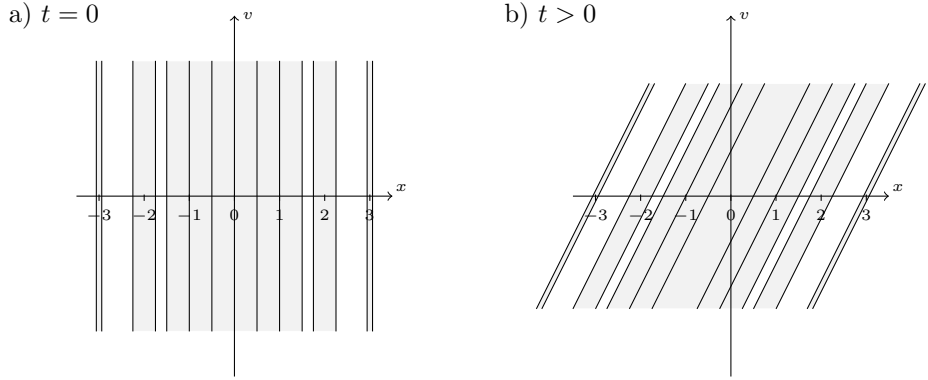
\begin{figure}[H]
\centering
\begin{tikzpicture}[scale=0.6]
\pgfmathsetmacro{\YY}{4}; 
\pgfmathsetmacro{\X}{4}; 
\pgfmathsetmacro{\T}{tan(45)}; 
\begin{scope}[xshift=-16cm] 
    \node at (-4,4)  {a) $t=0$}; 
	\filldraw[color=white, fill=black!5] (-1,3) rectangle (1,-3); 
	\foreach \x in {1,-1}{
		\filldraw[color=white, fill=black!5] (\x -0.5,3) rectangle (\x +0.5, -3); 
		}
	\foreach \x in {2,-2}{
		\filldraw[color=white, fill=black!5] (\x-0.25,3) rectangle (\x +0.25, -3); 
		}
	\foreach \x in {3,-3}{
		\filldraw[color=white, fill=black!5] (\x-0.06, 3) rectangle (\x+0.06, -3); 
		}
	\foreach \x in {1,-1}{
		\draw (\x,3) -- (\x, -3);
		}
	\foreach \x in {0.5, 1.5}{
		\draw (\x,3) -- (\x, -3);
		\draw (-\x, 3) -- (-\x, -3);
		}
	\foreach \x in {1.75, 2.25}{
		\draw (\x, 3) -- (\x, -3);
		\draw (-\x, 3) -- (-\x, -3);
		}
	\foreach \x in {2.94, 3.06}{
		\draw (\x, 3) -- (\x, -3); 
		\draw (-\x, 3) -- (-\x, -3);
		}
	\draw[->] (-3.5,0) -- coordinate (x axis mid) (3.5,0); 
	\draw (3.7,0.2) node {\tiny$x$};
	\foreach \x in {-3,-2,-1,0,1,2,3}{
		\draw (\x ,1pt) -- (\x ,-3pt) 
		node[anchor=north] {\tiny$\x$};}
	\draw[->] (0,-\YY) -- (0,\YY); 
	\draw (0.3, \YY) node {\tiny$v$};
\end{scope}
%
\begin{scope}[xshift=-5cm]
\node at (-4,4) {b) $t>0$};
\pgfmathsetmacro{\Y}{2.5};
\pgfmathsetmacro{\t}{\T * 0.5};
\pgfmathsetmacro{\TT}{\T-\t}; 
\filldraw[color=white, fill=black!5] (-1+\TT *\Y,\Y) -- (1+\TT *\Y ,\Y)-- (+1-\TT *\Y, -\Y) -- (-1-\TT * \Y, -\Y); 
\foreach \x in {1,-1}{ 
	\filldraw[color=white, fill=black!5] (\x -0.5+\TT *\Y,\Y) -- (\x +0.5 +\TT *\Y ,\Y)-- (\x +0.5-\TT *\Y, -\Y) -- (\x -0.5-\TT * \Y, -\Y);
}
\foreach \x in {2,-2}{
	\filldraw[color=white, fill=black!5] (\x-0.25+\TT *\Y,\Y) -- (\x +0.25+\TT *\Y ,\Y)-- (\x +0.25-\TT *\Y, -\Y) -- (\x -0.25-\TT * \Y, -\Y); 
}
\foreach \x in {-3,3}{ 
	\filldraw[color=white, fill=black!5] (\x-0.06+\TT *\Y, \Y) -- (\x +0.06+\TT *\Y ,\Y)-- (\x +0.06-\TT *\Y, -\Y) -- (\x -0.06-\TT * \Y, -\Y);
}
\foreach \x in {-1,1}{
	\draw (\x + \TT * \Y,\Y) -- (\x -\TT * \Y, -\Y);
}
\foreach \x in {1,-1}{
	\draw (\x -0.5+\TT * \Y,\Y) -- (\x -0.5-\TT * \Y, -\Y);
	\draw (\x +0.5+\TT * \Y, \Y) -- (\x +0.5-\TT * \Y, -\Y);
}
\foreach \x in {2,-2}{
	\draw (\x -0.25+\TT * \Y,\Y) -- (\x -0.25-\TT * \Y, -\Y);
	\draw (\x +0.25+\TT * \Y, \Y) -- (\x +0.25-\TT * \Y, -\Y);
}
\foreach \x in {3,-3}{
	\draw (\x -0.06+\TT * \Y,\Y) -- (\x -0.06-\TT * \Y, -\Y);
	\draw (\x +0.06+\TT * \Y, \Y) -- (\x +0.06-\TT * \Y, -\Y);
}
\draw[->] (-3.5,0) -- coordinate (x axis mid) (3.5,0); 
\draw (3.7,0.2) node {\tiny $x$};
\foreach \x in {-3,-2,-1,0,1,2,3}{
	\draw (\x ,1pt) -- (\x ,-3pt) 
	node[anchor=north] {\tiny$\x$};}
\draw[->] (0,-\YY) -- (0,\YY); 
\draw (0.3, \YY) node {\tiny$v$};
\end{scope}

\end{tikzpicture}
\caption{Motion of the union of strips with parameter $\varepsilon_{n} = 2^{-n}$, $n \geq 0$ 
under rotation of angle $\arctan t$ with respect to the $v$-axis.}\label{Fig1}
\end{figure}

\noindent
We now give a proof of assertion $(iii)$. Let $0<\theta_0< \frac{\pi}{2}$, $T>0$ and $r>0$.
We denote by 
\begin{equation} \label{def:D1}
\mathcal{D}_1=\{ (x,v)\in\rr^2: \ v=-x\tan \theta_0  \}
\text{  and  }
\mathcal{D}_0=\{ (x,0): \ x \in \mathbb{R}\},
\end{equation}
the straight lines composing the boundary of the conic subset $\omega$.
We notice that $(e^{-tB} \omega)_{t \in [0,T]}$ is an increasing family of cones whose boundaries are given by the straight lines $e^{-tB} \mathcal{D}_1$ and $e^{-tB} \mathcal{D}_0 = \mathcal{D}_0$.
For any $t>0$, the line 
\begin{equation} \label{exp(tB)D1}
e^{-tB} \mathcal{D}_1 = \Big\{ (x,v)\in \rr^2: \  x= \frac{-1+t \tan \theta_0}{\tan \theta_0} v \Big\},
\end{equation}
never coincides with $\mathcal{D}_0$.
In particular, the subset
$$\mathbb{R}^2 \setminus \underset{t \in[0,T]}{\cup} e^{-tB}\overline{\omega} = \mathbb{R}^2 \setminus e^{-TB}\overline{\omega},$$
is a non-empty open cone. By picking $X=(x,v) \in \rr^2$ with a sufficiently large norm in this set,
we can ensure that 
$$B_d(X,r) \cap e^{-tB} \omega=\emptyset,$$ 
for all $t\in[0,T]$. Such a point $X$ violates the integral thickness condition (\ref{lop2}) in time $T$. 
For any arbitrary $T>0$, the Kolmogorov equation (\ref{Kolm}) is therefore not null-controllable on $[0,T]$ from the fixed control subset $\omega$.

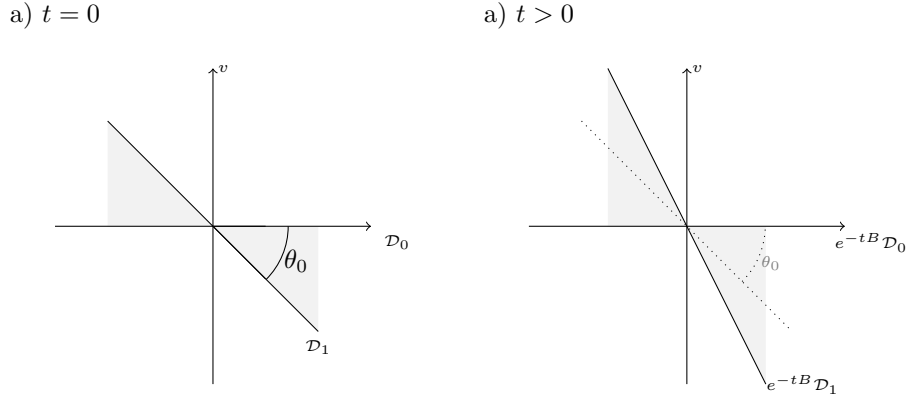
\begin{figure}[H]
	\begin{tikzpicture}[scale=0.7]
	\begin{scope}[xshift=-9cm]
	\node at (-3,4) {a) $t=0$};
	\filldraw[color=white, fill=black!5] (-2, 0) -- (0,0) -- (-2,2);
	\filldraw[color=white, fill=black!5] (2,-2) -- (0,0) -- (2,0);
	\draw
		(1,-1) coordinate (a) 
		-- (0,0) coordinate (b) 
		-- (1,0) coordinate (c) 
		pic["$\theta_0$", draw=black,  angle eccentricity=1.2, angle radius=1cm]
		{angle=a--b--c};
	\draw (-2,2) -- (2,-2);
	\node at (2, -2.3) {\tiny$\mathcal{D}_1$};
	\draw[->] (-3,0) -- (3,0); 
	\node at (3.5, -0.3) {\tiny$\mathcal{D}_0$};
	\draw[->] (0,-3) -- (0,3); 
	\node at (0.2, 3) {\tiny$v$};
	\end{scope}
	
	\begin{scope}
	\node at (-3,4) {a) $t>0$};
	\filldraw[color=white, fill=black!5] (-1.5, 0) -- (0,0) -- (-1.5,3);
	\filldraw[color=white, fill=black!5] (1.5,-3) -- (0,0) -- (1.5,0);
	\draw[ dotted] (1.05,-1.05) arc (315:360:1.5);
	\node[color=black!50] at (1.6, -0.7) {\tiny$\theta_0$};
	\draw[dotted] (-2,2) -- (2,-2);
	\draw (-1.5,3) -- (1.5,-3);
	\node at (2.2, -3) {\tiny$e^{-tB}\mathcal{D}_1$};
	\draw[->] (-3,0) -- (3,0); 
	\node at (3.5, -0.3) {\tiny$e^{-tB}\mathcal{D}_0$};
	\draw[->] (0,-3) -- (0,3); 
	\node at (0.2, 3) {\tiny$v$};
	\end{scope}
	\end{tikzpicture}
	\bigskip
	\caption{Motion of the cone under rotation when $t>0$.}\label{Fig2}
\end{figure}

\noindent
We now give a proof of assertion $(iv)$. Let $0<\theta_0 <\frac{\pi}{2}$. We denote  
$$\mathcal{D}_1=\big\{(x,v)\in\rr^2:\ v= -x \tan \theta_0\big\}, \qquad \mathcal{D}_2=\big\{(x,v)\in\rr^2:\ v= x \tan\theta_0\big\},$$
the straight lines composing the boundary of the conic subset $\omega$. We notice that $(e^{-tB} \omega)_{t \in [0,T]}$ is a family of cones whose boundaries are given by the straight lines $e^{-tB} \mathcal{D}_1$ and $e^{-tB} \mathcal{D}_2$. For any $t >0$, the straight line $e^{-tB} \mathcal{D}_1$ is given by (\ref{exp(tB)D1}), whereas 
$$e^{-tB} \mathcal{D}_2 =\Big\{(x,v)\in \rr^2: \ x=\frac{1+t \tan \theta_0}{\tan \theta_0}v \Big\},$$
is a straight line with positive slope with the $x$-axis and converging to the $x$-axis when $t \to +\infty$.

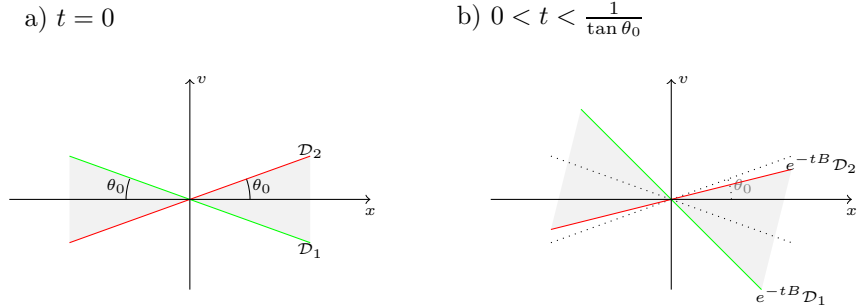
\begin{figure}[H]
\begin{tikzpicture}[scale=0.8]
\begin{scope}[xshift=-8cm]
\node at (-2,3) {a) $t=0$};
\filldraw[color=white, fill=black!5] (-2, -0.72) -- (0,0) -- (-2,0.72);
\filldraw[color=white, fill=black!5] (2,0.72) -- (0,0) -- (2,-0.72);1
\draw (1,0) arc (0:20:1); 
\node at (1.2, 0.22) {\tiny$\theta_0$};
\draw (-1,0.34) arc (160:180:1);
\node at (-1.22, 0.22) {\tiny$\theta_0$};
\draw[red] (-2, -0.72) -- (2,0.72);
\node at (2, 0.8) {\tiny$\mathcal{D}_2$};
\draw[green] (-2,0.72) -- (2,-0.72);
\node at (2, -0.85) {\tiny$\mathcal{D}_1$};
\draw[->] (-3,0) -- (3,0); 
\node at (3, -0.2) {\tiny$x$};
\draw[->] (0,-1.5) -- (0,2); 
\node at (0.2, 2) {\tiny$v$};
\end{scope}
\begin{scope}
\node at (-2,3) {b) $0<t<\frac{1}{\tan \theta_0}$};
\filldraw[color=white, fill=black!5] (2,0.5) -- (0,0) -- (1.5, -1.5);
\filldraw[color=white, fill=black!5] (-1.5, 1.5) -- (0,0) -- (-2,-0.5);
\draw[dotted] (-2, -0.72) -- (2,0.72);
\draw[dotted] (-2,0.72) -- (2,-0.72);
\draw[ dotted] (1,0) arc (0:20:1);
\node[color=black!50] at (1.2, 0.22) {\tiny$\theta_0$};
\draw[green] (1.5,-1.5) -- (-1.5,1.5);
\node at (2, -1.6) {\tiny$e^{-tB}\mathcal{D}_1$};
\draw[red] (-2,-0.5) -- (2, 0.5);
\node at (2.5, 0.6) {\tiny$e^{-tB}\mathcal{D}_2$};
\draw[->] (-3,0) -- (3,0); 
\node at (3, -0.2) {\tiny$x$};
\draw[->] (0,-1.5) -- (0,2); 
\node at (0.2, 2) {\tiny$v$};
\end{scope}
\end{tikzpicture}
\bigskip
\caption{Family of cones $(e^{-tB}\omega)_{t\in [0,T]}$.}\label{Fig3.1}
\end{figure}

\noindent
We observe that the set $\rr^2 \setminus e^{-tB} \omega$ is a non-empty cone for any time $t \geq 0$. As above, this implies in particular that for any $t \geq 0$, $e^{-tB} \omega$ is not a thick subset in $\rr^2$ and that assertion $(i)$ in Corollary~\ref{Cor_OUa} does not apply.

Let $0< T \leq \frac{2}{\tan \theta_0}$, $r>0$ and $y=(x,v) \in \rr^2$ be such that the open Euclidean ball $B_d(y,r)$ is contained in the cone $\rr^2 \setminus e^{-TB} \omega$ and tangent to $\mathcal{D}_2$ as in Figure~\ref{Fig3.2}.
This ball intersects the set $e^{-tB}\omega$  only on a time interval $t \in [0,\ell(y)) \subset [0,T]$ whose length $\ell(y)$ tends to zero when $|y| \rightarrow +\infty$. It follows that
\begin{multline*}
\int\limits_0^T \lambda \big(B_d(y, r) \cap e^{(t-T)B} \omega\big)dt=\int\limits_0^T \lambda \big(B_d(y, r) \cap e^{-tB} \omega\big)dt\\
= \int\limits_{0}^{\ell(y)}\lambda\big(B_d(y,r) \cap e^{-tB} \omega\big) dt 
\leq \lambda\big(B_d(0, r)\big) \ell(y)   \underset{|y| \rightarrow +\infty}{\longrightarrow} 0.
\end{multline*}
The integral thickness condition (\ref{lop2}) is therefore violated on $[0,T]$ when $0< T \leq \frac{2}{\tan \theta_0}$. It follows from the assertion $(ii)$ in Corollary~\ref{Cor_OUa} that the Kolmogorov equation (\ref{Kolm}) is not null-controllable on $[0,T]$ from the fixed control subset $\omega$, when $0< T \leq \frac{2}{\tan \theta_0}$.

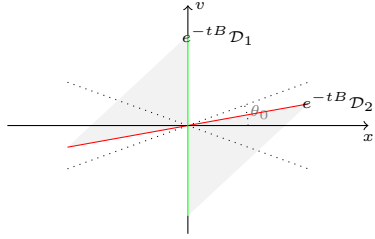
\begin{figure}[H]
\begin{tikzpicture}[scale=0.8]
\begin{scope}[xshift=-8cm, yshift=-6cm]
\node at (-2,3) {c) $t=\frac{1}{\tan\theta_0}$};
\filldraw[color=white, fill=black!5] (0,1.5) -- (0,0) -- (-2,-0.36);
\filldraw[color=white, fill=black!5] (0,-1.5) -- (0,0) -- (2,0.36);
\draw[dotted] (-2, -0.72) -- (2,0.72);
\draw[dotted] (-2,0.72) -- (2,-0.72);
\draw[ dotted] (1,0) arc (0:20:1);
\node[color=black!50] at (1.2, 0.22) {\tiny$\theta_0$};
\draw[red] (-2,-0.36) -- (2,0.36);
\node at (2.5, 0.4) {\tiny$e^{-tB}\mathcal{D}_2$};
\draw[->] (-3,0) -- (3,0); 
\draw[green] (0,-1.5) -- (0,1.5);
\node at (0.5, 1.5) {\tiny$e^{-tB}\mathcal{D}_1$};
\node at (3, -0.2) {\tiny$x$};
\draw (0,-1.5)-- (0,-1.8);
\draw[->] (0,1.5)-- (0,2); 
\node at (0.2, 2) {\tiny$v$};
\end{scope}
\begin{scope}[yshift=-6cm]
\node at (-2,3) {d) $\frac{1}{\tan\theta_0}<t\leq\frac{2}{\tan\theta_0}$};
\filldraw[color=white, fill=black!5] (3, 2.16) -- (0,0) -- (-2, -0.28);
\filldraw[color=white, fill=black!5] (3,0.42) -- (0,0) -- (-2,-1.44);
\draw[dotted] (-2, -0.72) -- (3,1.08);
\draw[dotted] (-2,0.72) -- (3,-1.08);
\draw[ dotted] (1,0) arc (0:20:1);
\node[color=black!50] at (1.2, 0.22) {\tiny$\theta_0$};
\draw[green] (3,2.16) -- (-2,-1.44);
\node at (3.2,2.25) {\tiny$e^{-tB}\mathcal{D}_1$}; 
\draw (2.5, 0.65) circle (0.23cm);
\filldraw (2.5,0.65) circle (0.02cm);
\draw (2.5, 0.65) -- (2.42, 0.871); 
\node at (3.3, 1.2) {\tiny$\mathcal{D}_2$};
\draw[red] (3,0.42) -- (-2, -0.28);
\node at (3.5, 0.3) {\tiny$e^{-tB}\mathcal{D}_2$};
\draw[->] (-3,0) -- (3.5,0); 
\node at (3, -0.2) {\tiny$x$};
\draw[->] (0,-2) -- (0,2); 
\node at (0.2, 2) {\tiny$v$};
\end{scope}
\end{tikzpicture}
\bigskip
\caption{Family of cones $(e^{-tB}\omega)_{t\in [0,T]}$.}\label{Fig3.2}
\end{figure}

\begin{figure}[H]
\begin{tikzpicture}[scale=0.8]
\begin{scope}[xshift=-8cm,yshift=-12cm]
\node at (-2,3) {e) $t>\frac{2}{\tan\theta_0}$};
\filldraw[color=white, fill=black!5] (5,1.3) -- (0,0) --(-2, -0.06);
\filldraw[color=white, fill=black!5] (-2,-0.52) -- (0,0) --(5,0.15);
\draw[dotted] (-2, -0.72) -- (2,0.72);
\draw[dotted] (-2,0.72) -- (2,-0.72);
\draw[ dotted] (1,0) arc (0:20:1);
\draw[red] (5,0.15) -- (-2,-0.06);
\node at (5.5, 0.16) {\tiny$e^{-tB}\mathcal{D}_2$};
\draw[green] (-2, -0.52) -- (5, 1.3);
\node at (5.5, 1.43) {\tiny$e^{-tB}\mathcal{D}_1$};
\node[color=black!50] at (1.2, 0.22) {\tiny$\theta_0$};
\draw[->] (-3,0) -- (7,0); 
\node at (7, -0.2) {\tiny$x$};
\draw[->] (0,-1.5) -- (0,2); 
\node at (0.2, 2) {\tiny$v$};
\end{scope}
\end{tikzpicture}
\bigskip
\caption{Family of cones $(e^{-tB}\omega)_{t\in [0,T]}$.}\label{Fig3.3}
\end{figure}

\noindent
Let $T>\frac{2}{\tan \theta_0}$. We aim at proving that the integral thickness condition (\ref{lop2}) holds on $[0,T]$.
Let $\frac{2}{\tan\theta_0}<T_1< T$, 
$\theta_1$ be the angle between the $x$-axis and the straight line $e^{-TB}\mathcal{D}_1$,
$\widetilde{\theta}_1$ be the angle between the $x$-axis and the straight line $e^{-T_1B}\mathcal{D}_1$. We observe that
 $0<\theta_1<\widetilde{\theta}_1<\theta_0$. We can choose the parameter $T_1$ close enough to $\frac{2}{\tan\theta_0}$ so that there exists
 $0<T_2<T$ such that 
$$e^{-T_2 B} \mathcal{D}_2 = e^{-T_1B}\mathcal{D}_1.$$
For any $X=(x,v)\in\mathbb{R}^2 \setminus\{0\}$ with $x\geq 0$, we introduce
$$\arg(X)=\left\lbrace\begin{array}{ll}
\arctan(\frac{v}{x}) \quad & \text{ if } x>0, \\
-\frac{\pi}{2} & \text{ if } x=0,\ v<0, \\
\frac{\pi}{2} & \text{ if } x=0,\ v>0.
\end{array}\right.$$
We notice that
\begin{itemize}
\item[-] if $\arg(X)\in[-\theta_0,\widetilde{\theta}_1]$, then $X$ belongs to $e^{-tB}\overline{\omega}$ at least for all $t\in[0,T_2]$
\item[-] if $\arg(X)\in[\widetilde{\theta}_1,\frac{\pi}{2}]$, then $X$ belongs to $e^{-tB}\overline{\omega}$ at least for all $t \in [T_1,T]$
\item[-] if $\arg(X)\in[-\frac{\pi}{2},-\theta_0]$, then $X$ belongs to $e^{-tB}\overline{\omega}$ at least for all $t \in [\frac{1}{\tan \theta_0},T]$
\end{itemize}
We set $T^*=\min\{ T_2, T-T_1, T-\frac{1}{\tan\theta_0} \}>0$. By changing $X$ to $-X$, we easily notice that any point $X=(x,v)\in\mathbb{R}^2 \setminus\{0\}$ belongs to $e^{-tB}\omega$ at least for all $t$ in an interval $\mathcal{I}(X) \subset [0,T]$ of length $|\mathcal{I}(X)|\geq T_*$. On the other hand, the point $X=0$ belongs to the closed cone given by the adherence of $e^{-tB}\omega$ for all $0 \leq t \leq T$.
By the conic structure of the set $e^{-tB} \omega$, we can find $\rho>0$ such that for all $X\in\mathbb{R}^2$ with $|X|>\rho$, and for all $t \in \mathcal{I}(X)$, $e^{-tB}\omega$ contains at least half of the ball $B_d(X,1)$. This implies in particular that for all $X\in\mathbb{R}^2$, $|X|>\rho$,
\begin{multline*}
\int_0^T \lambda\big(B_d(X,1) \cap e^{(t-T)B}\omega\big)dt=\int_0^T \lambda\big(B_d(X,1) \cap e^{-tB}\omega\big)dt \\
\geq \frac{1}{2} \lambda\big(B_d(0,1)\big) |\mathcal{I}(X)| \geq \frac{\pi T_*}{2}>0.
\end{multline*}
By continuity of the translation in $L^1$, we notice that the function
$$X \in\mathbb{R}^2 \mapsto \int_0^T \lambda \big( B_d(X,1) \cap e^{-tB}\omega \big) dt,$$
is continuous on $\rr^2$.
Furthermore, this function is positive as any point $X \in \mathbb{R}^2$ belongs to $e^{-tB}\omega$ at least for all $t$ in an interval $\mathcal{I}(X) \subset [0,T]$ of length $|\mathcal{I}(X)|\geq T_*$.
By compactness of the set $\overline{B_d(0,\rho)}$, the above function is therefore bounded from below by a positive constant on $\rr^2$,
$$\exists \delta>0, \forall X \in \rr^2, \quad  \int_0^T \lambda\big(B_d(X,1) \cap e^{(t-T)B}\omega\big)dt=\int_0^T \lambda \big( B_d(X,1) \cap e^{-tB}\omega \big) dt \geq \delta.$$
The integral thickness condition (\ref{lop2}) then holds on $[0,T]$ with $r=1$ and $\delta$. This ends the proof of Proposition~\ref{lop3}.
\end{proof}

\subsection{Rotation example with the Kolmogorov equation with a non-degenerate quadratic external potential}

The Kolmogorov equation with the non-degenerate quadratic external potential $V(x)=\frac{1}{2}x^2$, writes as
\begin{multline} \label{Kolm_xv}
(\partial_t + v \partial_x - \partial_xV(x) \partial_v - \partial_v^2 )f(t,x,v)\\
=(\partial_t + v \partial_x - x \partial_v - \partial_v^2 )f(t,x,v)= \un_{\omega}(x,v) u(t,x,v),
\end{multline}
with $(t,x,v) \in (0,T) \times\mathbb{R} \times\mathbb{R}$,
is an autonomous Ornstein-Uhlenbeck equation (\ref{OU_eq}) with the matrices
$$B=\left( \begin{array}{cc}
0 & -1 \\ 1 & 0
\end{array}\right), \quad
Q=\left( \begin{array}{cc}
0 & 0 \\ 0 & 2
\end{array}\right).$$
The Kalman condition holds since
$$[Q^{\frac{1}{2}},BQ^{\frac{1}{2}}]=
\left(\begin{array}{cccc}
0 & 0 & 0 & -\sqrt{2} \\
0 & \sqrt{2} & 0 & 0
\end{array}\right).$$
Notice that in this case the characteristics are given by rotations as
\begin{equation} \label{rotation_exp}
\forall  t>0 , \quad e^{- t B}=\left(\begin{array}{cc}
\cos t & \sin t \\
-\sin t & \cos t
\end{array}\right).
\end{equation}

\medskip

\begin{Prop}\label{fg123}\
\begin{enumerate}
\item[$(i)$] If the subset $\omega$ is thick in $\mathbb{R}^2$, then the Kolmogorov equation with external potential \emph{(\ref{Kolm_xv})} is null-controllable on $[0,T]$ from the fixed control subset $\omega$ for any positive time $T>0$
\item[$(ii)$] If $\omega$ is  made of a strip
$$\omega = \mathbb{R} \times (-L,L) \subset \rr^2_{x,v},$$
with $L>0$, then for any arbitrary $T>0$, the Kolmogorov equation with external potential \emph{(\ref{Kolm_xv})} is not null-controllable on $[0,T]$ from the fixed control subset $\omega$
\item[$(iii)$] If $\omega$ is a cone of the type
$$\omega = \big\{(x,v=\alpha x) \in \rr^2:\   0< \alpha <\tan \theta_0\big\},$$
with   $0<\theta_0< \frac{\pi}{4}$, 
then 
\begin{itemize}
\item[$a)$] The integral thickness condition \emph{(\ref{lop2})} on $[0,T]$ associated to the
Kolmogorov equation with external potential \emph{(\ref{Kolm_xv})} fails when $T < \pi-\theta_0$, and 
holds when $T > \pi-\theta_0$
\item[$b)$] The Kolmogorov equation with external potential \emph{(\ref{Kolm_xv})} is not null-controllable 
on $[0,T]$ from the fixed control subset $\omega$ when $T<\pi-\theta_0$
\item[$c)$] The null-controllability of the Kolmogorov equation with external potential \emph{(\ref{Kolm_xv})} 
on $[0,T]$ from the fixed control subset $\omega$ when $T \geq \pi-\theta_0$ is an open problem
\end{itemize}
\end{enumerate}
\end{Prop}

\medskip

\begin{proof}
The assertion $(i)$ is a direct consequence of assertion $(i)$ in Corollary~\ref{Cor_OUa}. We now give a proof of assertion $(ii)$. The strip-shaped control subset 
$$\omega = \mathbb{R} \times (-L,L) \subset \rr^2_{x,v},$$
is obviously not thick in $\mathbb{R}^2$. The sufficient condition $(i)$ in Proposition~\ref{fg123} therefore does not hold. According to (\ref{rotation_exp}), the subset $e^{-tB}\omega$ is a strip with width $2L$ and angle $-t$ with respect to the $x$-axis (see Figure~\ref{Fig4}).
Let $T>0$ and $N$ be the smallest integer satisfying $\frac{T}{\pi} \leq N$.
Let $r>0$ and $y \in \rr^2 \setminus \{0\}$.
For any $\mu>0$, the strip $e^{-tB}\omega$ intersects the ball $B_d(\mu y, r)$ on a union of at most $N$ time intervals 
$$t\in \mathcal{I}(\mu) = I_1(\mu)\cup...\cup I_{N}(\mu),$$
with length $| I_{k}(\mu) |$ converging to zero when $\mu \rightarrow +\infty$. It follows that
\begin{multline*}
\int_0^T \lambda \big( B_d(\mu y,r) \cap e^{(t-T)B}\omega \big) dt =\int_0^T \lambda \big( B_d(\mu y,r) \cap e^{-tB}\omega \big) dt \\
=\int_{t \in \mathcal{I}(\mu) } \lambda \big( B_d(\mu y,r) \cap e^{-tB}\omega \big) dt
\leq \lambda \big( B_d(0,r) \big) |\mathcal{I}(\mu)| \underset{\mu \rightarrow +\infty}{\longrightarrow} 0.
\end{multline*}
The integral thickness property (\ref{lop2}) therefore does not hold on $[0,T]$. The Kolmogorov equation with external potential (\ref{Kolm_xv}) is then not null-controllable on $[0,T]$ from the fixed control subset $\omega$.

\begin{figure}[H]
\centering
\begin{tikzpicture}[scale=0.7]
\begin{scope}[xshift=-10cm]
\node at (-3,4) {a) $t=0$};
\filldraw[color=white, fill=black!5] (-3,1) rectangle (3,-1); 
\draw (-3,-1) -- (3,-1); 
\node at (0.3,1.3) {$L$};
\draw (-3,1) -- (3, 1); 
\node at (0.3,-1.3) {$-L$};
\draw[->] (-3.5,0) -- (3.5,0); 
\node at (3.5, -0.2) {\tiny$x$};
\draw[->] (0,-3.5) -- (0,3.5); 
\node at (0.2, 3.5) {\tiny$v$};
\end{scope}

\begin{scope} 
\begin{scope}[rotate=140]  
\filldraw[color=white, fill=black!5] (-3,1) rectangle (3,-1);  
\draw (-3,-1) -- (3,-1); 
\draw (-3,1) -- (3, 1); 
\draw [<->] (3.3, -1) -- (3.3, 1); \node at (3.8, 0) {$2L$};
\end{scope}
\draw (2.45,0) arc (0:-40:0.9); \node at (2.6,-0.4) {$t$};
\node at (-3,4) {b) $t>0$};
\draw[->] (-3.5,0) -- (3.5,0); 
\node at (3.5, -0.2) {\tiny$x$};
\draw[->] (0,-3.5) -- (0,3.5); 
\node at (0.2, 3.5) {\tiny$v$};
\end{scope}

\end{tikzpicture}
\bigskip
\caption{Motion of the strip $\mathbb{R}\times (-L,L)$ under a rotation of angle $-t$ with respect to the $x$-axis.}\label{Fig4}
\end{figure}

\noindent
We now give a proof of assertion $(iii)$. Let $0<\theta_0 <\frac{\pi}{4}$.
For $t \geq 0$, the subset $e^{-tB}\omega$ is the cone $\omega$ rotated with angle $-t$, see Figure~\ref{Fig5}. 

\begin{figure}[H]
\begin{tikzpicture}[scale=0.7]
\begin{scope}[xshift=-9cm]
\node at (-3,4) {a) $t=0$};
\filldraw[color=white, fill=black!5] (2.5,2) -- (0,0) -- (2.5,0);
\filldraw[color=white, fill=black!5] (-2.5,-2) -- (0,0) -- (-2.5,0);
\draw (-2.5,-2) -- (2.5,2);
\draw
(1,0) coordinate (a) 
-- (0,0) coordinate (b) 
-- (1,0.8) coordinate (c) 
pic["$\theta_0$", draw=black,  angle eccentricity=1.2, angle radius=1cm]
{angle=a--b--c}; 
\draw[->] (-3,0) -- (3,0); 
\node at (3, -0.2) {\tiny$x$};
\draw[->] (0,-3) -- (0,3); 
\node at (0.2, 3) {\tiny$v$};
\end{scope}
\begin{scope}
\begin{scope}[rotate=-70]
\filldraw[color=white, fill=black!5] (2.5,2) -- (0,0) -- (2.5,0);
\filldraw[color=white, fill=black!5] (-2.5,-2) -- (0,0) -- (-2.5,0);
\draw (-2.5,-2) -- (2.5,2);
\draw (-2.5,0) -- (2.5,0);
\draw
(1,0) coordinate (a) 
-- (0,0) coordinate (b) 
-- (1,0.8) coordinate (c) 
pic["$\theta_0$", draw=black,  angle eccentricity=1.2, angle radius=1cm]
{angle=a--b--c};
\end{scope}
\node at (-3,4) {b) $t>0$};
\draw
(1,0) coordinate (d) 
-- (0,0) coordinate (e) 
-- (1,-2.7) coordinate (f) 
pic["$t$", draw=black,  angle eccentricity=1.2, angle radius=0.5cm]{angle=f--e--d};
\draw[->] (-3,0) -- (3,0); 
\node at (3, -0.2) {\tiny$x$};
\draw[->] (0,-3) -- (0,3); 
\node at (0.2, 3) {\tiny$v$};
\end {scope}
\end{tikzpicture}
\bigskip
\caption{Motion of the cone under rotation of angle $-t$.}\label{Fig5}
\end{figure}

\noindent
For any $t \geq 0$, this subset is never thick in $\rr^2$. The sufficient condition $(i)$ in Proposition~\ref{fg123} therefore does not hold. If $0<T<\pi-\theta_0$, $\mathbb{R}^2 \setminus \cup_{t \in [0,T]}   e^{-tB}\overline{\omega}$ is a non-empty open subset of~$\mathbb{R}^2$.
Let $r>0$. If the point $X=(x,v) \in \rr^2$ is chosen in this set with a sufficiently large norm, then 
$$B_d(X,r)\cap e^{-tB}\omega = \emptyset,$$ 
for all $t \in [0,T]$. The cone $\omega$ does not satisfy the integral thickness condition (\ref{lop2}) on $[0,T]$ associated to the Kolmogorov equation 
with the quadratic external potential (\ref{Kolm_xv}) when  $0<T < \pi-\theta_0$.  
We now consider the case when $T> \pi-\theta_0$ and set $T_*=\inf(T+\theta_0-\pi,\theta_0)>0$. 
There exists a positive constant $\rho>0$ such that for any $X=(x,v) \in \rr^2$ with $|X| >\rho$, there exists a time set $\mathcal{I}(X)$ with Lebesgue measure 
$\lambda(\mathcal{I}(X)) \geq T_*>0$ 
such that 
$$\forall t \in \mathcal{I}(X), \quad \lambda\big( B_d(X,1) \cap e^{-tB}\omega \big) \geq \frac{1}{2} \lambda \big( B_d(0,1) \big).$$
It implies that
$$\forall |X| >\rho, \quad \int_0^T \lambda \big( B_d(X,1) \cap e^{-tB}\omega \big) dt
 \geq \frac{1}{2} \lambda \big( B_d(0,1) \big)T_*>0.$$
 By continuity of the translation in $L^1$, we notice that the function
$$X \in\mathbb{R}^2 \mapsto \int_0^T \lambda \big( B_d(X,1) \cap e^{-tB}\omega \big) dt,$$
is continuous on $\rr^2$.
Furthermore, this function is positive as any point $X \in \mathbb{R}^2$ belongs to $e^{-tB}\omega$ at least for all $t$ in a time set $\tilde{\mathcal{I}}(X) \subset [0,T]$ with Lebesgue measure $\lambda(\tilde{\mathcal{I}}(X))\geq T_*>0$.
By compactness of the set $\overline{B_d(0,\rho)}$, the above function is therefore bounded from below by a positive constant on $\rr^2$,
$$\exists \delta>0, \forall X \in \rr^2, \quad  \int_0^T \lambda\big(B_d(X,1) \cap e^{(t-T)B}\omega\big)dt=\int_0^T \lambda \big( B_d(X,1) \cap e^{-tB}\omega \big) dt \geq \delta.$$
The integral thickness condition (\ref{lop2}) then holds on $[0,T]$ with $r=1$ and $\delta$. This ends the proof of Proposition~\ref{fg123}.

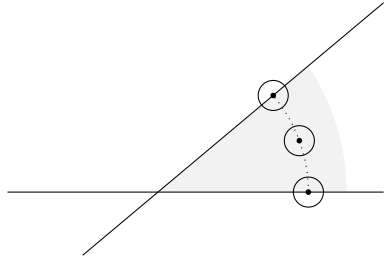
\begin{figure}[H]
\begin{tikzpicture}
\filldraw[color=white, fill=black!5] (2.5,0) -- (0,0) -- (2,1.68);
\filldraw[color=black!5,fill=black!5] (2.5,0) arc (0:34:3);
\draw (-2,0) -- (3,0);
\draw (2,0) circle (0.2cm);
\filldraw (2,0) circle (0.03cm);
\begin{scope}[rotate=20]
\draw (2,0) circle (0.2cm);
\filldraw (2,0) circle (0.03cm);
\end{scope}
\begin{scope}[rotate=40]
\draw (2,0) circle (0.2cm);
\filldraw (2,0) circle (0.03cm);
\end{scope}
\draw[dotted] (2,0) arc (0:40:2);
\draw (-1,-0.84) -- (3,2.52);
\end{tikzpicture}
\bigskip
\caption{Trace of the ball of radius one when the cone rotates.}\label{Fig6}
\end{figure}
\end{proof}

\subsection{Dilation example}

We consider the one-dimensional heat equation 
\begin{equation}\label{hjk10}
(\partial_t-\partial_x^2)f(t,x) = \un_{\omega(t)}(x)u(t,x),
\end{equation}
with the moving control support
$$\omega(t) = \omega  \sqrt{1+2\mu t}, \textrm{ where } \omega=[-1,1] \cup \underset{n \geq 1}{\bigcup} (n^2,n^2+n) \cup (-n^2-n,-n^2),$$
with $\mu>0$. The subset $\omega$ is not thick in $\mathbb{R}$. Indeed, for all $r>0$, there exists $n_0 \geq 1$ such that for all $n \geq n_0$,
$$\Big(n^2+\frac{3n}{2}-r,n^2+\frac{3n}{2}+r\Big) \cap \omega = \emptyset.$$
Equivalently, for any $t>0$, the subset $\omega(t)$ is not thick in $\mathbb{R}$. However, the following proposition shows that the heat equation (\ref{hjk10}) satisfies the integral thickness condition (\ref{lop2}) on $[0,T]$ with the moving control support $(\omega(t))_{t \in [0,T]}$, for any positive time $T>0$.
The null-controllability of the one-dimensional heat equation (\ref{hjk10}) with the moving control subset $(\omega(t))_{t \in [0,T]}$ is an open problem.

\begin{figure}[H]
\centering
 \begin{tikzpicture}[scale=0.8]
\pgfmathsetmacro{\m}{sqrt(3)}
 \begin{scope}
  \node at (-3,1.5) {a) $t=0$};
\draw[->] (-3,0) -- (10,0); 
  \foreach \x in {-3,...,7}{
	\draw[red] (0.15*\x^2 , 0) -- (0.15*\x^2 +0.2*\x, 0);
 } 
 \foreach \x in {-3,-2}{
 \node[red] at (0.15*\x^2 , 0) {$\big)$}; \node[red] at (0.15*\x^2 +0.2*\x, 0) {$\big($};
 \node at (0.15*\x^2-0.1 , -0.5) {\tiny$\x^2$};
 }
\foreach \x in {2,...,7}{
\node[red] at (0.15*\x^2 , 0) {$\big($}; \node[red] at (0.15*\x^2 +0.2*\x, 0) {$\big)$};
\node at (0.15*\x^2+0.1 , -0.5) {\tiny$\x^2$};
}
 \draw[red] (0.15*2,0) -- (-0.15*2, 0);
\node[red] at (0.15*2 , 0) {$\big)$}; \node at (0.15*2, 0.55) {\tiny$2$};
 \node[red] at (-0.15*2 , 0) {$\big($}; \node at (-0.15*2, 0.55) {\tiny$-2$};
 %
 \end{scope}
\begin{scope}[yshift=-3.5cm]
\node at (-1.3,1.5) {b) $t>0, \quad \alpha=\sqrt{1+2\mu t}$};
\draw[->] (-4,0) -- (10,0); 
\foreach \x in {-3,...,5}{
	\draw[red] (0.15*\m *\x^2 , 0) -- (0.15*\m *\x^2 +0.15*\m *\x, 0);
} 
 \foreach \x in {-3,-2}{
 \node[red] at (0.15*\m *\x^2 , 0) {$\big)$}; \node[red] at (0.15*\m *\x^2 +0.15*\m *\x, 0) {$\big($};
 \node at (0.15*\m *\x^2-0.1 , -0.5) {\tiny$\x^2\alpha$};
 }
\foreach \x in {2,...,5}{
\node[red] at (0.15*\m *\x^2 , 0) {$\big($}; \node[red] at (0.15*\m *\x^2 +0.15*\m *\x, 0) {$\big)$};
\node at (0.15*\m *\x^2+0.1 , -0.5) {\tiny$\x^2\alpha$};
}
 \draw[red] (0.15*\m *2,0) -- (-0.15*\m *2, 0);
\node[red] at (0.15*\m *2 , 0) {$\big)$}; \node at (0.15*\m *2-0.1, +0.55) {\tiny$2\alpha$};
 \node[red] at (-0.15*\m *2 , 0) {$\big($}; \node at (-0.15*\m *2+0.1, +0.55) {\tiny$-2\alpha$};
\end{scope}
\end{tikzpicture}
\caption{The control subsets $\omega(t)=\omega\sqrt{1+2\mu t}$ at time $t=0$ and $t>0$.}\label{Fig7}
\end{figure}

\medskip

\begin{Prop}\label{lop7}
For all $T>0$, there exist some positive constants $r,\delta>0$ such that 
\begin{equation}\label{lop6}
\forall x \in \rr, \quad \int_0^{T} \lambda \big((x-r,x+r) \cap \omega(t)   \big) dt \geq \delta >0.
\end{equation} 
\end{Prop}

\medskip

\begin{proof} 
Let  $T>0$. By symmetry, it is sufficient to establish (\ref{lop6}) when $x \geq 0$. We begin by studying the case when $x >\max(5,\sqrt{1+2\mu T})$. We first prove that 
for all $x>\max(5,\sqrt{1+2\mu T})$,
\begin{equation}\label{lop8}
\int_0^{T} \lambda \big((x-1,x+1) \cap \omega(t)   \big) dt \geq
\frac{(x-1)^2}{\mu}  \sum_{k \in \mathcal{I}(x)}  \frac{1}{k^4} \Big( 1- \Big(\frac{x+1}{x-1}\Big)^2 \frac{1}{(1+\frac{1}{k})^2}\Big),
\end{equation}
with
$$\mathcal{I}(x)= \Big\{ k \geq 2:\ \frac{\sqrt{x-1}}{(1+2\mu T)^{\frac{1}{4}}} \leq k \leq \frac{\sqrt{4x+5}-1}{2} \Big\}. $$
Let $x>\max(5,\sqrt{1+2\mu T})$. For any $0<t<T$ and $k \geq 2$, the set inclusion
$$(x-1,x+1) \subset \sqrt{1+2\mu t} (k^2,k^2+k),$$
is equivalent to the two estimates
$$\sqrt{1+2\mu t} k^2 \leq x-1, \qquad x+1 \leq \sqrt{1+2\mu t} ( k^2 + k),$$
or, equivalently 
$$\frac{1}{2\mu} \Big(\Big( \frac{x+1}{k^2+k} \Big)^2 -1\Big)  \leq t \leq \frac{1}{2\mu} \Big(\Big( \frac{x-1}{k^2} \Big)^2 -1\Big).$$
We observe that the set inclusion 
$$\Big(\frac{1}{2\mu} \Big(\Big( \frac{x+1}{k^2+k} \Big)^2 -1\Big), \frac{1}{2\mu} \Big(\Big( \frac{x-1}{k^2} \Big)^2 -1\Big)\Big) \subset (0,T),$$
is equivalent to the condition $k\in \mathcal{I}(x)$.
We finally obtain that
\begin{multline*}
\int_0^{T} \lambda \big((x-1,x+1) \cap \omega(t)   \big) dt \geq 
\int_0^{T} \lambda \Big((x-1,x+1) \cap \bigcup_{k \in \mathcal{I}(x)}\sqrt{1+2\mu t}(k^2,k^2+k)\Big) dt\\
=\sum_{k \in \mathcal{I}(x)}\int_0^{T} \lambda \big((x-1,x+1) \cap \sqrt{1+2\mu t}(k^2,k^2+k)\big) dt 
\end{multline*}
and 
$$\int_0^{T} \lambda \big((x-1,x+1) \cap \omega(t)   \big) dt \geq 
\sum\limits_{k \in \mathcal{I}(x)} \frac{1}{\mu} \Big(\Big(\frac{x-1}{k^2} \Big)^2 - \Big( \frac{x+1}{k^2+k} \Big)^2 \Big).$$
It proves the estimate (\ref{lop8}).
Next, we prove that there exist some positive constants $A>1$ and $C>0$ such that for all $x>A$,
\begin{equation}\label{lop9}
\int_0^{T} \lambda \big((x-1,x+1) \cap \omega(t)\big) dt \geq \frac{C}{4\mu}.
\end{equation}
We first notice that there exists $\eta>0$ such that for all $0<x<\eta$,
$$\frac{1}{(1+x)^2} \leq 1 - \frac{x}{2}.$$
There exists a positive constant $A'>\max(5,\sqrt{1+2\mu T})$ such that for all $x>A'$, 
$$\frac{\sqrt{x-1}}{(1+2\mu T)^{\frac{1}{4}}}>\frac{1}{\eta}.$$
We obtain that for all $x>A'$ and $k \in \mathcal{I}(x)$,
$$\frac{1}{k} \leq  \frac{(1+2\mu T)^{\frac{1}{4}}}{\sqrt{x-1}} < \eta.$$
It implies that for all $x>A'$ and $k \in \mathcal{I}(x)$,
\begin{equation}\label{sf100}
\frac{1}{(1+\frac{1}{k})^2} \leq 1 - \frac{1}{2k}.
\end{equation}
We can find some positive constants $C>0$, $C'>0$ and $A''>A'$ such that for all $x>A''$,
$$\sum\limits_{k \in \mathcal{I}(x)} \frac{1}{k^4} 
\leq \frac{1+2\mu T}{(x-1)^2}\# \mathcal{I}(x)
\leq \frac{1+2\mu T}{(x-1)^2} 
\Big(  \frac{\sqrt{4x+5}-1}{2}- \frac{\sqrt{x-1}}{(1+2\mu T)^{\frac{1}{4}}}   +1\Big)
\leq \frac{C'}{x^{3/2}}
$$ 
and
$$\sum\limits_{k \in \mathcal{I}(x)} \frac{1}{k^5}\geq \frac{2^5}{(\sqrt{4x+5}-1)^5} \# \mathcal{I}(x) \geq \frac{C}{x^2}.$$
According to (\ref{lop8}) and (\ref{sf100}), we obtain that for all $x>A''$,
\begin{align*}
\int_0^{T} \lambda \big((x-1,x+1) \cap \omega(t)\big) dt 
&  \geq 
\frac{(x-1)^2}{\mu}  \sum_{k \in \mathcal{I}(x)}  \frac{1}{k^4} \Big( 1- \Big(\frac{x+1}{x-1}\Big)^2  + \Big(\frac{x+1}{x-1}\Big)^2 \frac{1}{2k}\Big)
\\
&  \geq
- \frac{4x}{\mu} \Big( \sum\limits_{k \in \mathcal{I}(x)} \frac{1}{k^4}  \Big)
+ \frac{(x+1)^2}{2\mu}  \Big( \sum\limits_{k \in \mathcal{I}(x)} \frac{1}{k^5}\Big) 
\\
&  \geq
-\frac{4C'}{\mu\sqrt{x}} + \frac{C}{2\mu}.
\end{align*}
Finally, there exists a positive constant $A>A''$ such that for all $x>A$, 
$$\frac{4C'}{\mu\sqrt{x}} \leq \frac{C}{4\mu},$$
which establishes the estimate (\ref{lop9}). We now conclude with a compactness argument when $|x| \leq A$.
For any $x \in [-A,A]$, we notice that for all $0<t<T$,
$$(-1,1) \subset (x-2A,x+2A) \cap \omega(t).$$
It implies that
$$\forall |x| \leq A, \quad \int_0^{T} \lambda \big((x-2A,x+2A) \cap \omega(t)   \big) dt \geq 2T>0.$$
It gives the conclusion with $r=2A$ and $\delta=\min\{\frac{C}{4\mu}, 2T \}$, as when $x>A$, 
$$(x-1,x+1) \subset (x-2A,x+2A),$$
since $A>1$. This ends the proof of Proposition~\ref{lop7}.
\end{proof}

\section{An adapted Lebeau-Robbiano method with moving control subsets}\label{Sec:LR}

\medskip

\begin{Def}
Let $(X,\|.\|_X)$ be a Banach space.
A two-parameters family of bounded linear operators $\mathscr{U}(t,s)$, $0 \leq s \leq t \leq T$
is an evolution system of contractions of $X$ if the following three conditions are satisfied:
\begin{enumerate}
\item[$(i)$] $\forall 0 \leq s \leq T,  \ \mathscr{U}(s,s)=\emph{\textrm{Id}}_X$; $\forall 0 \leq s \leq r \leq t \leq T, \ \mathscr{U}(t,r)\mathscr{U}(r,s)=\mathscr{U}(t,s)$
\item[$(ii)$] $(t,s) \mapsto \mathscr{U}(t,s) \in \mathcal{L}_c(X)$ is strongly continuous for all $0 \leq s \leq t \leq T$
\item[$(iii)$] $\forall 0 \leq s \leq t \leq T, \ \|\mathscr{U}(t,s)\|_{\mathcal{L}_c(X)} \leq 1$
\end{enumerate}
where $\emph{\textrm{Id}}_X$ denotes the identity operator on $X$ and $\mathcal{L}_c(X)$ stands for the set of bounded linear operators on $X$.
\end{Def}

\medskip

The following result extends the abstract observability result established in~\cite{KB_KPS_JEP} (Theorem~2.1) to evolution systems and time-varying control subsets. It also allows some controlled blow-ups for small times in the dissipation estimates that is not covered by the result of~\cite{KB_KPS_JEP} (Theorem~2.1), and that is absolutely needed for various study cases and in particular for the proof of Theorem~\ref{Main_result_quad} in the present work. Despite these improvements, the following proof essentially follows the very same lines as the one given in~\cite{KB_KPS_JEP} (Theorem~2.1).

\medskip

\begin{thm} \label{Meta_thm_AdaptedLRmethod}
Let $\Omega$ be an open subset of $\mathbb{R}^d$; $T>0$; $(\omega(t))_{t \in [0,T]}$ be a moving control support in $\Omega$; $(\pi_k)_{k \geq 1}$ be a family of orthogonal projections defined on $L^2(\Omega)$; $\mathscr{U}(t,s)$, $0 \leq s \leq t \leq T$ be an evolution system of contractions on $L^2(\Omega)$; $\delta_1, m_2 \geq 0$; $\delta_2, c_1, c_1', c_2, c_2', a, b, m_1 >0$ be positive constants with $a<b$ and $0 \leq \delta_1 < \delta_2 \leq T$; $E$ be a measurable subset of $[\delta_1,\delta_2]$ with positive Lebesgue measure.
If the following uniform spectral estimates with respect to $t \in E$,
\begin{equation} \label{Meta_thm_IS}
\forall  t \in E, \forall g \in L^2(\Omega), \forall k \geq 1, \quad \|\pi_k g \|_{L^2(\Omega)} \leq c_1' e^{c_1 k^a} \|\pi_k g \|_{L^2(\omega(t))},
\end{equation}
and the following dissipation estimates with controlled blow-up
\begin{multline} \label{Meta_thm_dissip}
\forall g \in L^2(\Omega), \forall k \geq 1, \forall \delta_1 \leq s < t \leq \delta_2, \\
\|(1-\pi_k)\mathscr{U}(t,s) g \|_{L^2(\Omega)} \leq \frac{1}{c_2'(t-s)^{m_2}} e^{-c_2 (t-s)^{m_1} k^b} \|g\|_{L^2(\Omega)},
\end{multline}
hold, then there exists a positive constant $C>0$ such that the following observability estimate holds
\begin{multline} \label{meta_thm_IO}
\text{ for a.e. } t^* \in E, \forall R>0, \exists t_0 \in (t^*,t^*+R) \cap E, \forall g \in L^2(\Omega),
\\
\|\mathscr{U}(t_0,t^*) g \|_{L^2(\Omega)}^2 \leq \exp\Big(\frac{C}{(t_0-t^*)^{\frac{am_1}{b-a}}}\Big) \int\limits_{(t^*,t_0) \cap E} \|\mathscr{U}(t,t^*) g \|_{L^2(\omega(t))}^2  dt.
\end{multline}
In particular, there exists a positive constant $C_1>0$ such that
\begin{equation}\label{ty1}
\forall g \in L^2(\Omega), \quad \|\mathscr{U}(T,0) g \|_{L^2(\Omega)}^2 \leq C_1 \int\limits_{E} \|\mathscr{U}(t,0) g \|_{L^2(\omega(t))}^2  dt.
\end{equation}
More specifically, let $C>0$ be a positive constant verifying
\begin{equation} \label{meta_thm_cost_estim}
\displaystyle C >  c_1^{\frac{b}{b-a}}  c_2^{-\frac{a}{b-a}} \frac{2^\alpha \beta^\beta}{\gamma^\gamma (\beta-\gamma)^{\beta-\gamma}},
\end{equation}
with
$$\gamma= \frac{a m_1}{b-a}>0, \quad
\beta= 2\gamma + \frac{a b m_1}{(b-a)^2}>0, \quad
\alpha= \gamma + 1+ a + \frac{a^2}{b-a}>0.$$
Then, there exists $T_0>0$ such that, if $E=[\delta_1,\delta_2]=[0,T]$ with $0<T \leq T_0$, then
\begin{equation}\label{ty1_BIS}
\|\mathscr{U}(T,0) g \|_{L^2(\Omega)}^2 \leq \exp\left(\frac{C}{T^{\frac{am_1}{b-a}}}\right) \int_0^T \|\mathscr{U}(t,0) g \|_{L^2(\omega(t))}^2  dt.
\end{equation}
\end{thm}

\medskip

The proof of Theorem~\ref{Meta_thm_AdaptedLRmethod} relies on a telescopic series argument
due to~\cite{miller2010} (see also~\cite{miller_SMF}), and already used in~\cite{apraiz,phung,YubiaoZhang}.
The proof of the observability estimate (\ref{meta_thm_IO}) uses the following result which is already stated in~\cite{Gengsheng_Wang} (Lemma~2.3) and briefly proved in~\cite{JLL} (pp. 256-257). For the sake of completeness of the present work, a proof is given: 

\medskip

\begin{Prop} \label{Prop:theorie_mesure}
Let $0\leq \delta_1 < \delta_2<+\infty$ and $E$ be a measurable subset of $[\delta_1,\delta_2]$ with positive Lebesgue measure.
There exist some positive constants $0<\rho < 1$ and $C_0, C_0'>0$ such that 
for almost every $t^* \in E$ and for all $R>0$, there exists a decreasing sequence $(t_j)_{j \geq 0}$ 
of $(t^*,t^*+R) \cap E$ such that 
$$t_j \underset{j \rightarrow +\infty}{\longrightarrow} t^*,$$
\begin{equation} \label{TM:1}
\forall j \geq 1,\quad \lambda\Big(  \Big( t_{j},\frac{t_{j}+t_{j-1}}{2}  \Big) \cap E \Big) \geq \rho \Big(\frac{t_{j-1} - t_{j}}{2}\Big)>0,
\end{equation}
\begin{equation} \label{TM:2}
\forall j \geq 1, \quad t_{j}-t_{j+1} \geq C_0 (t_{j-1}-t_j)>0
\end{equation}
and
\begin{equation} \label{TM:3}
\quad t_0 - t_1 \geq C_0' (t_0 - t^*) > 0.
\end{equation}
\end{Prop}

\medskip

\begin{proof} 
The above result is proved here with the constants $\rho=\frac{3}{4}$, $C_0=\frac{1}{12}$ and $C_0'=\frac{1}{4}$. The proof can actually be performed with any $0<\rho<1$ and related constants $C_0(\rho)>0$ and $C'_0(\rho)>0$. For a given measurable subset $A$ of $\R$,
we recall that a Lebesgue, or density, point of $A$ is a point $x \in A$ satisfying
$$\frac{1}{2r} \lambda\big( [x-r,x+r] \cap A \big) = \frac{1}{2r} \int_{x-r}^{x+r} \un_A(y) dy  \underset{\substack{r \rightarrow 0\\ r>0}}{\longrightarrow} 1.$$
This property readily implies in particular that
\begin{equation}\label{bn1}
\frac{1}{r} \lambda\big( [x,x+r] \cap A \big) = \frac{1}{r} \int_{x}^{x+r} \un_A(y) dy  \underset{\substack{r \rightarrow 0\\ r>0}}{\longrightarrow} 1.
\end{equation}
By Lebesgue theorem, almost every point of $A$ is a density point of $A$.

\medskip

\noindent
\textit{Step 1}. We begin by constructing a particular subset $\widetilde{E}$ of $E$ satisfying $\lambda(\widetilde{E})=\lambda(E)$. It will be then sufficient to establish Proposition~\ref{Prop:theorie_mesure} for any $t^* \in \widetilde{E}$.
For any $m \geq 1$, we introduce the set 
$$E_m = \Big\{ \sigma \in E: \ \forall 0<r<\frac{1}{m}, \ \lambda\big( [\sigma,\sigma+r] \cap E \big) \geq \frac{3r}{4}\Big\} \subset E.$$
Let $D_m$, respectively $D$, be the set of Lebesgue points of $E_m$, respectively of $E$.
The sequence of subsets $(E_m)_{m \geq 1}$ is non-decreasing for the inclusion 
$$\forall m \geq 1, \quad E_m \subset E_{m+1}.$$
The sequence of subsets $(D_m)_{m \geq 1}$ is therefore also non-decreasing for the inclusion 
$$\forall m \geq 1, \quad D_m \subset D_{m+1}.$$
According to (\ref{bn1}), any Lebesgue point of $E$ belongs to a subset $E_m$ for some $m \geq 1$ sufficiently large 
$$D \subset \underset{m \geq 1}{\bigcup} E_m.$$
By Lebesgue theorem, we have $\lambda(D)=\lambda(E)$ and $\lambda(D_m)=\lambda(E_m)$ for all $m \geq 1$.
It follows that
$$\lambda(E)=\lambda(D) \leq 
\lambda \Big(\underset{m \geq 1}{\bigcup}  E_m\Big) 
= \lim_{m \rightarrow +\infty} \lambda(E_m)
= \lim_{m \rightarrow +\infty}  \lambda(D_m) = \lambda \Big( \underset{m \geq 1}{\bigcup} D_m \Big).$$
As 
$\underset{m \geq 1}{\bigcup} D_m \subset E,$
we obtain that
$$\lambda(E)=\lambda \Big( \underset{m \geq 1}{\bigcup} D_m \Big).$$
We set 
$\widetilde{E}= \underset{m \geq 1}{\bigcup} D_m.$
Let $t^* \in \widetilde{E}$ and $R>0$. There exists $m^*\geq 1$ such that $t^* \in D_{m^*}$. It follows from (\ref{bn1}) that
$$\frac{1}{r} \lambda\big( [t^*,t^*+r] \cap E_{m^*} \big)   \underset{\substack{r \rightarrow 0\\ r>0}}{\longrightarrow} 1.$$
There exists $0<r_0<\min\{ R,\frac{1}{m^*}\}$ such that
\begin{equation} \label{def:r0}
\forall 0 < r \leq r_0, \quad \lambda\big( [t^*,t^*+r] \cap E_{m^*} \big) \geq \frac{3r}{4}.
\end{equation}

\medskip

\noindent
\textit{Step 2}. Next, we construct by induction a decreasing sequence $(t_j)_{j \geq 0}$ of $E_{m^*} \cap (t^*,t^*+r_0)$ such that
\begin{equation}\label{lop11}
t_{j+1} \in  \Big[ t^*+\frac{t_j-t^*}{4} ,   t^*+\frac{3(t_j-t^*)}{4} \Big] \cap E_{m^*}  \subset  (t^*,t^*+r_0) \cap E_{m^*}.
\end{equation}
According to (\ref{def:r0}), the set $(t^*,t^*+r_0) \cap E_{m^*}$ is not empty. We choose an arbitrarily point $t_0 \in (t^*,t^*+r_0) \cap E_{m^*}$.
Let us now assume that the points $(t_l)_{0 \leq l \leq j}$, with $j \geq 0$, are already constructed and satisfy (\ref{lop11}) for all $0 \leq l \leq j-1$. We aim at finding $t_{j+1} \in (t^*,t^*+r_0)$ satisfying (\ref{lop11}).
Since $0<t_j-t^*< r_0$, we deduce from (\ref{def:r0}) that
$$\lambda\big( [t^*,t_j] \cap E_{m^*} \big) \geq \frac{3(t_j-t^*)}{4}.$$
The set $[ t^*+\frac{t_j-t^*}{4} ,   t^*+\frac{3(t_j-t^*)}{4}] \cap E_{m^*}$ is therefore not empty, since
otherwise the above measure would necessarily be less or equal than $\frac{t_j-t^*}{2}$. It is then sufficient to choose the point $t_{j+1}$ arbitrarily in this set. The sequence $(t_j)_{j \geq 0}$ is decreasing by construction of $(t^*,t^*+R) \cap E$.

\medskip

\noindent
\textit{Step 3}. Let $j \geq 1$. Since by construction $t_{j} \in E_{m^*}$ and $0<r:=\frac{t_{j-1} - t_{j}}{2} < \frac{r_0}{2} < \frac{1}{m^*}$, we deduce from the definition of the set $E_{m^*}$ that 
$$\lambda\big( (t_{j},t_{j}+r) \cap E \big) \geq \frac{3}{4} r,$$
which proves (\ref{TM:1}) with $\rho=\frac{3}{4}$.
Furthermore, it follows from (\ref{lop11}) that for all $j \geq 0$,   
$$\frac{t_j-t^*}{4} \leq t_{j+1}-t^* \leq \frac{3(t_j-t^*)}{4},$$
which implies that
$$t_0-t_1 \geq t_0-\Big(\frac{3}{4}t_0+\frac{1}{4}t^*\Big)=\frac{1}{4}(t_0-t^*).$$
This establishes (\ref{TM:3}) with $C_0'=\frac{1}{4}$.
We finally obtain by using anew (\ref{lop11}) that for all $j \geq 1$, 
\begin{multline*}
t_{j}-t_{j+1} \geq t_j-\Big(\frac{3t_j}{4}+\frac{t^*}{4}\Big)=\frac{1}{4}(t_j-t^*) \geq \frac{1}{16}(t_{j-1}-t^*) \\ \geq \frac{1}{16}\Big(t_{j-1}-\frac{4}{3}\Big(t_j-\frac{t_{j-1}}{4}\Big)\Big)=
\frac{1}{12}(t_{j-1}-t_j)>0.
\end{multline*}
This proves the estimates (\ref{TM:2}) with $C_0=\frac{1}{12}$.
\end{proof}

\medskip

We can now prove Theorem~\ref{Meta_thm_AdaptedLRmethod}.

\begin{proof}
For simplicity, the notation $\|\cdot\|$ refers in this proof to the norm $\|\cdot\|_{L^2(\Omega)}$. 
We aim at establishing (\ref{meta_thm_IO}) with the parameters $t^*$ and $t_0$ given by Proposition~\ref{Prop:theorie_mesure}.
We first observe that the estimate (\ref{ty1})  is a direct consequence of the contraction property of the evolution system and  (\ref{meta_thm_IO}).
Indeed, by applying  (\ref{meta_thm_IO}) to some allowed values of $t^* \in E$ and $t_0 \in (t^*,t^*+1) \cap E$, we obtain from the contraction evolution system properties that for all $g \in L^2(\Omega)$,
\begin{align*}
\|\mathscr{U}(T,0) g \|^2   = & \ \|\mathscr{U}(T,t_0)\mathscr{U}(t_0,t^*)\mathscr{U}(t^*,0)g\|^2  
 \leq \|\mathscr{U}(t_0,t^*)\mathscr{U}(t^*,0)g\|^2  \\
 \leq & \ \exp\Big(\frac{C}{(t_0-t^*)^{\frac{am_1}{b-a}}}\Big) \int_{(t^*,t_0)\cap E} \|\mathscr{U}(t,t^*)\mathscr{U}(t^*,0) g \|_{L^2(\omega(t))}^2 dt \\
 \leq & \ C_1 \int_{E} \|\mathscr{U}(t,0) g \|_{L^2(\omega(t))}^2 dt.
\end{align*}
with 
$$C_1=\exp\Big(\frac{C}{(t_0-t^*)^{\frac{am_1}{b-a}}}\Big)>1.$$
Let us now prove the estimate (\ref{meta_thm_IO}). 
Let $0<\rho<1$, $C_0, C_0'>0$ be as in Proposition~\ref{Prop:theorie_mesure}.
Let $0<\eps<2$. 
We consider the positive constants
\begin{equation} \label{def:gamma}
\gamma =\Big(  \frac{ (2+\eps) c_1 2^a }{ (2-\eps) c_2 C_0^{\frac{b m_1}{b-a}} } \Big)^{\frac{1}{b-a}}>0, \qquad 
M = (2+\eps) c_1(2\gamma)^a>0.
\end{equation}
By the definition (\ref{def:gamma}) of the constant $\gamma$, we observe that
\begin{equation} \label{def:M}
M=(2+\eps)c_1(2\gamma)^a = (2-\eps) c_2 \gamma^b  C_0^{\frac{b m_1}{b-a}}.
\end{equation}
Let $R>0$. There exists $0<\tilde{R} <R$ such that for all $0<\tau<\tilde{R}$,
\begin{equation} \label{def:R}
\frac{\gamma}{\tau^{\frac{m_1}{b-a}}} > 1, 
\
\frac{ \rho \tau }{2 c_1'^2}  \geq \exp\Big( - \eps c_1 \frac{(2\gamma)^a }{\tau^{\frac{a m_1}{b-a}} }   \Big),
\  
\frac{2}{C_0 c_2'^2 \tau^{2m_2-1}} \leq \exp\Big(  \eps c_2 \frac{\gamma^b C_0^{\frac{b m_1}{b-a}}}{\tau^{\frac{a m_1}{b-a}}}   \Big).
\end{equation}
Then, for all $0<\tau< \tilde{R}$, there exists an integer $k(\tau) \geq 1$ verifying
\begin{equation}\label{def:k}
1<\frac{\gamma}{\tau^{\frac{m_1}{b-a}}} \leq k(\tau) \leq  \frac{2 \gamma}{\tau^{\frac{m_1}{b-a}}},
\end{equation}
since according to (\ref{def:R}), the interval $(\gamma \tau^{-\frac{m_1}{b-a}}, 2 \gamma \tau^{-\frac{m_1}{b-a}})$  is of length greater than $1$, and is contained in $(1,+\infty)$.
Let $t^*$ and $(t_j)_{j \geq 0}$ be as in Proposition \ref{Prop:theorie_mesure} applied with the constant $\tilde{R}>0$ defined in (\ref{def:R}). We define
$$\forall j \geq 1, \quad \tau_j=\frac{t_{j-1}-t_{j}}{2}>0.$$
We observe from (\ref{TM:1}) and (\ref{TM:2}) that 
\begin{equation} \label{tau:R}
\forall j \geq 1, \quad 0<\tau_j < \tilde{R},
\end{equation}
\begin{equation} \label{tau:C0}
\forall j \geq 1, \quad   \tau_j \leq \frac{1}{C_0} \tau_{j+1},
\end{equation}
\begin{equation} \label{tau:rho}
\forall j \geq 1, \quad \lambda\left( (t_{j},t_{j}+\tau_j) \cap E \right) \geq \rho \tau_j.
\end{equation}
According to (\ref{def:k}) and (\ref{tau:R}), we can define $k_j = k(\tau_j) \geq 1$ such that
\begin{equation} \label{relation_k}
\forall j \geq 1, \quad \frac{\gamma}{\tau_j^{\frac{m_1}{b-a}}} \leq k_j \leq  \frac{2 \gamma}{\tau_j^{\frac{m_1}{b-a}}}.
\end{equation}

\medskip

\noindent
\emph{Step 1}. We begin by establishing the following estimate: $\forall j \geq 1$, $\forall g \in L^2(\Omega)$,
\begin{multline} \label{meta_thm_telescopic_formula}
f(\tau_j) \|\mathscr{U}(t_{j}+\tau_{j},t^*) g \|^2  - f(\tau_{j+1}) \| \mathscr{U}(t_{j+1}+\tau_{j+1},t^*) g \|^2 \\
\leq \int_{( t_{j} , t_{j}+\tau_{j} ) \cap E } \|\mathscr{U}(t,t^*) g \|_{L^2(\omega(t))}^2 dt,
\end{multline}
where
\begin{equation}\label{meta_thm_def_f}
f(s)=\exp\Big(-\frac{M}{s^{\frac{am_1}{b-a}}}\Big), \quad s>0.
\end{equation}
By using successively (\ref{tau:rho}), the contraction property of the evolution system, the Pythagorean identity, the spectral estimates (\ref{Meta_thm_IS}),  the triangular inequality and $\|\cdot\|_{L^2(\omega(t))} \leq \|\cdot\|$, we obtain that for all $j \geq 1$,
\begin{align*}
  & \ \rho \tau_j \frac{e^{- 2 c_1 k_j^a}}{2 c_1'^2}\|\mathscr{U}(t_j+\tau_j,t^*) g \|^2   \\   
   \leq &\  \frac{e^{- 2 c_1 k_j^a}}{2 c_1'^2} \int_{(t_j,t_j+\tau_j) \cap E} \|\mathscr{U}(t,t^*) g \|^2 dt \\ 
    \leq &\ \frac{e^{- 2 c_1 k_j^a}}{2 c_1'^2}  \int_{(t_j,t_j+\tau_j) \cap E}  \left( \| \pi_{k_j} \mathscr{U}(t,t^*) g \|^2 + \| (1-\pi_{k_j})\mathscr{U}(t,t^*)g \|^2 \right) dt \\ 
   \leq &\  \int_{(t_j,t_j+\tau_j) \cap E}  \left( \frac{1}{2} \| \pi_{k_j} \mathscr{U}(t,t^*) g \|_{L^2(\omega(t))}^2  +  \| (1-\pi_{k_j})\mathscr{U}(t,t^*) g \|^2  \right) dt
\end{align*}
that is
\begin{align*}
  & \ \rho \tau_j \frac{e^{- 2 c_1 k_j^a}}{2 c_1'^2}\|\mathscr{U}(t_j+\tau_j,t^*) g \|^2   \\ 
 \leq &\  \int_{(t_j,t_j+\tau_j) \cap E}  \Big(  \| \mathscr{U}(t,t^*) g \|_{L^2(\omega(t))}^2  + \| (1-\pi_{k_j} )\mathscr{U}(t,t^*) g \|_{L^2(\omega(t))}^2\\
 & \qquad \qquad  \qquad \qquad \qquad \qquad \qquad \qquad +  \| (1-\pi_{k_j})\mathscr{U}(t,t^*) g \|^2  \Big) dt   \\  
   \leq &\  \int_{(t_j,t_j+\tau_j) \cap E}  \left( \|\mathscr{U}(t,t^*) g \|_{L^2(\omega(t))}^2  + 2 \| (1-\pi_{k_j})\mathscr{U}(t,t^*) g \|^2  \right) dt. 
\end{align*}
By using successively the dissipation estimates (\ref{Meta_thm_dissip}) and the fact that
$$t \in (t_j,t_j+\tau_j) \quad \Longrightarrow \quad t-(t_{j+1}+\tau_{j+1}) \geq \tau_{j+1}>0,$$
for any $j \geq 1$, we obtain that for all $j \geq 1$,
\begin{align}
  & \  \int_{(t_j,t_j+\tau_j) \cap E}  2 \| (1-\pi_{k_j})\mathscr{U}(t,t^*)  g \|^2   dt \notag \\ \notag  
    \leq &\ \int_{(t_j,t_j+\tau_j)}  2 \| (1-\pi_{k_j})\mathscr{U}(t, t_{j+1}+\tau_{j+1})      \mathscr{U}(t_{j+1}+\tau_{j+1},t^*)  g \|^2   dt \notag \\ \notag  
   \leq &\ \int_{(t_j,t_j+\tau_j) } \frac{2e^{-2 c_2  [t-(t_{j+1}+\tau_{j+1})]^{m_1}  k_j^b }}{c_2'^2 [t-(t_{j+1}+\tau_{j+1})]^{2m_2}}     \| \mathscr{U}(t_{j+1}+\tau_{j+1},t^*) g \|^2 dt \\ \notag
   \leq &\  \frac{2 \tau_j}{c_2'^2 \tau_{j+1}^{2m_2}}  e^{-2c_2  \tau_{j+1}^{m_1}  k_j^b } \| \mathscr{U}(t_{j+1}+\tau_{j+1},t^*) g \|^2.  
\end{align}
We deduce from the previous two estimates that for all $j \geq 1$,
\begin{multline} \label{Telesc_1}
\frac{\rho \tau_j}{2 c_1'^2}  e^{- 2 c_1 k_j^a} \|\mathscr{U}(t_j+\tau_j,t^*) g \|^2 \leq 
  \int_{(t_j,t_j+\tau_j) \cap E}   \|\mathscr{U}(t,t^*) g \|_{L^2(\omega(t))}^2  dt   \\ 
+   \frac{2 \tau_j}{c_2'^2 \tau_{j+1}^{2m_2}}  e^{-2 c_2  \tau_{j+1}^{m_1}  k_j^b } \| \mathscr{U}(t_{j+1}+\tau_{j+1},t^*) g \|^2.
\end{multline}
By first using (\ref{def:R}), (\ref{tau:R}) and (\ref{relation_k}), and then (\ref{def:gamma}) and (\ref{meta_thm_def_f}), we obtain that for all $j \geq 1$,
\begin{multline*}
\frac{\rho \tau_j}{2 c_1'^2}  e^{- 2 c_1 k_j^a}
\geq 
\frac{\rho \tau_j}{2 c_1'^2} \exp\Big(-2 c_1 \frac{(2 \gamma)^a}{\tau_j^{\frac{a m_1}{b-a}}} \Big)
\geq\\ 
\exp\Big(-(2+\eps) c_1 \frac{(2 \gamma)^a}{\tau_j^{\frac{a m_1}{b-a}}} \Big) =  \exp\Big(- \frac{M}{\tau_j^{\frac{a m_1}{b-a}}} \Big) = f(\tau_j).
\end{multline*}
By using successively (\ref{relation_k}), (\ref{tau:C0}),  (\ref{tau:R}), (\ref{def:R}), (\ref{def:M})  and (\ref{meta_thm_def_f}), we obtain that for all $j \geq 1$,
\begin{align*}
 \frac{2 \tau_j}{c_2'^2 \tau_{j+1}^{2m_2}}  e^{ -2 c_2  \tau_{j+1}^{m_1}  k_j^b }  \leq & \ \frac{2 \tau_j}{c_2'^2 \tau_{j+1}^{2m_2}}   \exp\Big( -2 c_2  \tau_{j+1}^{m_1} \frac{\gamma^b}{\tau_j^{\frac{b m_1 }{b-a}}}\Big) \\
 \leq & \ \frac{2 }{C_0 c_2'^2 \tau_{j+1}^{2m_2-1}}  \exp\Big( -2 c_2  C_0^{\frac{b m_1}{b-a}} \frac{\gamma^b}{\tau_{j+1}^{\frac{a m_1 }{b-a}}}\Big) \\
 \leq  &\  \exp\Big( - (2-\eps) c_2  C_0^{\frac{b m_1}{b-a}} \frac{\gamma^b}{\tau_{j+1}^{\frac{a m_1 }{b-a}}}\Big)=\exp\Big(- \frac{M}{\tau_{j+1}^{\frac{a m_1}{b-a}}} \Big) = f(\tau_{j+1}).
\end{align*}
According to the two above estimates, we deduce (\ref{meta_thm_telescopic_formula}) from (\ref{Telesc_1}).

\medskip

\noindent
\textit{Step 2}. We can now derive the observability estimate (\ref{meta_thm_IO}).
Summing up the estimates  (\ref{meta_thm_telescopic_formula})  for all $j \geq 1$ provides that
\begin{equation}\label{kl1}
f(\tau_1) \|\mathscr{U}(t_1+\tau_1,t^*)g\|^2  \leq \int_{(t^*,t_0) \cap E} \|\mathscr{U}(t,t^*) g \|_{L^2(\omega(t))}^2 dt,
\end{equation}
since by the contractivity property of the evolution system and (\ref{meta_thm_def_f}),
$$f(\tau_j) \|\mathscr{U}(t_j+\tau_j,t^*) g \|^2 \leq \exp\Big(- \frac{M}{\tau_j^{\frac{am_1}{b-a}}}\Big) \|g\|^2 \underset{j \rightarrow +\infty}{\longrightarrow} 0,$$
as $\tau_j \to 0$ when $j \to+\infty$, because $t_j \to t^*$ when $j \to+\infty$; 
and since the intervals $(t_j,t_j+\tau_j)$ are disjoint and included in $(t^*,t_0)$.
According to (\ref{TM:3}), the estimate
$$\tau_1 = \frac{t_0-t_1}{2} \geq\frac{C_0'}{2} (t_0-t^*)>0,$$
implies that
$$f(\tau_1)=\exp\Big( - \frac{M}{\tau_1^{\frac{a m_1}{b-a}}}\Big) \geq \exp\Big( - \frac{C}{(t_0-t^*)^{\frac{a m_1}{b-a}}} \Big),$$
with  $C=M \big( \frac{2}{C_0'} \big)^{\frac{a m_1}{b-a}}$.
We finally obtain from the contractivity property of the evolution system that
\begin{multline}\label{hel1}
\|\mathscr{U}(t_0,t^*)g\|^2  \leq \|\mathscr{U}(t_1+\tau_1,t^*)g\|^2 \\
\leq \exp\Big( \frac{C}{(t_0-t^*)^{\frac{a m_1}{b-a}}} \Big) \int_{(t^*,t_0) \cap E} \|\mathscr{U}(t,t^*) g \|_{L^2(\omega(t))}^2 dt.
\end{multline}
This ends the proof of (\ref{meta_thm_IO}).

\medskip

\noindent
\emph{Step 3: Estimate of the observability cost when $E=[0,T]$.}
We consider the case when  $E=[\delta_1,\delta_2]=[0,T]$. Let $\mu>0$. We define the sequence $(t_j)_{j \geq 0}$ by
$$\forall j \geq 0, \quad t_j=\frac{T}{2^{(j+1)\mu}} \in (0,T).$$
We notice that on one hand 
$$t_j \underset{j \rightarrow +\infty}{\longrightarrow} 0$$ 
and that on the other hand, the conclusions of Proposition~\ref{Prop:theorie_mesure} hold with $t^*=0$ and $R=T$, since assertion (\ref{TM:1}) holds with $\rho=1$, assertion (\ref{TM:2}) holds with $C_0=\frac{1}{2^\mu}$ and assertion (\ref{TM:3}) holds with $C_0'=1-\frac{1}{2^\mu}$. We can therefore resume the previous proof and deduce from (\ref{hel1}) the following observability estimate holds for all $g \in L^2(\Omega)$,
\begin{align}\label{hel2}
& \ \|\mathscr{U}(T,0)g\|^2  \leq \Big\|\mathscr{U}\Big(\frac{T}{2^{\mu}},0\Big)g\Big\|^2 \\ \notag
\leq & \ \exp\Big( \frac{2^{\mu\frac{a m_1}{b-a}}C}{T^{\frac{a m_1}{b-a}}} \Big) \int_{(0,\frac{T}{2^{\mu}})} \|\mathscr{U}(t,0) g \|_{L^2(\omega(t))}^2 dt\\ \notag
\leq & \ \exp\Big( \frac{\tilde{C}}{T^{\frac{a m_1}{b-a}}} \Big) \int_0^T \|\mathscr{U}(t,0) g \|_{L^2(\omega(t))}^2 dt.
\end{align}
by the contractivity property of the evolution system with the constants
$$\gamma =\Big(\frac{ (2+\eps) c_1  }{ (2-\eps) c_2 }  2^{a + \mu \frac{b m_1}{b-a}}\Big)^{\frac{1}{b-a}}, \quad 
M = (2+\eps) c_1 (2\gamma)^a, \quad \tilde{C}= M \Big( \frac{2^{2\mu+1}}{2^\mu -1} \Big)^{ \frac{a m_1}{b-a} }.$$
By choosing the parameter $0<\eps <2$ small enough, the positive constant $\tilde{C}>0$ appearing in the observability estimate (\ref{hel2}) may be chosen arbitrarily as
\begin{equation} \label{cost:mu}
\tilde{C}> \Big( \frac{2^{2\mu+1}}{2^\mu -1} \Big)^{ \frac{a m_1}{b-a} }   c_1 2^{1+a}  \Big(  \frac{ c_1  }{ c_2 }  2^{a + \mu \frac{b m_1}{b-a}  } \Big)^{\frac{a}{b-a}}, 
\end{equation}
as long as the parameter $R=T>0$ is chosen sufficiently small to ensure (\ref{def:R}) with this choice of the parameter $0<\eps <2$.
The estimate (\ref{cost:mu}) can be equivalently written as
\begin{equation} \label{cost:mu_1}
\tilde{C}> 2^{\tilde{\alpha}} c_1^{\frac{b}{b-a}}  c_2^{-\frac{a}{b-a}} h(\mu),
\quad \text{ with } \quad
\tilde{\alpha}= \frac{a m_1}{b-a} + 1+ a + \frac{a^2}{b-a}>0
\end{equation}
and
$$h(\mu)=\frac{2^{\mu \tilde{\beta}}}{(2^\mu-1)^{\tilde{\gamma}}},
\quad \text{with}  \quad
\tilde{\gamma}= \frac{a m_1}{b-a}>0, \quad
\tilde{\beta}= 2\tilde{\gamma} + \frac{a b m_1}{(b-a)^2}>0.$$
The function $h$ takes its minimum value on $(0,\infty)$ at $\mu_*$ such that $2^{\mu_*}=\frac{\tilde{\beta}}{\tilde{\beta}-\tilde{\gamma}}$ and
$$h(\mu_*)=\frac{\tilde{\beta}^{\tilde{\beta}}}{\tilde{\gamma}^{\tilde{\gamma}} (\tilde{\beta}-\tilde{\gamma})^{\tilde{\beta}-\tilde{\gamma}}}.$$
If $C$ is a positive constant such that (\ref{meta_thm_cost_estim}) holds, we conclude that there exists $T_0>0$ such that (\ref{ty1_BIS}) holds for all $0<T\leq T_0$.
\end{proof}

\section{Proofs of the main results for Ornstein-Uhlenbeck equations}\label{Sec:Proof}

This section is devoted to the proof of Theorem~\ref{thm_OU_CS_CN}. Let $T>0$ and $(\omega(t))_{t \in [0,T]}$ be a moving control support in $\mathbb{R}^d$.
We assume that the generalized Kalman rank condition (\ref{Kalman_time}) holds at time $T$.
We begin by noticing that the null-controllability on $[0,T]$ from the moving control support $(\omega(t))_{t \in [0,T]}$ of the equation (\ref{syst_general}) is equivalent to the null-controllability of the following equation
\begin{equation} \label{syst_LR1.z}
\left\lbrace \begin{array}{l}
\partial_t f - \frac{1}{2}\textrm{Tr}\big(A(t)A(t)^T\nabla_x^2 f\big) - \big\langle B(t)x, \nabla_x f\big\rangle-\frac{1}{2}\textrm{Tr}\big(B(t)\big)f=\un_{\omega(t)}(x)u, \\
f|_{t=0}=f_0 \in L^2(\rr_x^d),
\end{array}\right.
\end{equation}
as it is sufficient to change the unknown function $f$ to $f\exp(\frac{1}{2}\int_0^t\textrm{Tr}(B(s))ds)$.
This reduction ensures that the adjoint system (\ref{adj_general1}) defined below generates a contraction evolution system on $L^2(\rr_x^d)$.
Next, we observe that the $L^2(\rr^d_x)$-adjoint of the operator
$$\frac{1}{2}\textrm{Tr}\big(A(t)A(t)^T\nabla_x^2\big)+\big\langle B(t)x, \nabla_x\big\rangle+\frac{1}{2}\textrm{Tr}\big(B(t)\big),$$
is given by
\begin{multline*}
\Big(\frac{1}{2}\textrm{Tr}\big(A(t)A(t)^T\nabla_x^2\big)+\big\langle B(t)x, \nabla_x\big\rangle+\frac{1}{2}\textrm{Tr}\big(B(t)\big)\Big)^*\\
=\frac{1}{2}\textrm{Tr}\big(A(t)A(t)^T\nabla_x^2\big)-\big\langle B(t)x, \nabla_x\big\rangle-\frac{1}{2}\textrm{Tr}\big(B(t)\big).
\end{multline*}
By using the Hilbert uniqueness method, see Proposition~\ref{Prop:HUM} in Section~\ref{Appendix:HUM},
the null-controllability on $[0,T]$ of the equation (\ref{syst_LR1.z})
from the moving control set $(\omega(t))_{t \in [0,T]}$ is equivalent to the following observability estimate 
with respect to the moving observation support $(\omega(T-t))_{t \in [0,T]}$,
\begin{equation}\label{dfg2}
\exists C>0, \forall g_0 \in L^2(\rr^d), \quad \|g(T)\|_{L^2(\rr^d)}^2 \leq C\int_{0}^T\|g(t)\|_{L^2(\omega(T-t))}^2dt,
\end{equation}
where $g$ is the mild solution to the Cauchy problem
\begin{equation} \label{adj_general1}
\left\lbrace \begin{array}{l}
\partial_t g(t,x) - \tilde{P}(t)g(t,x)=0, \quad x \in \rr^d,\\
g|_{t=0}=g_0 \in L^2(\rr^d),
\end{array}\right.
\end{equation}
with
$$\tilde{P}(t)=\frac{1}{2}\textrm{Tr}\big(A(T-t)A(T-t)^T\nabla_x^2\big)-\big\langle B(T-t)x, \nabla_x\big\rangle-\frac{1}{2}\textrm{Tr}\big(B(T-t)\big).$$
We deduce from Proposition~\ref{Prop:WP_Hom} and (\ref{iop1}) that there exists an evolution system $(\mathscr{U}(t,s))_{0 \leq s \leq t \leq T}$ of contractions of $L^2(\rr^d)$ such that the mild solution to the Cauchy problem (\ref{adj_general1}) is given by 
\begin{equation}\label{ty30}
\forall g_0 \in L^2(\rr^d), \forall 0 \leq t \leq T, \quad g(t)=\mathscr{U}(t,0)g_0 \in  L^2(\rr^d).
\end{equation}

\subsection{Sufficient condition for null-controllability}

Let $\pi_k: L^2(\rr^d) \rightarrow E_k$ be the orthogonal projections onto $E_k$ the closed subspace of $L^2(\rr^d)$-functions whose Fourier transforms are
supported in the cube $[-k,k]^d$,
\begin{equation}\label{asd100}
E_k=\big\{f \in L^2(\rr^d) : \textrm{supp }\widehat{f} \subset [-k,k]^d\big\}, \quad k \geq 1.
\end{equation}
The normalization of the Fourier transform used throughout this work is the one given by (\ref{norma}).
Instrumental in the proof of assertion $(i)$ in Theorem~\ref{thm_OU_CS_CN} are the following dissipation estimates similar to the ones established 
in~\cite{KB_KPS_JEP} (Proposition~3.2) in the autonomous case and
in~\cite{KB_KPS_JMAA} (Proposition 2.2) in the non-autonomous one. However, the result of~\cite{KB_KPS_JMAA} (Proposition 2.2) cannot be directly applied here and its proof needs to be revisited as follows: 

\medskip

\begin{Prop} \label{Prop_dissip_OU_non_auton}
Let $T>0$. 
We assume that the generalized Kalman rank condition \emph{(\ref{Kalman_time})} holds at time $T$.
Then, there exist some positive constants $0<\tilde{\eps}<T$, $m_1>0$ and $c_2>0$ such that 
\begin{multline}\label{eq6.z}
\forall g_0 \in L^2(\rr^d), \forall 0 \leq s \leq t \leq \tilde{\eps}, \forall k \geq 1, \\  \| (1-\pi_k)\mathscr{U}(t,s) g_0\|_{L^2(\rr^d)} \leq e^{- c_2(t-s)^{m_1}k^2} \|g_0\|_{L^2(\rr^d)},
\end{multline}
where $\mathscr{U}$ is the contraction evolution system on $L^2(\rr^d)$ associated to the adjoint system \eqref{adj_general1}
given by Proposition~\emph{\ref{Prop:WP_Hom}}.
\end{Prop}

\medskip

\begin{proof}
There exists a positive constant $0<\eps <T$ such that $(T-\eps,T+\eps) \subset I$, since $I$ is an open interval and $T \in I$.
We first aim at establishing that there exist some positive constants $c>0$, $0< \tilde{\eps} < \eps$ and a positive integer $m_1 \geq 1$ such that 
\begin{equation}\label{hop1}
\forall \xi \in \rr^d, \forall 0 \leq t \leq \tau \leq \tilde{\eps}, \quad \int_{t}^{\tau}|A(T-s)^TR(t,s)^T\xi|^2ds \geq c  (\tau-t)^{m_1}|\xi|^2,
\end{equation}
with $|\cdot|$ being the Euclidean norm on $\rr^{d}$ and $R$ being the resolvent of the time-varying linear system 
\begin{equation}\label{ty32}
\dot{X}(t)=B(T-t)X(t).
\end{equation}
We recall for instance from~\cite{coron_book} (Proposition~1.5) that this resolvent satisfies the following properties 
\begin{equation}\label{asd1}
\forall t, \tau \in [0,T], \quad R(t,\tau)R(\tau,t)=I_d, \quad  (\partial_2 R)(t,\tau)=-R(t,\tau)B(T-\tau).
\end{equation}
We notice from (\ref{asd2}), (\ref{asd3}) and (\ref{asd1}) that for all $k \geq 0$ and $t, \tau \in [0,T]$, 
\begin{equation}\label{asd4}
 \frac{d^{k}}{d\tau^{k}}\big(A(T-\tau)^TR(t,\tau)^T\big)=(-1)^k\tilde{A}_k(T-\tau)^TR(t,\tau)^T.
\end{equation}
We consider the function
$$f_{\xi}(t,\tau)=\int_{t}^{\tau}|A(T-s)^TR(t,s)^T\xi|^2ds, \quad t, \tau \in [0,T],$$
depending on the parameter $\xi \in \rr^{d}$.
According to (\ref{asd4}), we easily check by the Leibniz formula that 
\begin{multline}\label{re1}
\forall k \geq 0, \forall t, \tau \in [0,T], \\ (\partial_2^{k+1}f_{\xi})(t,\tau)=(-1)^k\sum_{j=0}^k\binom{k}{j} \langle \tilde{A}_{k-j}(T-\tau)^TR(t,\tau)^T\xi,
\tilde{A}_j(T-\tau)^TR(t,\tau)^T\xi\rangle,
\end{multline}
where $\langle \cdot,\cdot \rangle$ denotes the Euclidean scalar product on $\rr^{d}$. 
The generalized Kalman rank condition (\ref{Kalman_time}) holding at time $T$ implies that
\begin{equation}\label{re2}
\bigcap_{k=0}^{+\infty}\textrm{Ker}\big(\tilde{A}_{k}(T)^T\big) \cap \rr^{d}=\{0\}.
\end{equation}
By induction, we easily check from (\ref{re1}) that for all $k \geq 0$, 
\begin{equation}\label{re3}
\forall 0 \leq l \leq 2k+1, \ (\partial_2^lf_{\xi})(0,0)=0 \Longleftrightarrow  \xi \in \bigcap_{j=0}^{k}\textrm{Ker}\big(\tilde{A}_{j}(T)^T\big) \cap \rr^{d}.
\end{equation}
According to (\ref{re1}), (\ref{re2}) and (\ref{re3}), it follows that for all $\xi \in \rr^{d} \setminus \{0\}$, there exists $k_{\xi} \geq 0$ such that 
\begin{equation}\label{re0}
\forall 0 \leq j \leq 2k_{\xi}, \quad (\partial_2^jf_{\xi})(0,0)=0
\end{equation}
and
\begin{equation}\label{re0b}
(\partial_2^{2k_{\xi}+1}f_{\xi})(0,0)=\binom{2k_{\xi}}{k_{\xi}}|\tilde{A}_{k_{\xi}}(T)^T\xi|^2>0.
\end{equation}
We aim at proving that for all $\xi \in \mathbb{S}^{d-1}$ (the unit sphere), there exist positive constants $c_{\xi}>0$, $0< \tilde{\eps}_{\xi} < \eps $ and an open neighborhood $V_{\xi}$ of $\xi$ in $\mathbb{S}^{d-1}$ such that 
\begin{equation}\label{re4}
\forall 0 \leq t < \tau \leq \tilde{\eps}_{\xi}, \forall \eta \in V_{\xi}, \quad \int_{t}^{\tau}|A(T-s)^TR(t,s)^T\eta|^2ds \geq c_{\xi}(\tau-t)^{2k_{\xi}+1},
\end{equation}
By analogy with~\cite{joh} (Proposition~3.2), we proceed by contradiction. If the assertion (\ref{re4}) does not hold, there exist sequences of non-negative real numbers $(t_l)_{l \geq 0}$, $(\tau_l)_{l \geq 0}$ satisfying
\begin{equation}\label{so0}
\forall l \geq 0, \quad 0 \leq t_l < \tau_l < \eps, \qquad \lim_{l \to +\infty}t_l=\lim_{l \to +\infty}\tau_l=0,
\end{equation}
and a sequence $(\eta_l)_{l \geq 0}$ of elements in $\mathbb{S}^{d-1}$ so that 
\begin{equation}\label{re5}
\lim_{l \to +\infty}\eta_l=\xi, 
\end{equation}
and 
\begin{equation}\label{re5b}
\lim_{l \to +\infty}\frac{1}{(\tau_l-t_l)^{2k_{\xi}+1}}\int_{t_l}^{\tau_l}|A(T-s)^TR(t_l,s)^T\eta_l|^2ds=0.
\end{equation}
We deduce from (\ref{re5b}) that 
\begin{equation}\label{re6}
\lim_{l \to +\infty}\frac{1}{(\tau_l-t_l)^{2k_{\xi}+1}}\sup_{0 \leq t \leq \tau_l-t_l}\int_{t_l}^{t_l+t}|A(T-s)^TR(t_l,s)^T\eta_l|^2ds=0.
\end{equation}
Setting
\begin{equation}\label{re7}
u_l(x)=\frac{1}{(\tau_l-t_l)^{2k_{\xi}+1}}\int_{t_l}^{t_l+x(\tau_l-t_l)}|A(T-s)^TR(t_l,s)^T\eta_l|^2ds \geq 0, 
\end{equation}
for any $0 \leq x \leq 1$,
we can reformulate (\ref{re6}) as 
\begin{equation}\label{re8}
\lim_{l \to +\infty}\sup_{0 \leq x \leq 1}|u_l(x)|=0.
\end{equation}
By writing that
\begin{multline}\label{so1}
f_{\eta_l}(t_l,\tau)=\int_{t_l}^{\tau}|A(T-s)^TR(t_l,s)^T\eta_l|^2ds=\sum_{j=0}^{2k_{\xi}+1}a_l^{(j)}(\tau-t_l)^{j}\\
+\frac{(\tau-t_l)^{2k_{\xi}+2}}{(2k_{\xi}+1)!}\int_0^1(1-s)^{2k_{\xi}+1}(\partial_2^{2k_{\xi}+2}f_{\eta_l})\big(t_l,t_l+s(\tau-t_l)\big)ds,
\end{multline}
with $a_l^{(j)}=(\partial_2^jf_{\eta_l})(t_l,t_l)(j!)^{-1}$, we notice from (\ref{re1}) that there exists a positive constant $M>0$ such that 
\begin{multline}\label{so2}
\forall l \geq 0,\forall \tau \in [0,T], \\ \Big|\frac{1}{(2k_{\xi}+1)!}\int_0^1(1-s)^{2k_{\xi}+1}(\partial_2^{2k_{\xi}+2}f_{\eta_l})\big(t_l,t_l+s(\tau-t_l)\big)ds\Big| \leq M.
\end{multline}
It follows from (\ref{re7}), (\ref{so1}) and (\ref{so2}) that 
\begin{equation}\label{re9}
\forall 0 \leq x \leq 1, \forall l \geq 0, \quad \Big|u_l(x)-\sum_{j=0}^{2k_{\xi}+1}\frac{a_l^{(j)}}{(\tau_l-t_l)^{2k_{\xi}+1-j}}x^{j}\Big| \leq M(\tau_l-t_l) x^{2k_{\xi}+2}.
\end{equation}
It follows from (\ref{so0}), (\ref{re8}) and (\ref{re9}) that 
\begin{equation}\label{re10}
\lim_{l \to +\infty}\sup_{0 \leq x \leq 1}|p_l(x)|=0,
\end{equation}
with 
\begin{equation}\label{re11}
p_l(x)=\sum_{j=0}^{2k_{\xi}+1}\frac{a_l^{(j)}}{(\tau_l-t_l)^{2k_{\xi}+1-j}}x^{j}.
\end{equation}
By using the equivalence of norms in finite-dimensional vector spaces, we deduce from (\ref{re10}) that 
\begin{equation}\label{re12}
\forall 0 \leq j \leq 2k_{\xi}+1, \quad \lim_{l \to +\infty}\frac{a_l^{(j)}}{(\tau_l-t_l)^{2k_{\xi}+1-j}}=0.
\end{equation}
We obtain in particular that 
\begin{equation}\label{re13}
\lim_{l \to +\infty}a_l^{(2k_{\xi}+1)}=0.
\end{equation}
According to (\ref{re0b}), this is in contradiction with the fact that 
\begin{equation}\label{re14}
\lim_{l \to +\infty}a_l^{(2k_{\xi}+1)}=\lim_{l \to +\infty}\frac{(\partial_2^{2k_{\xi}+1}f_{\eta_l})(t_l,t_l)}{(2k_{\xi}+1)!}=\frac{1}{(2k_{\xi}+1)!}(\partial_2^{2k_{\xi}+1}f_{\xi})(0,0)>0.
\end{equation}
By covering the compact set $\mathbb{S}^{d-1}$ by finitely many open neighborhoods of the form $(V_{\xi_j})_{1 \leq j \leq N}$, and letting $c=\inf_{1 \leq j \leq N}c_{\xi_j}>0$, $0<\tilde{\eps}=\inf_{1 \leq j \leq N}\{\tilde{\eps}_{\xi_j},1\} < \eps $, $m_1=1+\sup_{1 \leq j \leq N}k_{\xi_j}\geq 1$, we conclude that 
$$\forall \xi \in \rr^d, \forall 0 \leq t \leq \tau \leq \tilde{\eps}, \quad \int_{t}^{\tau}|A(T-s)^TR(t,s)^T\xi|^2ds \geq c  (\tau-t)^{m_1}|\xi|^2.$$
It ends the proof of the estimate (\ref{hop1}). We can now derive the estimates (\ref{eq6.z}).
Let $\mathscr{U}$ be the contraction evolution system associated to the adjoint system \eqref{adj_general1}
given by Proposition~\ref{Prop:WP_Hom}.
We deduce from (\ref{k_explicit}) in Proposition~\ref{Prop:WP_Hom} that for all $0 \leq t \leq \tau \leq T$ and $g_0 \in L^2(\rr^d)$,
\begin{multline}\label{sdf1.h}
\widehat{(\mathscr{U}(\tau,t)g_0)}(\xi) 
= \exp\Big(\frac{1}{2}\int_{t}^{\tau}\textrm{Tr}\big(B(T-s)\big)ds\Big)\\
\times \widehat{g_0}\big(R(\tau,t)^T\xi\big)\exp\Big(-\frac{1}{2}\int_{t}^{\tau}|A(T-s)^TR(\tau,s)^T\xi|^2 ds\Big),
\end{multline}
where the resolvent $R$ is defined in (\ref{ty32}).
We deduce from (\ref{hop1}) and (\ref{sdf1.h}) that for all $0 \leq t \leq \tau \leq \tilde{\eps}$, $k \geq 1$ and $g_0 \in L^2(\rr^d)$, 
\begin{align}\label{sdf3}
& \  \| (1-\pi_k)\mathscr{U}(\tau,t)g_0\|_{L^2(\rr^d)}^2\\ \notag
= & \ \frac{1}{(2\pi)^d}e^{\int_t^{\tau}\textrm{Tr}(B(T-s))ds}
\int_{\xi \in \rr^d \setminus [-k,k]^d}\hspace{-0.5cm}\big|\widehat{g_0}\big(R(\tau,t)^T\xi\big)\big|^2e^{-\int_t^{\tau} |A(T-s)^TR(\tau,s)^T \xi|^2 ds}d\xi\\ \notag
= & \ \frac{1}{(2\pi)^d}\big|\text{det}\big(R(t,\tau)\big)\big|e^{\int_t^{\tau}\textrm{Tr}(B(T-s))ds}
\int_{R(t,\tau)^T\xi \in \rr^d \setminus [-k,k]^d}\hspace{-2.2cm}|\widehat{g_0}(\xi)|^2e^{-\int_t^{\tau} |A(T-s)^TR(t,s)^T \xi|^2 ds}d\xi\\ \notag
\leq & \ \frac{1}{(2\pi)^d}
\int_{|\xi| \geq k \|R(t,\tau)^T\|^{-1}}|\widehat{g_0}(\xi)|^2e^{-c(\tau-t)^{m_1}|\xi|^2}d\xi,
\end{align}
since
$$\forall t_1,t_2,t_3 \in [0,T], \quad R(t_1,t_2)R(t_2,t_3)=R(t_1,t_3),$$
see e.g.~\cite{coron_book} (Proposition~1.5), and by Liouville formula 
$$\forall t, \tau \in [0,T], \quad \text{det}\big(R(\tau,t)\big)=\exp\Big(\int_{t}^{\tau}\textrm{Tr}\big(B(T-s)\big)ds\Big),$$
see e.g.~\cite{cartan} (Proposition II.2.3.1).
We deduce from (\ref{sdf3}) that there exists a positive constant $c_2>0$ such that for all $0 \leq t \leq \tau \leq \tilde{\eps}$, $k \geq 1$ and $g_0 \in L^2(\rr^d)$, 
\begin{equation}\label{sdf5}
\| (1-\pi_k)\mathscr{U}(\tau,t)g_0\|_{L^2(\rr^d)}^2 \leq e^{-2c_2(\tau-t)^{m_1}k^2}\|g_0\|_{L^2(\rr^d)}^2.
\end{equation}
It proves the estimate (\ref{eq6.z}) and ends the proof of Proposition~\ref{Prop_dissip_OU_non_auton}.
\end{proof}

We can now resume the proof of assertion $(i)$ in Theorem~\ref{thm_OU_CS_CN}.
Let $\delta>0$, $\alpha \in (0,+\infty)^d$ and $E$ be a measurable subset of $[0,T]$ satisfying
\begin{equation} \label{hyp:E1}
\exists 0<r_0 \leq T,  \forall  0< r \leq r_0 , \quad \lambda(E\cap[T-r,T])>0.
\end{equation}
We assume that $\omega(t)$ is a $(\delta,\alpha)$-thick subset in $\mathbb{R}^d$ for all $t \in E$.
In order to establish assertion $(i)$ in Theorem~\ref{thm_OU_CS_CN}, we prove that the observability estimate \eqref{dfg2} holds by applying
Theorem~\ref{Meta_thm_AdaptedLRmethod} (formula (\ref{ty1})) with $\mathscr{U}$ the contraction evolution system on $L^2(\rr^d)$ associated to the adjoint system \eqref{adj_general1}
given by Proposition~\ref{Prop:WP_Hom}, $\delta_1=0$, $0<\delta_2=\min(\tilde{\eps},r_0)<T$,  $\tilde{\eps}$ given by Proposition \ref{Prop_dissip_OU_non_auton} and
the moving control support $(\widetilde{\omega}(t))_{t\in [0,T]}$ defined as
$$\widetilde{\omega}(t)=\omega(T-t), \quad t \in [0,T].$$
Proposition~\ref{Prop_dissip_OU_non_auton} shows that the dissipation estimates (\ref{Meta_thm_dissip}) hold with $c_2'=1$, $m_2=0$ and $b=2$. Regarding the uniform spectral estimates (\ref{Meta_thm_IS}), we notice that $\widetilde{\omega}(t)$ is a $(\delta,\alpha)$-thick subset in $\rr^d$ for all $t \in T-E=\{ T-s : s\in E\}$. 
It follows from (\ref{hyp:E1}) that the subset $\widetilde{E}= (T-E) \cap [0,\delta_2]$ is therefore a measurable subset of $[0,\delta_2]$ with positive Lebesgue measure as
$$\lambda(\widetilde{E})=\lambda((T-E)\cap [0,\delta_2])=\lambda(E \cap [T-\delta_2,T])>0.$$ 
We deduce from Theorem~\ref{thm_PLSK} that the uniform spectral estimates (\ref{Meta_thm_IS}) hold with $\widetilde{E} \subset [0,\delta_2]$ and the parameters
\begin{equation}\label{fk1}
a=1, \qquad c_1=2\mathcal{C}|\alpha|\ln\Big(\frac{\mathcal{C}^d}{\delta}\Big)>0, \qquad c_1'=\Big(\frac{\mathcal{C}^d}{\delta}\Big)^{\mathcal{C}d}>0.
\end{equation}
We can therefore deduce (\ref{dfg2}) from Theorem~\ref{Meta_thm_AdaptedLRmethod}. This ends the proof of assertion~$(i)$ in Theorem~\ref{thm_OU_CS_CN}.

\medskip

\begin{rk} \label{O-Uh:cost}
Let $C>0$ be a positive constant verifying 
$$C> \frac{32}{c_2}\Big(\mathcal{C}|\alpha| \ln\Big(\frac{\mathcal{C}^d}{\delta}\Big)\Big)^{2}\Big(\frac{8}{3}\Big)^{3m_1},$$
where the positive constants $c_2>0$ and $m_1>0$ are given by Proposition~\ref{Prop_dissip_OU_non_auton}, and $\mathcal{C}>1$ is the universal constant given by Theorem~\ref{thm_PLSK}. We deduce from Theorem~\ref{Meta_thm_AdaptedLRmethod} that there exists a positive constant $0<T_0<T$ such that, if $E=[\tilde{T},T]$ with $T_0 \leq \tilde{T} <T$, then the following quantitative observability estimate holds
$$ \forall g_0 \in L^2(\rr^d), \quad \|g(T)\|_{L^2(\rr^d)}^2 \leq \exp\Big( \frac{C}{(T-\tilde{T})^{m_1}}\Big)\int_{0}^{T-\tilde{T}}\|g(t)\|_{L^2(\omega(T-t))}^2dt,$$
where $g$ is the mild solution to the Cauchy problem \emph{(\ref{adj_general1})}.
This provides a quantitative estimate of the observability cost with respect to the characteristic parameters related to the thickness property of the moving control support. 
On the other hand, let us notice that the positive parameter $m_1>0$ is actually independent on the control support.
\end{rk}

\medskip

\subsection{Necessary condition  for null-controllability}

This section is devoted to give a proof of assertion $(ii)$ in Theorem~\ref{thm_OU_CS_CN}.
Let $T>0$. We assume that the non-autonomous Ornstein-Uhlenbeck equation (\ref{syst_LR1.z}) is null-controllable on $[0,T]$ from $(\omega(t))_{t \in [0,T]}$,
or equivalently, that the adjoint system (\ref{adj_general1}) is observable on $[0,T]$ from $(\omega(T-t))_{0 \leq t \leq T}$, that is, there exists a positive constant $C>0$ such that any solution of (\ref{adj_general1}) satisfies the observability estimate (\ref{dfg2}).
We aim at finding some positive constants  $r>0$ and $\delta>0$ such that
\begin{equation} \label{Integral_relative_density_BIS}
\forall z \in \mathbb{R}^d, \quad \int_0^T \lambda \big(B_d(z,r) \cap R(0,T-t) \omega(t) \big) dt \geq \delta>0,
\end{equation}
where $B_d(z,r)$ denotes the open Euclidean ball centered at $z$ with radius $r$, and $R$ stands for the resolvent of the time-varying linear system $\dot{X}(t)=B(T-t)X,$
given by the solution of (\ref{Resolvante}). In order to derive this necessary condition, we try out the observability estimate (\ref{dfg2}) with explicit gaussian solutions of (\ref{adj_general1}) centered at~$z\in \rr^d$.
The computation of the Fourier transform of the solution of (\ref{adj_general1}) with respect to the space variable is performed in appendix (Proposition~\ref{Prop:WP_Hom}), see also \cite{KB_KPS_JMAA} (Appendix A.2),
$$\widehat{g}(t,\xi)=
\text{exp}\Big(\frac{1}{2}\int_0^t \text{Tr}\big(B(T-s)\big) ds \Big) \widehat{g}_0(R(t,0)^T\xi) e^{-\frac{1}{2}\int_{0}^t |A(T-s)^T R(t,s)^T\xi|^2 ds},$$
when $0 \leq t \leq T$ and $\xi \in \rr^d$,
or equivalently
\begin{equation}\label{ty0}
\widehat{g}(t,\xi)= \text{exp}\Big(\frac{1}{2}\int_0^t \text{Tr}\left(B(T-s) \right) ds \Big) \widehat{g}_0(R(t,0)^T\xi) e^{-\frac{1}{2}| \sqrt{Q_t} R(t,0)^T \xi|^2},
\end{equation}
where
\begin{equation}\label{ty11}
Q_t= \int_{0}^t R(0,s) A(T-s) A(T-s)^T R(0,s)^T ds,
\end{equation}
is a symmetric positive semidefinite matrix for all $0 \leq t \leq T$.
Let $z\in\rr^d$. We consider the gaussian initial datum
$$g_0(x)= e^{-\frac{|x-z|^2}{2\alpha}},$$
with $\alpha>0$. Since
$$\widehat{g_{0}}(\xi)= (2\pi\alpha)^{\frac{d}{2}} e^{-\frac{\alpha |\xi|^2}{2}} e^{-i z \cdot \xi},$$
we deduce from (\ref{ty0}) that for all $0 \leq t \leq T$ and $\xi \in \rr^d$,
$$\widehat{g}(t,\xi)= (2\pi\alpha)^{\frac{d}{2}} \text{exp}\Big(\frac{1}{2}\int_0^t \text{Tr}\big(B(T-s) \big) ds \Big)
                        e^{-\frac{1}{2} | \sqrt{M_t} R(t,0)^T \xi |^2}  e^{-i z \cdot R(t,0)^T \xi},$$
with $M_t=Q_t+\alpha I_d$.
We deduce from (\ref{ty11}) that $M_t$ is a symmetric positive definite matrix for all $0 \leq t \leq T$, and that the estimate
$$\forall 0 \leq t \leq T, \forall x \in \rr^d, \quad \alpha|x|^2 \leq |\sqrt{M_t}x|^2 =(M_tx) \cdot x \leq (\alpha+\|Q_T\|)|x|^2 ,$$
implies that
\begin{equation}\label{ty13}
\forall 0 \leq t \leq T, \quad 0< \alpha^d \leq \det M_t
\end{equation}
and for all $0 \leq t \leq T$ and $x \in \rr^d$, 
\begin{equation}\label{ty10}
 \frac{1}{\sqrt{\alpha+\|Q_T\|}}|x|  \leq  |\sqrt{M_t^{-1}}x|=|\sqrt{M_t}^{-1}x| \leq \frac{1}{\sqrt{\alpha}}|x|.
\end{equation}
By using Liouville formula
\begin{equation}\label{ty12}
\forall t, \tau \in [0,T], \quad \det(R(\tau,t))=\exp\Big(\int_t^{\tau}\textrm{Tr}(B(T-s))ds\Big),
\end{equation}
see e.g.~\cite{cartan} (Proposition II.2.3.1),
and computing the inverse Fourier transform, we obtain that
$$g(t,x)= \frac{\alpha^{\frac{d}{2}}}{\sqrt{\det M_t}} \text{exp}\Big(-\frac{1}{2}\int_0^t \text{Tr}\big(B(T-s) \big) ds \Big)e^{-\frac{1}{2} |\sqrt{M_t^{-1}} ( R(0,t)x - z )|^2 }.$$
By using the substitution rule with $y=R(0,T)x-z$, and anew the Liouville formula, we notice that the left-hand-side term in the observability estimate (\ref{dfg2}) is a positive constant independent of the parameter $z$,
\begin{equation}\label{lop1}
C_1=\int_{\mathbb{R}^d} |g(T,x)|^2 dx =\frac{\alpha^d}{\det M_T}  \int_{\mathbb{R}^d} e^{- |\sqrt{M_T^{-1}}  y|^2}dy >0.
\end{equation}
According to (\ref{ty13}), (\ref{ty10}), (\ref{ty12}) and the substitution rule with $y=R(0,t)x-z$, the right-hand-side of the observability estimate (\ref{dfg2}) can be bounded from above up to the positive constant $C>0$ as 
\begin{align*}
  & \int_0^T \int_{\omega(T-t)}  |g(t,x)|^2 dx dt
\\
= & \int_0^T \frac{\alpha^{d}}{ \det M_t } \text{exp}\Big(-\int_0^t \text{Tr}\big(B(T-s)\big)ds \Big)\Big( \int_{\omega(T-t)}  e^{- |\sqrt{M_t^{-1}}(R(0,t)x - z)|^2 } dx\Big) dt
\\
= & \int_0^T \frac{\alpha^{d}}{ \det M_t }
\Big(\int_{R(0,t)\omega(T-t)-z}  e^{-|\sqrt{M_t^{-1}} y|^2 } dy\Big) dt
\\
\leq & \int_0^T   \int_{R(0,t)\omega(T-t)-z} e^{-\frac{|y|^2}{\alpha+\|Q_T\|}} dy  dt
\\
\leq &  \int_0^T \Big(\int_{B_d(0,r) \cap [R(0,t)\omega(T-t)-z]}e^{-\frac{|y|^2}{\alpha+\|Q_T\|}} dy\Big) dt +  \int_0^T \Big(\int_{|y|>r} e^{-\frac{|y|^2}{\alpha+\|Q_T\|}}dy\Big) dt
\\
\leq &  \int_0^T \lambda \big( B_d(0,r) \cap [R(0,t) \omega(T-t)-z] \big) dt + T\int_{|y|>r} e^{-\frac{|y|^2}{\alpha+\|Q_T\|}} dy.
\end{align*}
for any $r>0$.
Since by dominated convergence theorem
$$\lim_{r \to +\infty}\int_{|y|>r} e^{-\frac{|y|^2}{\alpha+\|Q_T\|}} dy=0,$$
there exists a positive constant $r_0>1$ such that 
$$CT \int_{|y|>r_0} e^{-\frac{|y|^2}{\alpha+\|Q_T\|}} dy \leq \frac{C_1}{2}.$$
It then follows from (\ref{dfg2}) and (\ref{lop1}) that for all $z \in \rr^d$,
$$0<\frac{C_1}{2} \leq C \int_0^T \lambda \big(B_d(0,r_0) \cap [R(0,t) \omega(T-t)-z] \big) dt,$$
or equivalently by translation invariance of the Lebesgue measure that for all $z \in \rr^d$,
$$0<\frac{C_1}{2 C} \leq  \int_0^T \lambda \big( B_d(z,r_0) \cap R(0,t) \omega(T-t) \big) dt=\int_0^T \lambda \big( B_d(z,r_0) \cap R(0,T-t) \omega(t) \big) dt.$$
It establishes the integral thickness condition (\ref{Integral_relative_density_BIS}) with $\delta = \frac{C_1}{2 C}$. This ends the proof of Theorem~\ref{thm_OU_CS_CN} as we have already checked that assertion $(iii)$ is a direct consequences of assertions $(i)$ and $(ii)$ in Theorem~\ref{thm_OU_CS_CN}.

\section{Proof of the null-controllability for quadratic equations with zero singular spaces}
\label{sec:quad}

This section is devoted to the proof of Theorem~\ref{Main_result_quad}. Let $q : \rr_{x}^{d} \times \rr_{\xi}^d \rightarrow \cc$ be a complex-valued quadratic form with a non-negative real part $\textrm{Re }q \geq 0$, and a zero singular space $S=\{0\}$.
Let $T>0$, $\delta>0$, $\alpha \in (0,+\infty)^d$, $(\omega(t))_{t \in [0,T]}$ be a moving control support in $\mathbb{R}^d$ 
and $E$ be a measurable subset of $[0,T]$ with positive Lebesgue measure $\lambda(E)>0$. We assume that $\omega(t)$ is a $(\delta,\alpha)$-thick subset of $\mathbb{R}^d$ for all $t \in E$. By the Hilbert uniqueness method \cite{coron_book} (Theorem~2.44) and its extension to the moving control support case as presented in Proposition~\ref{Prop:HUM} for Ornstein-Uhlenbeck equations\footnote{The very same proof can be readily adapted to the quadratic case.}, the result of null-controllability on $[0,T]$ from the moving control support $(\omega(t))_{t \in [0,T]}$ given by Theorem \ref{Main_result_quad} is equivalent to the following observability estimate 
$$\int\limits_{\mathbb{R}^d} |g(T,x)|^2 dx \leq C_T \int\limits_0^T \int\limits_{\omega(T-t)} |g(t,x)|^2 dx dt,$$
for the adjoint system
\begin{equation} \label{syst_adj.1}
\left\lbrace \begin{array}{ll}
\partial_t g(t,x) +q^w(x,D_x)^*g(t,x)=0, \quad & (t,x) \in (0,+\infty) \times \mathbb{R}^d, \\
g|_{t=0} = g_0 \in L^2(\mathbb{R}^d),
\end{array} \right.
\end{equation}
that is, there exists a positive constant $C_T>0$ such that
$$\forall g_0 \in L^2(\rr^d), \quad \|e^{-T(q^w)^*}g_0\|_{L^2(\rr^d)}^2 \leq C_T\int_0^T\|e^{-t(q^w)^*}g_0\|_{L^2(\omega(T-t))}^2dt.$$
The $L^2(\rr^d)$-adjoint of the quadratic operator $(q^w,D(q^w))$ defined in (\ref{3}) is given by the quadratic operator $(\overline{q}^w, D(\overline{q}^w))$, 
whose Weyl symbol is the complex conjugate $\overline{q}$ of the symbol $q$. 
We notice that the assumptions of Theorem~\ref{Main_result_quad} hold for the quadratic operator $P = q^w(x, D_x)$ 
if and only if they hold as well for its $L^2(\rr^d)$-adjoint operator $P^* = \overline{q}^w(x,D_x)$. In order to prove Theorem~\ref{Main_result_quad}, it is therefore sufficient to prove that, if $q : \rr_{x}^{d} \times \rr_{\xi}^d \rightarrow \cc$ is a complex-valued quadratic form with a non-negative real part $\textrm{Re }q \geq 0$ and a zero singular space $S=\{0\}$, and if $E$ is a measurable subset of $[0,T]$ with positive Lebesgue measure $\lambda(E)>0$, such that $\omega(t)$ is a $(\delta,\alpha)$-thick subset of $\mathbb{R}^d$ for all $t \in E$, then there exists a positive constant $C_T>0$ such that
\begin{equation}\label{up3}
\forall g_0 \in L^2(\rr^d), \quad \|e^{-Tq^w}g_0\|_{L^2(\rr^d)}^2 \leq C_T\int_0^T\|e^{-tq^w}g_0\|_{L^2(\omega(t))}^2dt.
\end{equation}
We establish the observability estimate (\ref{up3}) by applying Theorem~\ref{Meta_thm_AdaptedLRmethod} (formula (\ref{ty1}))
with the contraction evolution system $\mathscr{U}(t,s)=e^{-(t-s)q^w}$, for $0 \leq s \leq t$, defined by the contraction semigroup $(e^{-tq^w})_{t \geq 0}$ on $L^2(\rr^d)$. As in the proof of Theorem~\ref{thm_OU_CS_CN}, we consider anew the orthogonal projections $\pi_k: L^2(\rr^d) \rightarrow E_k$ onto $E_k$ the closed subspace of $L^2(\rr^d)$-functions whose Fourier transforms are supported in the cube $[-k,k]^d$, with $k\geq 1$.
We deduce from Theorem~\ref{thm_PLSK} that the uniform spectral estimates (\ref{Meta_thm_IS}) hold with the subset $E$ and the parameters
\begin{equation}\label{fk1.h}
a=1, \qquad c_1=2\mathcal{C}|\alpha|\ln\Big(\frac{\mathcal{C}^d}{\delta}\Big)>0, \qquad c_1'=\Big(\frac{\mathcal{C}^d}{\delta}\Big)^{\mathcal{C}d}>0.
\end{equation}
The next result establishes the dissipation estimates (\ref{Meta_thm_dissip}) with $\delta_1=0$, $\delta_2=t_0$, $m_1=2k_0+1$, $m_2=(2k_0+1)(d+1)$ and $b=2$. Let us point out that Proposition~\ref{Prop:dissip_quad} does not rule out some blow-up for small times in the dissipation estimates. The assumption of dissipation estimates with controlled blow-up (\ref{Meta_thm_dissip}) is then really essential and the result of~\cite{KB_KPS_JEP} (Theorem~2.1) does not apply even in the case of fixed control subsets. The observability estimate (\ref{up3}) is then deduced by applying Theorem~\ref{Meta_thm_AdaptedLRmethod} (formula (\ref{ty1})). Up to the proof of Proposition~\ref{Prop:dissip_quad}, the one of Theorem~\ref{Main_result_quad} is then complete.

\medskip

\begin{Prop} \label{Prop:dissip_quad}
There exist some positive constants $c_2>0$, $c_2'>0$, $0<t_0<T$ such that for all $0 < t \leq t_0$, $g \in L^2(\rr^d)$, $k \geq 1$,
$$\|  (1-\pi_k) e^{-t q^w} g \|_{L^2(\rr^d)} \leq \frac{1}{c_2'} \frac{1}{t^{(2k_0+1)(d+1)}}e^{-c_2 t^{2k_0+1} k^2}  \|  g \|_{L^2(\rr^d)},$$
where $0 \leq k_0 \leq 2d-1$ denotes the smallest integer satisfying  
$$\Big( \bigcap_{j=0}^{k_0}\emph{\textrm{Ker}}\big[\emph{\textrm{Re }}F(\emph{\textrm{Im }}F)^j \big]\Big)\cap \rr^{2d}=\{0\}.$$
\end{Prop}

\medskip

\begin{proof}
The proof of Proposition~\ref{Prop:dissip_quad} relies on basic estimates on Hermite functions recalled in Lemma~\ref{ge3bis}, see appendix (Section~\ref{6.sec.harmo}). 
In the work~\cite{HPSVII} (Theorem~1.2), Hitrik, Viola and the third author have shown that the contraction semigroup $(e^{-tq^w})_{t \geq 0}$ is smoothing for any positive time in the Gelfand-Shilov space $S_{1/2}^{1/2}(\rr^d)$,
$$\forall g \in L^2(\rr^d), \forall t>0, \quad e^{-tq^w}g \in S_{1/2}^{1/2}(\rr^d).$$
We refer the reader to the appendix (Section~\ref{GSreg}) for the definition and some characterizations of the Gelfand-Shilov regularity. 
More specifically, we deduce from~\cite{HPSVII} (Proposition~4.1) that there exist some positive constants $C_0>1$ and $t_0>0$ such that
\begin{equation}\label{eq3}
\forall 0 \leq t \leq t_0, \quad \Big\| e^{\frac{t^{2k_0+1}}{C_0}(-\Delta_x+|x|^2)}e^{-tq^w} \Big\|_{\mathcal{L}_c(L^2(\rr^d))} \leq C_0,
\end{equation} 
with $\mathcal{L}_c(L^2(\rr^d))$ the space of bounded linear operators on $L^2(\rr^d)$, that is,
\begin{multline}\label{eq3..1}
\forall 0 \leq t \leq t_0, \forall g \in L^2(\rr^d), \\ \sum_{\gamma \in \nn^d}|(e^{-tq^w}g,\Phi_{\gamma})_{L^2(\rr^d)}|^2e^{\frac{t^{2k_0+1}}{C_0}(4|\gamma|+2d)} \leq C_0^2 \|g\|_{L^2(\rr^d)}^2,
\end{multline} 
where $(\Phi_{\gamma})_{\gamma \in \nn^d}$ denotes the Hermite basis whose definition is recalled in appendix (Section~\ref{6.sec.harmo}).
We deduce from (\ref{eq3..1}) that 
\begin{multline}\label{eq3..2}
\forall 0 \leq t \leq t_0, \forall \gamma \in \nn^d, \forall g \in L^2(\rr^d),  \\ |(e^{-tq^w}g,\Phi_{\gamma})_{L^2(\rr^d)}|\leq C_0e^{-\frac{t^{2k_0+1}}{C_0}(2|\gamma|+d)}  \|g\|_{L^2(\rr^d)}.
\end{multline} 
For any $0<t \leq t_0$ and $g \in L^2(\rr^d)$, the exponential decay of the Hermite coefficients (\ref{eq3..2}) ensures that the function 
$$e^{-tq^w}g=\sum_{\gamma \in \nn^d}(e^{-tq^w}g,\Phi_{\gamma})_{L^2(\rr^d)}\Phi_{\gamma} \in \mathscr{S}(\rr^d),$$ 
belongs to the Schwartz space with convergence of the above series in the Schwartz space $\mathscr{S}(\rr^d)$, see for instance~\cite{toft} (Proposition~1.2). By continuity of the operator $x^{\alpha}\partial_x^{\beta} : \mathscr{S}(\rr^d) \rightarrow \mathscr{S}(\rr^d)$, with $\alpha, \beta \in \nn^d$, we obtain that 
$$x^{\alpha}\partial_x^{\beta}(e^{-tq^w}g)=\sum_{\gamma \in \nn^d}(e^{-tq^w}g,\Phi_{\gamma})_{L^2(\rr^d)}x^{\alpha}\partial_x^{\beta}\Phi_{\gamma} \in \mathscr{S}(\rr^d),$$
with convergence of the series in $\mathscr{S}(\rr^d)$, and a fortiori in $L^2(\rr^d)$. We deduce from (\ref{eq3..2}) and Lemma~\ref{ge3bis} with $r=1/2$ that for all $\alpha, \beta \in \nn^d$, $\eps>0$, $0<t \leq t_0$, $g \in L^2(\rr^d)$,
\begin{multline}\label{jj1}
\|x^{\alpha}\partial_x^{\beta}(e^{-tq^w}g)\|_{L^2(\rr^d)} \leq \sum_{\gamma \in \nn^d}|(e^{-tq^w}g,\Phi_{\gamma})_{L^2(\rr^d)}|\|x^{\alpha}\partial_x^{\beta}\Phi_{\gamma}\|_{L^2(\rr^d)}\\
\leq 2^{\frac{d}{2}}C_0\Big(\frac{2^{\frac{3}{2}}e^{\frac{1}{2}}}{\inf(\eps^{\frac{1}{2}},1)}\Big)^{|\alpha|+|\beta|}\sqrt{\alpha!}\sqrt{\beta!}\|g\|_{L^2(\rr^d)} 
\sum_{\gamma \in \nn^d}e^{-\frac{t^{2k_0+1}}{C_0}(2|\gamma|+d)}  
e^{\frac{\eps}{2}d |\gamma|}.
\end{multline}
With the choice $\eps=\frac{2t^{2k_0+1}}{dC_0}>0$, it follows from (\ref{jj1}) that for all $\alpha, \beta \in \nn^d$, $0<t \leq t_0$, $g \in L^2(\rr^d)$,
\begin{multline}\label{jj2}
\|x^{\alpha}\partial_x^{\beta}(e^{-tq^w}g)\|_{L^2(\rr^d)} \\
\leq 2^{\frac{d}{2}}C_0\Big(\frac{2^{\frac{3}{2}}e^{\frac{1}{2}}}{\inf((\frac{2t^{2k_0+1}}{dC_0})^{\frac{1}{2}},1)}\Big)^{|\alpha|+|\beta|}\sqrt{\alpha!}\sqrt{\beta!}\|g\|_{L^2(\rr^d)} 
\sum_{\gamma \in \nn^d}e^{-\frac{t^{2k_0+1}}{C_0}|\gamma|}.
\end{multline}
We notice that
\begin{align}\label{jj6}
& \ \sum_{\gamma \in \nn^d}e^{-\frac{t^{2k_0+1}}{C_0}|\gamma|}=\sum_{k=0}^{+\infty}\binom{k+d-1}{k}e^{-\frac{t^{2k_0+1}}{C_0}k} \\ \notag
\leq & \  \frac{1}{(d-1)!}\sum_{k=0}^{+\infty}(k+d-1)^{d-1}e^{-\frac{t^{2k_0+1}}{C_0}k}\\ \notag
\leq & \ \frac{1}{(d-1)!}\sum_{k=0}^{+\infty}\frac{1}{(k+1)^2}(k+d)^{d+1}e^{-\frac{t^{2k_0+1}}{C_0}(k+d)}e^{\frac{t^{2k_0+1}}{C_0}d}.
\end{align}
It follows from (\ref{jj6}) that for all $0<t \leq t_0$, 
\begin{multline}\label{jj6hj}
 \sum_{\gamma \in \nn^d}e^{-\frac{t^{2k_0+1}}{C_0}|\gamma|}
\leq \frac{1}{(d-1)!}\Big(\sup_{x \geq 1}x^{d+1}e^{-\frac{t^{2k_0+1}}{C_0}x}\Big)\Big(\sum_{k=1}^{+\infty}\frac{1}{k^2}\Big)e^{\frac{t_0^{2k_0+1}}{C_0}d}\\
\leq \frac{d(d+1) C_0^{d+1}}{t^{(2k_0+1)(d+1)}}\Big(\sum_{k=1}^{+\infty}\frac{1}{k^2}\Big)e^{\frac{t_0^{2k_0+1}}{C_0}d},
\end{multline}
since
$$\forall x \geq 1, \quad x^{d+1}=\frac{1}{(d+1)!}\Big(\frac{t^{2k_0+1}x}{C_0}\Big)^{d+1} \frac{(d+1)! C_0^{d+1}}{t^{(2k_0+1)(d+1)}} \leq \frac{(d+1)! C_0^{d+1}}{t^{(2k_0+1)(d+1)}}e^{\frac{t^{2k_0+1}}{C_0}x}.$$
We deduce from (\ref{jj2}) and (\ref{jj6hj}) that there exists a positive constant $C_1>1$ such that for all $\alpha, \beta \in \nn^d$, $0<t \leq t_0$, $g \in L^2(\rr^d)$,
\begin{equation}\label{jj3}
\|x^{\alpha}\partial_x^{\beta}(e^{-tq^w}g)\|_{L^2(\rr^d)} 
\leq \frac{C_1^{1+|\alpha|+|\beta|}}{t^{\frac{2k_0+1}{2}(|\alpha|+|\beta|+2d+2)}}\sqrt{\alpha!}\sqrt{\beta!}\|g\|_{L^2(\rr^d)}.
\end{equation}
By using the Parseval formula, we deduce from (\ref{jj3}) that for all $\beta \in \nn^d$, $0<t \leq t_0$, $g \in L^2(\rr^d)$,
\begin{multline}\label{jj4}
\|\xi^{\beta}\widehat{e^{-tq^w}g}\|_{L^2(\rr^d)}=\|\mathcal{F}(\partial_x^{\beta}(e^{-tq^w}g))\|_{L^2(\rr^d)}=(2\pi)^{\frac{d}{2}}\|\partial_x^{\beta}(e^{-tq^w}g)\|_{L^2(\rr^d)}\\
\leq (2\pi)^{\frac{d}{2}}\frac{C_1^{1+|\beta|}}{t^{\frac{2k_0+1}{2}(|\beta|+2d+2)}}\sqrt{\beta!}\|g\|_{L^2(\rr^d)},
\end{multline}
where $\mathcal{F}$ denotes the Fourier transform on $\rr^d$. We obtain from (\ref{jj4}) that 
 for all $0<t \leq t_0$, $g \in L^2(\rr^d)$,
\begin{align}\label{jj5}
& \ \Big\|\exp\Big(\frac{t^{2k_0+1}}{8C_1^2}|\xi|^2\Big)\widehat{e^{-tq^w}g}\Big\|_{L^2(\rr^d)}
\leq \sum_{\beta \in \nn^d}   \frac{1}{\beta!} \Big(\frac{t^{2k_0+1}}{8 C_1^2}\Big)^{|\beta|} \|\xi^{2\beta}\widehat{e^{-tq^w}g}\|_{L^2(\rr^d)}\\ \notag
\leq & \ (2\pi)^{\frac{d}{2}}\sum_{\beta \in \nn^d}   \frac{1}{\beta!} \Big(\frac{t^{2k_0+1}}{8 C_1^2}\Big)^{|\beta|} \frac{C_1^{2|\beta|+1}}{t^{(2k_0+1)(|\beta|+d+1)}}\sqrt{(2\beta)!}  \|g\|_{L^2(\rr^d)}\\  \notag
\leq & \ (2\pi)^{\frac{d}{2}}\frac{C_1}{t^{(2k_0+1)(d+1)}}\Big(\sum_{\beta \in \nn^d}\frac{1}{4^{|\beta|}}\Big) \|g\|_{L^2(\rr^d)},
\end{align}
since $\sqrt{(2\beta)!} \leq 2^{|\beta|} \beta!$. By using that
$$\sum_{\beta \in \nn^d}\frac{1}{4^{|\beta|}}=\sum_{k=0}^{+\infty}\binom{k+d-1}{k}\frac{1}{4^k} \leq 2^{d-1}\sum_{k=0}^{+\infty}\frac{1}{2^k}=2^d,$$
we obtain that for all $0<t \leq t_0$, $g \in L^2(\rr^d)$,
\begin{equation}\label{jj5hj}
 \Big\|\exp\Big(\frac{t^{2k_0+1}}{8C_1^2}|\xi|^2\Big)\widehat{e^{-tq^w}g}\Big\|_{L^2(\rr^d)}
\leq \frac{2^d(2\pi)^{\frac{d}{2}} C_1}{t^{(2k_0+1)(d+1)}}\|g\|_{L^2(\rr^d)}.
\end{equation}
It follows that there exists a positive constant $C_2>1$ such that for all $0<t \leq t_0$, $g \in L^2(\rr^d)$,
\begin{equation}\label{jj7}
\Big\|\exp\Big(\frac{t^{2k_0+1}}{8C_1^2}|\xi|^2\Big)\widehat{e^{-tq^w}g}\Big\|_{L^2(\rr^d)} \leq \frac{C_2}{t^{(2k_0+1)(d+1)}}\|g\|_{L^2(\rr^d)}.
\end{equation}
We obtain that for all $0<t \leq t_0$, $g \in L^2(\rr^d)$, $k \geq 0$,
\begin{align}\label{jj8}
& \ \|(1-\un_{[-k,k]^d}(\xi))\widehat{e^{-tq^w}g}\|_{L^2(\rr^d)}\\ \notag
=& \ \Big\|(1-\un_{[-k,k]^d}(\xi))\exp\Big(-\frac{t^{2k_0+1}}{8C_1^2}|\xi|^2\Big)\exp\Big(\frac{t^{2k_0+1}}{8C_1^2}|\xi|^2\Big)\widehat{e^{-tq^w}g}\Big\|_{L^2(\rr^d)}\\ \notag
\leq & \ \exp\Big(-\frac{t^{2k_0+1}}{8C_1^2}k^2\Big)\Big\|\exp\Big(\frac{t^{2k_0+1}}{8C_1^2}|\xi|^2\Big)\widehat{e^{-tq^w}g}\Big\|_{L^2(\rr^d)}\\
 \leq & \ \frac{C_2}{t^{(2k_0+1)(d+1)}}\exp\Big(-\frac{t^{2k_0+1}}{8C_1^2}k^2\Big)\|g\|_{L^2(\rr^d)}, \notag
\end{align}
where $\un_{[-k,k]^d}$ denotes the characteristic function of the set $[-k,k]^d$. This ends the proof of Proposition~\ref{Prop:dissip_quad}.
\end{proof}

\section{Appendix}\label{appendix}

\subsection{Well-posedness of the homogeneous and inhomogeneous Cauchy problems}\label{appensix:well-pos}

This subsection is devoted to recall from~\cite{KB_KPS_JMAA} (Appendix A.1) the well-posedness of the homogeneous and inhomogeneous Cauchy problems for non-autonomous Ornstein-Uhlenbeck equations.

We first study the well-posedness of the homogeneous equation 
\begin{equation} \label{CY_Hom}
\left\lbrace \begin{array}{l}
\partial_t k(t,x) - P_0(t) k(t,x)=0, \quad  (t,x) \in (t_0,T_1) \times \mathbb{R}^d, \\
k|_{t=t_0} = k_0 \in L^2(\mathbb{R}^d),
\end{array} \right.
\end{equation}
associated to the non-autonomous Ornstein-Uhlenbeck operator
\begin{equation} \label{def:P_appendix}
P_0(t)=\frac{1}{2} \text{Tr}\big(A_0(t)A_0(t)^T\nabla_x^2\big)-\big\langle B_0(t)x,\nabla_x\big\rangle-\frac{1}{2}\text{Tr}\big(B_0(t)\big), \quad t \in (T_0,T_1),
\end{equation}
with $T_0 \leq t_0 < T_1$ and $A_0, B_0 \in C^0([T_0,T_1],M_d(\mathbb{R}))$.
In order to define the concept of weak solutions, we introduce the space $E(t_0,T_1)$ of functions $\varphi \in C^0([t_0,T_1],L^2(\mathbb{R}^d))$ satisfying

\medskip

\begin{itemize}
\item[$(i)$] $\varphi(\cdot,x) \in C^1((t_0,T_1),\cc)$ for all $x \in \mathbb{R}^d$
\item[$(ii)$] $\varphi(t,\cdot) \in C^2(\mathbb{R}^d,\cc)$ for all $t \in (t_0,T_1)$
\item[$(iii)$] The functions $\partial_t \varphi+\langle B_0(t) x , \nabla_x \varphi \rangle$, $\nabla_x^2 \varphi$, $\varphi$ belong to  $L^2((t_0,T_1)\times\mathbb{R}^d)$
\end{itemize}

\medskip

We consider the following notion of weak solutions:

\medskip

\begin{Def}
Let $T_0 \leq t_0 < T_1$, $A_0, B_0 \in C^0([T_0,T_1],M_d(\mathbb{R}))$ and $k_0 \in L^2(\mathbb{R}^d)$.
A weak solution to the Cauchy problem \emph{(\ref{CY_Hom})} is a function $k \in C^0([t_0,T_1],L^2(\mathbb{R}^d))$ such that
$k(t_0)=k_0$ in $L^2(\mathbb{R}^d)$, and satisfying for all $\varphi \in E(t_0,T_1)$, $t_* \in (t_0,T_1)$,
$$\int\limits_{\mathbb{R}^d} \big(k(t_*,x) \varphi(t_*,x)-k_0(x)\varphi(t_0,x)\big) dx
= \int\limits_{t_0}^{t_*} \int\limits_{\mathbb{R}^d} k(t,x) \big(\partial_t \varphi(t,x)+P_0(t)^*\varphi(t,x) \big) dx dt,$$
with 
\begin{equation} \label{def:P*_appendix}
P_0(t)^*=\frac{1}{2}\emph{\text{Tr}}\big(A_0(t)A_0(t)^T\nabla_x^2\big)+\big\langle B_0(t)x,\nabla_x\big\rangle+\frac{1}{2}\emph{\text{Tr}}\big(B_0(t)\big), \quad t \in (T_0,T_1).
\end{equation}
\end{Def}

\medskip

We establish the following result:

\medskip

\begin{Prop} \label{Prop:WP_Hom}
Let $T_0<T_1$, $A_0, B_0 \in C^0([T_0,T_1],M_d(\mathbb{R}))$ and $\mathscr{T}=\{ (t,t_0) : T_0 \leq t_0 \leq t \leq T_1 \}$. 
There exists a strongly continuous mapping
$$\begin{array}{cccl}
\mathscr{U}: & \mathscr{T} & \rightarrow & \mathcal{L}_c(L^2(\mathbb{R}^d)), \\
   & (t,t_0)     & \mapsto     & \mathscr{U}(t,t_0),
\end{array}$$
with $\mathcal{L}_c(L^2(\mathbb{R}^d))$ denoting the space of bounded linear operators on $L^2(\mathbb{R}^d)$, satisfying

\medskip

\begin{itemize}
\item[$(i)$] $\forall T_0 \leq t \leq T_1, \quad \mathscr{U}(t,t)=\emph{\textrm{Id}}_{L^2(\mathbb{R}^d)}$
\item[$(ii)$] $\forall T_0 \leq t_0 \leq t_1 \leq t_2 \leq T_1, \quad \mathscr{U}(t_2,t_1)\mathscr{U}(t_1,t_0) = \mathscr{U}(t_2,t_0)$
\item[$(iii)$] For all $k_0 \in L^2(\mathbb{R}^d)$, the function $k(t)=\mathscr{U}(t,t_0)k_0$ is the unique weak solution to the Cauchy problem \emph{(\ref{CY_Hom})}
\end{itemize}

\medskip

\noindent
Furthermore, the Fourier transform of the function $k(t)=\mathscr{U}(t,t_0)k_0$ is given by
\begin{equation}\label{k_explicit}
\widehat{k}(t,\xi) = \widehat{k}_0\big(R_0(t,t_0)^T \xi \big)  e^{\frac{1}{2} \int_{t_0}^t\emph{\text{Tr}}(B_0(s))ds}e^{-\frac{1}{2} \int_{t_0}^t |A_0(s)^TR_0(t,s)^T \xi|^2ds},
\end{equation}
where $R_0$ denotes the resolvent associated to the linear time-varying system $\dot{X}(t)=B_0(t) X(t)$, that is, for all $T_0 \leq t_0, t \leq T_1$,
\begin{equation}\label{rt9}
\left\lbrace \begin{array}{ll}
\frac{\partial R_0}{\partial t} (t,t_0)=B_0(t) R_0(t,t_0), \\
R_0(t_0,t_0)=I_d.
\end{array} \right.
\end{equation}
\end{Prop}

\medskip

In the above statement, the normalization of the Fourier transform with respect to the space variable is given by
\begin{equation}\label{norma}
\widehat{k}(t,\xi)=\int_{\mathbb{R}^d} k(t,x) e^{-i x \cdot \xi} dx.
\end{equation}
Following~\cite{pazy} (Chapter 5, Section~5.1, Definition~5.3, p.~129), the two parameter family of bounded linear operators $(\mathscr{U}(t_1,t_2))_{(t_1,t_2) \in \mathscr{T}}$ is called the evolution system associated to the homogeneous equation (\ref{CY_Hom}). More specifically, we shall say that the mapping $\mathscr{U}(t,t_0)$ is the evolution mapping associated to the family of operators $s \in [t_0,t] \mapsto P_0(s)$.

\begin{proof} Let $T_0 \leq t_0 \leq T_1$ and $k_0 \in L^2(\mathbb{R}^d)$.

\medskip

\noindent
\textit{Step 1.} We first derive heuristically an explicit expression of the Fourier transform~$\widehat{k}$.
To that end, we consider $k$ a smooth solution to the Cauchy problem (\ref{CY_Hom}) and define the function $K : [t_0,T_1] \times \mathbb{R}^d \rightarrow \mathbb{C}$ by
\begin{equation} \label{CVAR:k/K}
k(t,x)=K\big(t, R_0(t_0,t)x\big).
\end{equation}
We recall for instance from~\cite{coron_book} (Proposition~1.5) that for all $T_0 \leq t_1,t_2,t_3 \leq T_1$, 
\begin{equation}\label{rt1}
 R_0(t_1,t_2)R_0(t_2,t_1)=I_d, \quad R_0(t_1,t_2)R_0(t_2,t_3)=R_0(t_1,t_3)
\end{equation}
and
\begin{equation}\label{rt2}
\forall T_0 \leq t_1,t_2 \leq T_1, \quad (\partial_2 R_0)(t_1,t_2)= - R_0(t_1,t_2)B_0(t_2).
\end{equation}
According to (\ref{rt1}) and (\ref{rt2}), the function $K$ is well-defined and a direct computation provides that 
\begin{equation}\label{rt3}
(\partial_t k)(t,x) + \big\langle B_0(t)x,(\nabla_x k)(t,x)\big\rangle = (\partial_t K)\big(t,R_0(t_0,t)x\big).
\end{equation}
It follows from (\ref{CY_Hom}) and (\ref{rt3}) that 
$$\left\lbrace \begin{array}{l}
\partial_t K(t,y) - \frac{1}{2} \text{Tr}\big(R_0(t_0,t) A_0(t) A_0(t)^T R_0(t_0,t)^T  \nabla_y^2 K(t,y) \big) \\
\qquad \qquad \qquad \qquad \qquad \qquad \qquad \qquad + \frac{1}{2} \text{Tr}\big(B_0(t)\big) K(t,y)=0, \\
K|_{t=t_0} = k_0 \in L^2(\mathbb{R}^d).
\end{array} \right.$$
By taking the Fourier transform, we deduce that 
$$\left\lbrace \begin{array}{l}
\partial_t \widehat{K}(t,\eta) + \frac{1}{2}|A_0(t)^TR_0(t_0,t)^T \eta|^2 \widehat{K}(t,\eta)  + \frac{1}{2} \text{Tr}\big(B_0(t)\big) \widehat{K}(t,\eta)=0,  \\
\widehat{K}(t_0,\eta) = \widehat{k}_0(\eta).
\end{array} \right.$$
It leads to the following explicit expression 
\begin{equation} \label{defK}
\widehat{K}(t,\eta)=\widehat{k}_0(\eta) e^{- \frac{1}{2} \int_{t_0}^t |A_0(s)^T R_0(t_0,s)^T \eta |^2ds} e^{- \frac{1}{2} \int_{t_0}^t \text{Tr}(B_0(s))ds},
\end{equation}
for all $(t,\eta) \in [t_0,T_1]\times\mathbb{R}^d$.
By using the Liouville formula 
\begin{equation}\label{iop1}
\forall t_1, t_2 \in [T_0,T_1], \quad \text{det}\big(R_0(t_2,t_1)\big)=\exp\Big(\int_{t_1}^{t_2}\textrm{Tr}\big(B_0(s)\big)ds\Big),
\end{equation}
see e.g.~\cite{cartan} (Proposition II.2.3.1), and the substitution rule with $y=R_0(t_0,t)x$, it follows that
$$\widehat{k}(t,\xi) 
=  \int_{\mathbb{R}^d} K\big(t,R_0(t_0,t)x\big) e^{-i x \cdot \xi} dx 
=  |\text{det}(R_0(t,t_0))| \int_{\mathbb{R}^d} K(t,y) e^{-i (R_0(t,t_0)y) \cdot \xi} dy,$$ 
that is
\begin{multline*}
\widehat{k}(t,\xi) 
= e^{\int_{t_0}^t \text{Tr}(B_0(s))ds} \widehat{K}\big(t,R_0(t,t_0)^T \xi \big) \\
= \widehat{k}_0 \big(R_0(t,t_0)^T \xi \big)  e^{\frac{1}{2} \int_{t_0}^t\text{Tr}(B_0(s))ds } e^{-\frac{1}{2} \int_{t_0}^t |A_0(s)^T R_0(t,s)^T\xi|^2ds}, 
\end{multline*}
since $R_0(t_0,s)^T R_0(t,t_0)^T = \big(R_0(t,t_0)R_0(t_0,s)\big)^T=R_0(t,s)^T$.
We obtain the formula (\ref{k_explicit}).

\medskip

\noindent
\textit{Step 2.} We prove that the $L^2$-function $k$ whose Fourier transform is given by (\ref{k_explicit}), is a weak solution to the Cauchy problem (\ref{CY_Hom}).
We easily notice that $k(t_0)=k_0$ and $k \in C^0([t_0,T_1],L^2(\mathbb{R}^d))$. Then, we use the change of function
\begin{equation}\label{rt5}
\varphi(t,x)=\psi\big(t,R_0(t_0,t)x\big)\big|\text{det}\big(R_0(t_0,t)\big)\big|.
\end{equation}
According to (\ref{rt1}), the function $\psi$ is well-defined.
It follows from the Liouville formula (\ref{iop1}) that 
\begin{multline}\label{rt6}
(\partial_t \varphi)(t,x) + \big\langle B_0(t)x,(\nabla_x\varphi)(t,x)\big\rangle\\
=\big|\text{det}\big(R_0(t_0,t)\big)\big|\big(\partial_t\psi\big(t,R_0(t_0,t)x\big)-\text{Tr}\big(B_0(t)\big)\psi\big(t,R_0(t_0,t)x\big)\big),
\end{multline}
since $\text{det}\big(R_0(t_0,t)\big) \in \rr_+^*$ for all $T_0 \leq t \leq T_1$.
According to (\ref{CVAR:k/K}), (\ref{rt5}) and (\ref{rt6}), it is sufficient to prove that for all $\psi \in \widetilde{E}(t_0,T_1)$, $t_0<t_*<T_1$,
\begin{multline*}
\int_{\mathbb{R}^d} \big(K(t_*,y)\psi(t_*,y) - k_0(y)\psi(t_0,y)\big) dy \\
  =\int_{t_0}^{t_*} \int_{\mathbb{R}^d} K(t,y) \Big( \partial_t \psi 
+ \frac{1}{2} \text{Tr}\big(R_0(t_0,t)A_0(t)A_0(t)^TR_0(t_0,t)^T\nabla_y^2\psi\big) \\ -\frac{1}{2}\text{Tr}\big(B_0(t)\big)\psi\Big)(t,y)dydt.
\end{multline*}
where $\widetilde{E}(t_0,T_1)$ stands for the space of functions $\psi \in C^0([t_0,T_1],L^2(\mathbb{R}^d))$ satisfying
\begin{itemize}
\item[$(i)$] $\psi(\cdot,y) \in C^1((t_0,T_1),\cc)$ for all $y \in \mathbb{R}^d$
\item[$(ii)$] $\psi(t,\cdot) \in C^2(\mathbb{R}^d,\cc)$ for all $t \in (t_0,T_1)$
\item[$(iii)$] The functions $\partial_t \psi$, $\nabla_y^2 \psi$, $\psi$ belong to $L^2((t_0,T_1)\times\mathbb{R}^d)$.
\end{itemize}

\medskip

\noindent
For all $\psi \in \widetilde{E}(t_0,T_1)$ and $t_0<t_*<T_1$, it follows from the Plancherel theorem, (\ref{k_explicit}) and (\ref{CVAR:k/K}) that
\begin{align*}
  & \int\limits_{t_0}^{t_*} \int\limits_{\mathbb{R}^d} K(t,y) 
    \Big( \partial_t \psi + \frac{1}{2} \text{Tr}\big(R_0(t_0,t)A_0(t)A_0(t)^TR_0(t_0,t)^T \nabla_y^2\psi\big)
   - \frac{1}{2}\text{Tr}\big(B_0(t)\big) \psi \Big) dy dt  \\
= & \frac{1}{(2\pi)^d} \int\limits_{t_0}^{t_*} \int\limits_{\mathbb{R}^d} \widehat{K}(t,\eta) 
    \Big(\overline{ \partial_t \widehat{\overline{\psi}} - \frac{1}{2} |A_0(t)^TR_0(t_0,t)^T  \eta|^2 \widehat{\overline{\psi}}-\frac{1}{2}\text{Tr}\big(B_0(t)\big) \widehat{\overline{\psi}} }\Big)(t,\eta)  d\eta dt  \\
= & \frac{1}{(2\pi)^d}  \int\limits_{t_0}^{t_*} \int\limits_{\mathbb{R}^d} \frac{\partial}{\partial t} \Big[ \widehat{K}(t,\eta) \overline{\widehat{\overline{\psi}}(t,\eta)} \Big] d\eta dt 
=  \int\limits_{\mathbb{R}^d} \Big( \widehat{K}(t_*,\eta) \overline{\widehat{\overline{\psi}}(t_*,\eta)} - \widehat{K}(t_0,\eta) \overline{\widehat{\overline{\psi}}(t_0,\eta)} \Big)  \frac{d\eta}{(2\pi)^d} \\
= & \int\limits_{\mathbb{R}^d} \big( K(t_*,y) \psi(t_*,y) - K(t_0,y) \psi(t_0,y) \big) dy=\int\limits_{\mathbb{R}^d} \big( K(t_*,y) \psi(t_*,y) - k_0(y) \psi(t_0,y) \big) dy.
\end{align*}

\bigskip

\noindent
\textit{Step 3: Definition and properties of the evolution system.}
For all $(t,t_0) \in \mathscr{T}$ and $k_0 \in L^2(\mathbb{R}^d)$, we define $\mathscr{U}(t,t_0)k_0$ as the $L^2$-function $k(t)$ whose Fourier transform is given by (\ref{k_explicit}). With this definition, we easily check that $\mathscr{U}(t,t)=\textrm{Id}_{L^2(\rr^d)}$ for all $T_0 \leq t \leq T_1$, and that the mapping $\mathscr{U}$ is strongly continuous from $\mathscr{T}$ to $\mathcal{L}_c(L^2(\mathbb{R}^d))$ thanks to Plancherel theorem. On the other hand, with $k_1=\mathscr{U}(t_1,t_0)k_0$, $k_2=\mathscr{U}(t_2,t_0)k_0$ and $k_3=\mathscr{U}(t_2,t_1)k_1$, it follows from (\ref{k_explicit}) that for all $T_0 \leq t_0 \leq t_1 \leq t_2 \leq T_1$, $k_0 \in L^2(\mathbb{R}^d)$, 
$$\widehat{k}_1(\xi) = \widehat{k}_0\big(R_0(t_1,t_0)^T\xi\big)e^{\frac{1}{2}\int_{t_0}^{t_1}\text{Tr}(B_0(s))ds}e^{-\frac{1}{2}\int_{t_0}^{t_1}|A_0(s)^TR_0(t_1,s)^T \xi|^2ds},$$
$$\widehat{k}_2(\xi) = \widehat{k}_0\big(R_0(t_2,t_0)^T\xi\big)e^{\frac{1}{2}\int_{t_0}^{t_2}\text{Tr}(B_0(s))ds}e^{-\frac{1}{2}\int_{t_0}^{t_2}|A_0(s)^TR_0(t_2,s)^T\xi|^2ds}$$
and
\begin{align*}
\widehat{k}_3(\xi)=& \ \widehat{k}_1\big(R_0(t_2,t_1)^T\xi\big)e^{\frac{1}{2} \int_{t_1}^{t_2}\text{Tr}(B_0(s))ds}e^{-\frac{1}{2} \int_{t_1}^{t_2}|A_0(s)^TR_0(t_2,s)^T \xi|^2ds} \\
= & \  \widehat{k}_0\big(R_0(t_1,t_0)^TR_0(t_2,t_1)^T\xi\big)
 e^{\frac{1}{2} \int_{t_0}^{t_2}\text{Tr}(B_0(s))ds}\\
 & \ \times e^{-\frac{1}{2} \int_{t_0}^{t_1}|A_0(s)^TR_0(t_1,s)^T R_0(t_2,t_1)^T\xi|^2ds} 
 e^{-\frac{1}{2} \int_{t_1}^{t_2}|A_0(s)^TR_0(t_2,s)^T\xi|^2ds}\\
 = & \ \widehat{k}_2(\xi),
\end{align*}
since 
$$R_0(t_1,s)^TR_0(t_2,t_1)^T=\big(R_0(t_2,t_1)R_0(t_1,s)\big)^T=R_0(t_2,s)^T$$
and
$$R_0(t_1,t_0)^TR_0(t_2,t_1)^T=\big(R_0(t_2,t_1)R_0(t_1,t_0)\big)^T=R_0(t_2,t_0)^T.$$
It proves that for all $T_0 \leq t_0 \leq t_1 \leq t_2 \leq T_1$,
$$\mathscr{U}(t_2,t_1)\mathscr{U}(t_1,t_0)=\mathscr{U}(t_2,t_0).$$

\bigskip

\noindent
\textit{Step 4: Uniqueness of the weak solution to the Cauchy problem \emph{(\ref{CY_Hom})}.}
Let $k$ be a  weak solution to the Cauchy problem (\ref{CY_Hom}) associated with the initial datum $k_0=0$. It follows that for all $\varphi \in E(t_0,T_1)$, $t_0<t_*<T_1$,
\begin{equation} \label{WS_uniqueness}
\int\limits_{\mathbb{R}^d} k(t_*,x) \varphi(t_*,x) dx
= \int\limits_{t_0}^{t_*} \int\limits_{\mathbb{R}^d} k(t,x) \big( \partial_t \varphi(t,x) + P_0(t)^* \varphi(t,x) \big) dx dt.
\end{equation}
Let $t_0 \leq t_* \leq T_1$ be fixed. We aim at proving that $k(t_*)= 0$. To that end,
we consider a sequence $(g_p)_{p \geq 1}$ of $C_0^{\infty}(\rr^d)$ functions satisfying 
$$\lim_{p \to +\infty}\|\widehat{g}_p - \widehat{k}(t_*)\|_{L^2(\mathbb{R}^d)}.$$
By Plancherel theorem, we observe that 
\begin{equation}\label{rt7}
\lim_{p \to +\infty}\|g_p - k(t_*)\|_{L^2(\mathbb{R}^d)}=0.
\end{equation}
Following the very same strategy as in the two first steps, we build a weak solution $\varphi_p:(t_0,t_*) \times \mathbb{R}^d \rightarrow \mathbb{C}$ to the Cauchy problem
$$\left\lbrace \begin{array}{ll}
\partial_t \varphi_p(t,x) + P_0(t)^* \varphi_p(t,x)=0, \\
\varphi_p|_{t=t_*}=\overline{g_p}.
\end{array}\right.$$
By deriving a similar formula as in (\ref{k_explicit}), we notice that the function $\varphi_p$ is smooth in the space variable as its Fourier transform in the space variable is compactly supported. This similar formula as in (\ref{k_explicit}) also shows that the function $\varphi_p$ is smooth in the time variable. It follows that the function $\varphi_p$ is a pointwise solution of the equation
$$\forall (t,x) \in (t_0,t_*)\times\mathbb{R}^d, \quad \partial_t \varphi_p(t,x) + P_0(t)^* \varphi_p(t,x)=0.$$
Furthermore, we easily check that $\varphi_p$ is an admissible test function.
Then, we deduce from (\ref{WS_uniqueness}) and (\ref{rt7}) that
$$\forall p \geq 1, \quad \int\limits_{\mathbb{R}^d} k(t_*,x)\overline{g_p(x)}dx=(k(t_*),g_p)_{L^2(\rr^d)}=0,$$
implying that $k(t_*)=0$, when passing to the limit $p \to +\infty$. It ends the proof of the uniqueness of the weak solution to the Cauchy problem (\ref{CY_Hom}).
\end{proof}

\medskip

Regarding the inhomogeneous equation
\begin{equation} \label{CY_Nonhom}
\left\lbrace \begin{array}{ll}
\partial_t h(t,x) - P_0(t) h(t,x)=\un_{\omega(t)}(x)u(t,x), \quad & (t,x) \in (t_0,T_1) \times \mathbb{R}^d, \\
h|_{t=t_0} = h_0 \in L^2(\mathbb{R}^d),
\end{array} \right.
\end{equation}
with $(\omega(t))_{t \in [T_0,T_1]}$ being a moving control support in $\rr^d$,
we use the notion of mild solutions defined in~\cite{pazy} (Chapter~5, Section~5.5, Definition~5.1, p.~146):

\medskip

\begin{Def}
Let $T_0 \leq t_0 < T_1$, $A_0, B_0 \in C^0([T_0,T_1],M_d(\mathbb{R}))$, $k_0 \in L^2(\mathbb{R}^d)$, $u \in L^1((t_0,T_1),L^2(\mathbb{R}^d))$ and $(\omega(t))_{t \in [T_0,T_1]}$ a moving control support in $\rr^d$.
The mild solution to the Cauchy problem \emph{(\ref{CY_Nonhom})} is the function $h \in C^0([t_0,T_1],L^2(\mathbb{R}^d))$  given by
$$h(t)=\mathscr{U}(t,t_0)h_0 + \int_{t_0}^t \mathscr{U}(t,s) [\un_{\omega(s)}u(s)] ds,$$
with equality in $L^2(\mathbb{R}^d)$ for all $t \in [t_0,T_1]$, where $\mathscr{U}$ stands for the evolution system given by 
Proposition~\emph{\ref{Prop:WP_Hom}}.
\end{Def}

\subsection{Hilbert uniqueness method} \label{Appendix:HUM}

This section is devoted to the proof of following characterization of null-controllability. This result extends the one established in~\cite{KB_KPS_JMAA} (Proposition~2.8) to the case of moving control support:

\medskip

\begin{Prop}\label{Prop:HUM}
Let $T_0<0<T_1$, $A_0, B_0 \in C^0([T_0,T_1],M_d(\mathbb{R}))$, 
$P_0(t)$ be the non-autonomous Ornstein-Uhlenbeck operator defined in \emph{(\ref{def:P_appendix})} and $P_0(t)^*$ its $L^2(\rr^d)$-adjoint given in \emph{(\ref{def:P*_appendix})}.
Let $T \in (0,T_1)$ and $(\omega(t))_{t \in [0,T]}$ be a moving control support on $[0,T]$ in $\mathbb{R}^d$.
The null-controllability on $[0,T]$ from the moving control support $(\omega(t))_{t \in [0,T]}$ of the system
\begin{equation} \label{syst_controle}
\left\lbrace \begin{array}{ll}
\partial_t f(t,x) - P_0(t) f(t,x)=\un_{\omega(t)}(x)u(t,x), \quad & (t,x) \in (0,T_1) \times \mathbb{R}^d, \\
f|_{t=0} = f_0 \in L^2(\mathbb{R}^d),
\end{array} \right.
\end{equation}
is equivalent to the existence of an observability constant $C>0$ such that,
for all $g_0 \in L^2(\mathbb{R}^d)$, the mild solution to the Cauchy problem
\begin{equation} \label{syst_adj}
\left\lbrace \begin{array}{ll}
\partial_t g(t,x) - P_0(T-t)^* g(t,x)=0, \quad & (t,x) \in (0,T-T_0) \times \mathbb{R}^d, \\
g|_{t=0} = g_0 \in L^2(\mathbb{R}^d),
\end{array} \right.
\end{equation}
satisfies
$$\int\limits_{\mathbb{R}^d} |g(T,x)|^2 dx \leq C \int\limits_0^T \int\limits_{\omega(T-t)} |g(t,x)|^2 dx dt.$$
\end{Prop}

\medskip

\noindent
Instrumental in the proof of Proposition~\ref{Prop:HUM} is the following result:

\medskip

\begin{Lem} \label{Prop:U*}
If $\mathscr{U}$ stands for the evolution system given by Proposition~\emph{\ref{Prop:WP_Hom}}, then the $L^2$-adjoint $\mathscr{U}(t,t_0)^*$ of the evolution mapping $\mathscr{U}(t,t_0)$ is  equal to the evolution mapping $\tilde{\mathscr{U}}(t-t_0,0)$ associated to the family of operators $s \in [0,t-t_0] \mapsto P_0(t-s)^*$.
\end{Lem}

\medskip

\begin{proof} Let $T_0 \leq t_0 < t \leq T_1$ and $g_0 \in L^2(\mathbb{R}^d)$. Setting $g(t-t_0)=\tilde{\mathscr{U}}(t-t_0,0)g_0$, we deduce from 
Proposition~\ref{Prop:WP_Hom} with the suitable substitutions that 
\begin{multline}\label{rt11}
\widehat{g}(t-t_0,\xi) =\\
  \widehat{g}_0 \big( \mathcal{R}(t-t_0,0)^T \xi \big)  
e^{-\frac{1}{2} \int_{0}^{t-t_0}\text{Tr}(B_0(t-s))ds}e^{-\frac{1}{2} \int_{0}^{t-t_0}|A_0(t-s)^T \mathcal{R}(t-t_0,s)^T \xi|^2ds},
\end{multline}
where $\mathcal{R}$ stands for the resolvent associated to the system $\dot{X}(s)=-B_0(t-s)X(s)$, that is,
$$\left\lbrace \begin{array}{l}
\frac{\partial \mathcal{R}}{\partial s}(s,s_0)=-B_0(t-s)\mathcal{R}(s,s_0), \\
\mathcal{R}(s_0,s_0)=I_d.
\end{array}\right.$$
According to (\ref{rt9}), we notice that
\begin{equation}\label{rt10}
\mathcal{R}(s,s_0)=R_0(t-s,t-s_0).
\end{equation}
It follows from (\ref{rt11}) and (\ref{rt10}) that 
\begin{equation}\label{rt13}
\widehat{g}(t-t_0,\xi)  = \widehat{g_0} \big(R_0(t_0,t)^T \xi \big) 
e^{-\frac{1}{2} \int_{t_0}^{t}\text{Tr}(B_0(s))ds}e^{-\frac{1}{2} \int_{t_0}^{t}|A_0(s)^TR_0(t_0,s)^T \xi|^2ds}.
\end{equation}
We deduce from Plancherel theorem, (\ref{k_explicit}), the substitution rule with $\eta= R_0(t,t_0)^T \xi$,
the Liouville formula (\ref{iop1}) and (\ref{rt13}) that for all $k_0$, $g_0 \in L^2(\mathbb{R}^d)$,
\begin{align*}
  & \big(\mathscr{U}(t,t_0)k_0,g_0\big)_{L^2(\mathbb{R}^d)}
= \frac{1}{(2\pi)^d}\int_{\mathbb{R}^d} \widehat{k}(t,\xi) \overline{\widehat{g}_0(\xi)} d\xi 
\\
= & \int_{\mathbb{R}^d} \widehat{k}_0 \big(R_0(t,t_0)^T \xi \big)  
  e^{\frac{1}{2} \int_{t_0}^t\text{Tr}(B_0(s))ds}
  e^{-\frac{1}{2} \int_{t_0}^t|A_0(s)^TR_0(t,s)^T\xi|^2ds} \overline{\widehat{g}_0(\xi)}  \frac{d\xi}{(2\pi)^d}
\\
= & \int_{\mathbb{R}^d} \widehat{k}_0(\eta)   
e^{-\frac{1}{2} \int_{t_0}^t\text{Tr}(B_0(s))ds}
  e^{-\frac{1}{2} \int_{t_0}^t|A_0(s)^TR_0(t,s)^TR_0(t_0,t)^T\eta|^2ds} 
 \overline{\widehat{g}_0\big( R_0(t_0,t)^T \eta\big)}  \frac{d\eta}{(2\pi)^d}
\\
= & \int_{\mathbb{R}^d} \widehat{k}_0(\eta)\overline{\widehat{g}(t-t_0,\eta)}  \frac{d\eta}{(2\pi)^d}
=\int_{\mathbb{R}^d} k_0(x)\overline{g(t-t_0,x)}dx
=\big(k_0,\tilde{U}(t-t_0,0)g_0\big)_{L^2(\rr^d)},
\end{align*}
since
$$R_0(t,s)^TR_0(t_0,t)^T=R_0(t_0,s)^T, \quad \big|\text{det}\big(R_0(t_0,t)\big)\big|=\exp\Big(-\int_{t_0}^{t}\textrm{Tr}\big(B_0(s)\big)ds\Big).$$
This ends the proof of Lemma~\ref{Prop:U*}.
\end{proof}

\medskip

\begin{proof} We can now resume the proof of Proposition~\ref{Prop:HUM}. We consider the following linear mappings
$$
\begin{array}{crcl}
C_2: & L^2(\mathbb{R}^d) & \rightarrow & L^2(\mathbb{R}^d), \\
     &     f_0           & \mapsto     & \mathscr{U}(T,0)f_0
\end{array}$$
and
$$\begin{array}{crcl}
C_3: & L^2((0,T)\times\mathbb{R}^d) & \rightarrow & L^2(\mathbb{R}^d), \\
     &     u                        & \mapsto     & \int_0^T \mathscr{U}(T,s)[ \un_{\omega(s)} u(s)  ]  ds.
\end{array}$$
For any $f_0 \in L^2(\mathbb{R}^d)$, the function $C_2(f_0)=k(T)$ is the weak solution at time $T$ to the Cauchy problem (\ref{CY_Hom}) associated to the initial datum $k_0=f_0$, with $t_0=0$. On the other hand, for any $u \in L^2((0,T)\times\mathbb{R}^d)$, 
the function $C_3(u)=h(T)$ is the mild solution to the Cauchy problem (\ref{CY_Nonhom}) at time $T$ associated to the initial datum $h_0=0$, with $t_0=0$.

The null-controllability of the non-autonomous Ornstein-Uhlenbeck equation (\ref{syst_controle}) on $[0,T]$ from the moving control support $(\omega(t))_{t \in [0,T]}$ is equivalent to the set inclusion 
$$C_2(L^2(\mathbb{R}^d)) \subset C_3(L^2((0,T)\times\mathbb{R}^d)),$$ 
since  the mild solution at time $T$ is given by $f(T)=C_2(f_0)+C_3(u)$.
According to~\cite{coron_book} (Lemma~2.48), this set inclusion is also equivalent to the existence of a positive constant $M>0$ such that for all $g_0 \in L^2(\mathbb{R}^d)$,
\begin{equation}\label{rt20}
\| C_2^*g_0 \|_{L^2(\mathbb{R}^d)} \leq M \|C_3^* g_0 \|_{L^2((0,T)\times\mathbb{R}^d)},
\end{equation}
where $C_2^*$ and $C_3^*$ denote the adjoint operators.
We deduce from Lemma~\ref{Prop:U*} that 
\begin{equation}\label{rt21}
C_2^* g_0 = \mathscr{U}(T,0)^* g_0= \tilde{\mathscr{U}}(T,0)g_0 = g(T),
\end{equation}
where $g$ is the weak solution of (\ref{syst_adj}). On the other hand, it follows from Lemma~\ref{Prop:U*} that for all $u \in L^2((0,T)\times\mathbb{R}^d)$,
\begin{multline}\label{rt22}
 (u,C_3^*g_0)_{L^2((0,T)\times \mathbb{R}^d)}=(C_3 u,g_0)_{L^2(\mathbb{R}^d)}\\
=  \int_0^T \big(\mathscr{U}(T,s)[\un_{\omega(s)}u(s)],g_0\big)_{L^2(\mathbb{R}^d)} ds 
= \int_0^T \big(u(s),\un_{\omega(s)} \mathscr{U}(T,s)^*g_0\big)_{L^2(\mathbb{R}^d)} ds\\ 
=  \int_0^T \big(u(s),\un_{\omega(s)} \tilde{\mathscr{U}}(T-s,0)g_0\big)_{L^2(\mathbb{R}^d)} ds 
= \int_0^T \big(u(s),\un_{\omega(s)} g(T-s)\big)_{L^2(\mathbb{R}^d)} ds,
\end{multline}
where $g$ is the weak solution of (\ref{syst_adj}). It follows from (\ref{rt21}) and (\ref{rt22}) that the estimate (\ref{rt20}) reads as 
$$\int\limits_{\mathbb{R}^d} |g(T,x)|^2 dx \leq M^2 \int\limits_0^T \int\limits_{\omega(T-t)} |g(t,x)|^2 dx dt.$$
This ends the proof of Proposition~\ref{Prop:HUM}.
\end{proof}

\subsection{Miscellaneous facts about thick subsets in $\rr^d$}\label{Sec:rds}

In the works~\cite{Panejah-61,Panejah-62}, Panejah addressed the problem of characterizing the subsets $S \subset \rr^d$ for which the semi-norm $\|\cdot\|_{L^2(S)} $ defines a norm on specific vector subspaces of $L^2(\mathbb{R}^d)$. For the vector subspace given by $L^2(\rr^d)$-functions whose Fourier transforms are supported in a fixed compact set of $\rr^d$,
the thickness property of the subset $S$ turns out to be a necessary and sufficient condition for the semi-norm $\|\cdot\|_{L^2(S)} $ to be a norm on this vector subspace, see Definition~\ref{Def:thick}. This result was established by Panejah in the one-dimensional setting. In the multidimensional case, Panejah only proved that the thickness property is a necessary condition for the semi-norm $\|\cdot\|_{L^2(S)} $ to be a norm. The full equivalence in the multidimensional setting was then established independently by Logvinenko and Sereda~\cite{Logvinenko_Sereda}, and by Katsnelson~\cite{Katsnelson-73}.

\medskip

\begin{thm}[Logvinenko-Sereda~\cite{Logvinenko_Sereda}, Katsnelson~\cite{Katsnelson-73}]
Let $S$ be a measurable subset of $\mathbb{R}^d$.
The two following assertions are equivalent:
\begin{enumerate}
\item[$(i)$] $S$ is a thick subset in $\mathbb{R}^d$
\item[$(ii)$] For all bounded subsets $\Sigma \subset \mathbb{R}^d$,
there exists a positive constant $C(\Sigma,S)>0$ such that
$$\forall f  \in L^2(\mathbb{R}^d),\ \emph{\text{supp }}\widehat{f} \subset \Sigma, \quad \|f\|_{L^2(\mathbb{R}^d)} \leq C(\Sigma,S) \|f\|_{L^2(S)}.$$
\end{enumerate}
\end{thm}

\medskip

We refer the reader to the monographs \cite{Havin_book,MuscaluSchlag} for detailed introductions on this topic. Let us only notice for now  that the estimate in assertion $(ii)$ is actually a spectral estimate in the common terminology used in control theory. In order to use this kind of the spectral estimates in control theory, it is then essential to understand how the positive constants $C(\Sigma,S)>0$ do depend on the bounded subsets $\Sigma$. This problem was addressed by Kovrijkine who established in~\cite{Kovrijkine} (Theorem~3) the following quantitative version of the estimate $(ii)$:

\medskip

\begin{thm}[Kovrijkine~\cite{Kovrijkine}] \label{thm_PLSK}
There exists a universal constant $\mathcal{C}>1$ such that for all $d \geq 1$, $0<\delta \leq 1$, $\alpha=(\alpha_1,...,\alpha_d) \in (0,+\infty)^d$, $\beta=(\beta_1,...,\beta_d) \in (0,+\infty)^d$,
for all $(\delta,\alpha)$-thick set $S$ in $\mathbb{R}^d$
and for all $f \in L^2(\mathbb{R}^d)$ whose Fourier transform is supported in a parallelepiped with sides parallel to the coordinate axes with length respectively $\beta_1$,..., $\beta_d$
then
\begin{equation}\label{ty3}
\|f\|_{L^2(\mathbb{R}^d)} \leq \left( \frac{\mathcal{C}^d}{\delta} \right)^{\mathcal{C}(\alpha \cdot \beta+d)} \|f\|_{L^2(S)},
\end{equation}
where $\alpha \cdot \beta=\sum_{j=1}^d\alpha_j \beta_j$ denotes the Euclidean dot product.
\end{thm}

\medskip

\subsection{Hermite functions}\label{6.sec.harmo}
This section is devoted to set some notations and recall basic facts about Hermite functions.
 The standard Hermite functions $(\phi_{n})_{n\geq 0}$ are defined for $x \in \rr$,
 \begin{equation}\label{Ar1}
 \phi_{n}(x)=\frac{(-1)^n}{\sqrt{2^n n!\sqrt{\pi}}} e^{\frac{x^2}{2}}\frac{d^n}{dx^n}(e^{-x^2})
 =\frac{1}{\sqrt{2^n n!\sqrt{\pi}}} \Bigl(x-\frac{d}{dx}\Bigr)^n(e^{-\frac{x^2}{2}})=\frac{ a_{+}^n \phi_{0}}{\sqrt{n!}},
 \end{equation}
where $a_{+}$ is the creation operator
$$a_{+}=\frac{1}{\sqrt{2}}\Big(x-\frac{d}{dx}\Big).$$
The family $(\phi_{n})_{n\geq 0}$ is an orthonormal basis of $L^2(\R)$
composed by the eigenfunctions of the harmonic oscillator
$$\mathcal{H}_1=-\Delta_x+x^2=\sum_{n\ge 0}(2n+1)\mathbb P_{n},\quad 1=\sum_{n \ge 0}\mathbb P_{n},$$
where $\mathbb P_{n}$ stands for the orthogonal projection
$$\mathbb P_{n}f=(f,\phi_n)_{L^2(\rr)}\phi_n.$$
It satisfies the identities
\begin{equation}\label{ge1}
\forall n \geq 0, \quad a_+\phi_n=\sqrt{n+1}\phi_{n+1}, \quad a_-\phi_n=\sqrt{n}\phi_{n-1} \textrm{ (}=0\textrm{ si }n=0\textrm{)} ,
\end{equation}
where
\begin{equation}\label{ge2}
a_{\pm}=\frac{1}{\sqrt{2}}\Big(x\mp\frac{d}{dx}\Big).
\end{equation}
Instrumental in the core of this work are the following estimates on the Hermite functions given in the next lemma. They are an adaptation in a simpler setting of the analysis led in the work~\cite{langen} (Lemma~3.2). The same estimates were also established in~\cite{karelJFA} (Lemma A.1) while using a different normalization for the harmonic oscillator and the Hermite functions. For the sake of completeness, we adapt the proof given in~\cite{karelJFA} with the normalization used in the present work:

\medskip

\begin{Lem}\label{ge3}
We have
\begin{equation}\label{ge4}
\forall n, k,l \geq 0, \quad \|x^k\partial_x^l\phi_n\|_{L^2(\rr)} \leq 2^{\frac{k+l}{2}}\sqrt{\frac{(k+l+n)!}{n!}},
\end{equation}
\begin{multline}\label{ge5}
\forall r \geq \frac{1}{2}, \forall \eps >0, \forall n, k,l \geq 0,  \\ \|x^k\partial_x^l\phi_n\|_{L^2(\rr)} \leq \sqrt{2}\big((1-\delta_{n,0})\exp(\eps r n^{\frac{1}{2r}})+\delta_{n,0}\big)
\Big(\frac{2^{r+1}e^r}{\inf(\eps^r,1)}\Big)^{k+l}(k!)^r(l!)^r,
\end{multline}
where $\delta_{n,0}$ stands for the Kronecker symbol, i.e., $\delta_{n,0}=1$ if $n=0$, $\delta_{n,0}=0$ if $n \neq 0$.  
\end{Lem}

\medskip

\begin{proof}
The estimate (\ref{ge4}) is trivial if $k=l=0$, since the family $(\phi_{n})_{n \geq 0}$ is an orthonormal basis of $L^2(\R)$. We notice from (\ref{ge1}) and (\ref{ge2}) that
\begin{equation}\label{ge6}
x\phi_n=\frac{1}{\sqrt{2}}(a_++a_-)\phi_n=\sqrt{\frac{n+1}{2}}\phi_{n+1}+\sqrt{\frac{n}{2}}\phi_{n-1},
\end{equation}
\begin{equation}\label{ge7}
\partial_x\phi_n=\frac{1}{\sqrt{2}}(a_--a_+)\phi_n=\sqrt{\frac{n}{2}}\phi_{n-1}-\sqrt{\frac{n+1}{2}}\phi_{n+1}.
\end{equation}
This implies that
$$\|x\phi_n\|_{L^2(\rr)}=\sqrt{\frac{2n+1}{2}}, \quad \|\partial_x\phi_n\|_{L^2(\rr)}=\sqrt{\frac{2n+1}{2}},$$
since $(\phi_{n})_{n \geq 0}$ is an orthonormal basis of $L^2(\R)$.
It follows that the estimate (\ref{ge4}) holds as well when $(k,l)=(1,0)$ or $(k,l)=(0,1)$. We complete the proof of the estimate (\ref{ge4}) by induction. We assume that the estimate holds for any $k,l \geq 0$, $k+l \leq m$, with $m \geq 1$. Let $k,l \geq 0$ such that $k+l=m$. It follows from (\ref{ge6}) and (\ref{ge7}) that
$$x^{k+1}\partial_x^l\phi_n=\sqrt{\frac{n+1}{2}}x^{k}\partial_x^l\phi_{n+1}+\sqrt{\frac{n}{2}}x^{k}\partial_x^l\phi_{n-1}-lx^{k}\partial_x^{l-1}\phi_{n},$$
$$x^k\partial_x^{l+1}\phi_n=\sqrt{\frac{n}{2}}x^{k}\partial_x^l\phi_{n-1}-\sqrt{\frac{n+1}{2}}x^{k}\partial_x^l\phi_{n+1}.$$
We deduce from the induction hypothesis that
\begin{align*}
& \ \|x^{k+1}\partial_x^l\phi_n\|_{L^2(\rr)} \\
\leq & \ \sqrt{\frac{n+1}{2}}\|x^{k}\partial_x^l\phi_{n+1}\|_{L^2(\rr)}+\sqrt{\frac{n}{2}}\|x^{k}\partial_x^l\phi_{n-1}\|_{L^2(\rr)}
 +l\|x^{k}\partial_x^{l-1}\phi_{n}\|_{L^2(\rr)}  \\
\leq & \ 2^{\frac{k+l-1}{2}}\sqrt{\frac{(k+l+n+1)!}{n!}}\Big(1+\frac{n+l}{\sqrt{(k+l+n+1)(k+l+n)}}\Big)\\
 \leq &\  2^{\frac{k+l+1}{2}}\sqrt{\frac{(k+l+n+1)!}{n!}},
\end{align*}
\begin{align*}
\|x^k\partial_x^{l+1}\phi_n\|_{L^2(\rr)} \leq & \ \sqrt{\frac{n}{2}}\|x^{k}\partial_x^l\phi_{n-1}\|_{L^2(\rr)}+\sqrt{\frac{n+1}{2}}\|x^{k}\partial_x^l\phi_{n+1}\|_{L^2(\rr)}\\ 
\leq &\ 2^{\frac{k+l-1}{2}}\sqrt{\frac{(k+l+n+1)!}{n!}}\Big(1+\frac{n}{\sqrt{(k+l+n+1)(k+l+n)}}\Big) \\
\leq & \ 2^{\frac{k+l+1}{2}}\sqrt{\frac{(k+l+n+1)!}{n!}}.
\end{align*}
This ends the proof of the estimate (\ref{ge4}). 
We now prove the estimates (\ref{ge5}).
When $n=0$, we deduce from (\ref{ge4}) that
$$\forall k,l \geq 0, \quad \|x^k\partial_x^l\phi_0\|_{L^2(\rr)} \leq 2^{\frac{k+l}{2}}\sqrt{{(k+l)!}}\le 2^{k+l}\sqrt{{k!}}\sqrt{{l!}},$$
since 
\begin{equation}\label{gh1}
\binom{k+l}{k}=\frac{(k+l)!}{k!l!}\leq \sum_{j=0}^{k+l}\binom{k+l}{j} = 2^{k+l}.
\end{equation}
It follows that 
$$\forall r \geq \frac{1}{2}, \forall \eps >0, \forall k,l \geq 0, \quad 
\|x^k\partial_x^l\phi_0\|_{L^2(\rr)} \leq \sqrt{2}\Big(\frac{2^{r+1}e^r}{\inf(\eps^r,1)}\Big)^{k+l}(k!)^r(l!)^r,$$
since
$$\forall r \geq \frac{1}{2}, \forall \eps >0, \quad 2^{k+l}\sqrt{{k!}}\sqrt{{l!}} \leq 2^{k+l}({k!})^r({l!})^r \leq \sqrt{2}\Big(\frac{2^{r+1}e^r
}{\inf(\eps^r,1)}\Big)^{k+l}({k!})^r({l!})^r.$$
The estimates (\ref{ge5}) therefore hold when $n=0$.
When $k=l=0$ and $n \geq 1$, the estimates (\ref{ge5}) also hold since $\|\phi_n\|_{L^2(\rr)} =1$.
From now, we may therefore assume that $k+l \geq 1$ and $n\ge 1$.
We notice that for all $n\ge 1$,
\begin{multline*}
n!=\Gamma(n+1)=\int_{0}^{+\infty} e^{-t} t^{n}dt=\Big(\frac n e\Big)^{n}\int_{0}^{+\infty} ne^{-(s-1)n} s^{n}ds \\ \ge \Big(\frac n e\Big)^{n}\int_{1}^{2} ne^{-(s-1)n} ds 
=\Big(\frac n e\Big)^{n}(1-e^{-n})\ge\frac12\Big(\frac n e\Big)^{n},
\end{multline*}
so that 
$$\forall n \geq 1, \quad n^{\frac n2}\le \sqrt{2} \sqrt{n!} e^{\frac{n}{2}}.$$
It follows that 
\begin{equation}\label{vbn1}
\forall r\ge \frac{1}{2}, \forall n\ge 1,\quad n^{\frac n2}\le \sqrt{2} \sqrt{n! e^{n}}\le  \sqrt{2}\bigl({n! e^{n}}\bigr)^{r}.
\end{equation}
We distinguish two cases. When $1\le k+l \leq n$, we deduce from (\ref{ge4}) that for all $r \geq 1/2$ and $\eps >0$,
\begin{multline}\label{ge9}
\|x^k\partial_x^l\phi_n\|_{L^2(\rr)} \leq 2^{\frac{k+l}{2}}\sqrt{\frac{(k+l+n)!}{n!}} \leq 2^{\frac{k+l}{2}}(k+l+n)^{\frac{k+l}{2}}\leq 2^{\frac{k+l}{2}}(2n)^{\frac{k+l}{2}}\\
\leq \Big(\frac{2}{\eps^r}\Big)^{k+l}((k+l)!)^r\Big(\frac{(\eps n^{\frac{1}{2r}})^{k+l}}{(k+l)!}\Big)^r
\leq \Big(\frac{2}{\eps^r}\Big)^{k+l}\exp(\eps r n^{\frac{1}{2r}})((k+l)!)^r.
\end{multline}
When $k+l > n\ge 1$, we deduce from (\ref{ge4}) and (\ref{vbn1}) that for all $r \geq 1/2$,
\begin{multline}\label{ge10}
\|x^k\partial_x^l\phi_n\|_{L^2(\rr)} \leq 2^{\frac{k+l}{2}}\sqrt{\frac{(k+l+n)!}{n!}} \leq 2^{\frac{k+l}{2}}(k+l+n)^{\frac{k+l}{2}}\\ \leq 2^{\frac{k+l}{2}}(2k+2l)^{\frac{k+l}{2}}
\leq  \sqrt 2 (2e^{r})^{k+l}((k+l)!)^r.
\end{multline}
It follows from (\ref{ge9}) and (\ref{ge10}) that for all $r \geq 1/2$, $\eps >0$, $n \geq 1$, $k+l \geq 1$, 
$$\|x^k\partial_x^l\phi_n\|_{L^2(\rr)} \leq \sqrt 2 \Big(\frac{2e^r}{\inf(\eps^r,1)}\Big)^{k+l}\exp(\eps r n^{\frac{1}{2r}})((k+l)!)^r.$$
By using anew (\ref{gh1}), we finally obtain that for all $r \geq 1/2$, $\eps >0$, $n \ge 1$, $k+l \geq 1$,
$$\|x^k\partial_x^l\phi_n\|_{L^2(\rr)} \leq \sqrt 2 \Big(\frac{2^{r+1}e^r}{\inf(\eps^r,1)}\Big)^{k+l}\exp(\eps r n^{\frac{1}{2r}})(k!)^r(l!)^r.$$
The estimates (\ref{ge5}) therefore hold when $k+l \geq 1$ and $n \geq 1$. The proof of Lemma~\ref{ge3} is complete.
\end{proof}

The $d$-dimensional Hermite functions $(\Phi_{\alpha})_{\alpha \in \nn^d}$,
$$\Phi_{\alpha}(x)=\prod_{j=1}^d\phi_{\alpha_j}(x_j), \quad x=(x_1,...,x_d) \in \rr^d, \ \alpha=(\alpha_1,...,\alpha_d) \in \nn^d,$$ 
is an orthonormal basis of $L^2(\R^d)$ composed by the eigenfunctions of the $d$-dimensional harmonic oscillator
$$\mathcal{H}_d=-\Delta_x+|x|^2=\sum_{\alpha \in \nn^d}(2|\alpha|+d)\mathbb P_{\alpha},\quad 1=\sum_{\alpha \in \nn^d}\mathbb P_{\alpha},$$
with $|\alpha|=\alpha_1+...+\alpha_d$, when $\alpha=(\alpha_1,...,\alpha_d) \in \nn^d$, where $\mathbb P_{\alpha}$ stands for the orthogonal projection
$$\mathbb P_{\alpha}f=(f,\Phi_{\alpha})_{L^2(\rr^d)}\Phi_{\alpha}.$$
It satisfies the identities
\begin{equation}\label{ge1b}
a_{+,j}\Phi_{\alpha}=\sqrt{\alpha_j+1}\Phi_{\alpha+e_j}, \quad a_{-,j}\Phi_{\alpha}=\sqrt{\alpha_j}\Phi_{\alpha-e_j} \textrm{ (}=0\textrm{ si }\alpha_j=0\textrm{)},
\end{equation}
with $(e_1,....,e_d)$ the canonical basis of $\rr^d$, where
\begin{equation}\label{ge2b}
a_{\pm,j}=\frac{1}{\sqrt{2}}\Big(x_j\mp\frac{d}{dx_j}\Big).
\end{equation}
By tensorization, the estimates proved in Lemma~\ref{ge3} can be readily extended to the multidimensional setting as follows:
 
\bigskip

\begin{Lem}\label{ge3bis}
We have
\begin{equation}\label{ge4bis}
\forall \alpha, \beta, \gamma \in \nn^d, \quad \|x^{\alpha}\partial_x^{\beta}\Phi_{\gamma}\|_{L^2(\rr^d)} \leq 2^{\frac{|\alpha|+|\beta|}{2}}\sqrt{\frac{(\alpha+\beta+\gamma)!}{\gamma!}},
\end{equation}
with $\alpha!=\prod_{j=1}^d\alpha_j!$, when $\alpha=(\alpha_1,...,\alpha_d) \in \nn^d$ and $|\alpha|=\alpha_1+...+\alpha_d$,
\begin{multline}\label{ge5bis}
\forall r \geq \frac{1}{2}, \forall \eps >0, \forall \alpha, \beta, \gamma \in \nn^d,  \\ \|x^{\alpha}\partial_x^{\beta}\Phi_{\gamma}\|_{L^2(\rr^d)} \leq 2^{\frac{d}{2}}\big((1-\delta_{\gamma,0})\exp(\eps r d |\gamma|^{\frac{1}{2r}})+\delta_{\gamma,0}\big)
\Big(\frac{2^{r+1}e^r}{\inf(\eps^r,1)}\Big)^{|\alpha|+|\beta|}(\alpha!)^r(\beta!)^r,
\end{multline}
where $\delta_{\gamma,0}=\prod_{j=1}^d\delta_{\gamma_j,0}$ stands for the Kronecker symbol, i.e., $\delta_{\gamma,0}=1$ if $\gamma=(0,...,0)$, or $\delta_{\gamma,0}=0$ if $\gamma \neq (0,...,0)$.  
\end{Lem}

\medskip

\subsection{Gelfand-Shilov regularity}\label{GSreg}

We refer the reader to the works~\cite{gelfand,rodino1,rodino,toft} and the references herein for extensive expositions of the Gelfand-Shilov regularity theory.
The Gelfand-Shilov spaces $S_{\nu}^{\mu}(\rr^d)$, with $\mu,\nu>0$, $\mu+\nu\geq 1$, are defined as the spaces of smooth functions $f \in C^{\infty}(\rr^d)$ satisfying the estimates
$$\exists A,C>0, \quad |\partial_x^{\alpha}f(x)| \leq C A^{|\alpha|}(\alpha !)^{\mu}e^{-\frac{1}{A}|x|^{1/\nu}}, \quad x \in \rr^d, \ \alpha \in \mathbb{N}^d,$$
or, equivalently
$$\exists A,C>0, \quad \sup_{x \in \rr^d}|x^{\beta}\partial_x^{\alpha}f(x)| \leq C A^{|\alpha|+|\beta|}(\alpha !)^{\mu}(\beta !)^{\nu}, \quad \alpha, \beta \in \mathbb{N}^d.$$
These Gelfand-Shilov spaces  $S_{\nu}^{\mu}(\rr^d)$ may also be characterized as the spaces of Schwartz functions $f \in \mathscr{S}(\rr^d)$ satisfying the estimates
$$\exists C>0, \exists \eps>0, \quad |f(x)| \leq C e^{-\eps|x|^{1/\nu}}, \ x \in \rr^d, \qquad |\widehat{f}(\xi)| \leq C e^{-\eps|\xi|^{1/\mu}}, \ \xi \in \rr^d.$$
In particular, we notice that Hermite functions belong to the symmetric Gelfand-Shilov space  $S_{1/2}^{1/2}(\rr^d)$. More generally, the symmetric Gelfand-Shilov spaces $S_{\mu}^{\mu}(\rr^d)$, with $\mu \geq 1/2$, can be nicely characterized through the decomposition into the Hermite basis $(\Phi_{\alpha})_{\alpha \in \mathbb{N}^d}$, see e.g. \cite{toft} (Proposition~1.2),
\begin{multline*}
f \in S_{\mu}^{\mu}(\rr^d) \Leftrightarrow f \in L^2(\rr^d),  \exists t_0>0, \ \big\|\big((f,\Phi_{\alpha})_{L^2}\exp({t_0|\alpha|^{\frac{1}{2\mu}})}\big)_{\alpha \in \mathbb{N}^d}\big\|_{l^2(\mathbb{N}^d)}<+\infty\\
\Leftrightarrow f \in L^2(\rr^d),  \exists t_0>0, \quad  \|e^{t_0\mathcal{H}_d^{\frac{1}{2\mu}}}f\|_{L^2(\rr^d)}<+\infty,
\end{multline*}
where $\mathcal{H}_d=-\Delta_x+|x|^2$ stands for the harmonic oscillator.

\end{document}